\documentclass[leqno,12pt]{article}
\usepackage{amsmath,amssymb,rotating,a4wide,color}
\usepackage{wasysym}
\usepackage{hyperref}
\usepackage{enumerate}
\usepackage{mathtools}
\usepackage{float}
\usepackage{colortbl}
\usepackage{makecell}
\usepackage[T1]{fontenc}
\usepackage[font=footnotesize,labelfont=bf,margin=2cm]{caption}
\usepackage[isbn=false]{biblatex}
\usepackage{array}
\newcolumntype{M}{>{\centering\arraybackslash}m{\dimexpr.25\linewidth-2\tabcolsep}}

\usepackage{tikz}
\usetikzlibrary{circuits.logic.US,circuits.logic.IEC,fit}
\usetikzlibrary{cd}
\usetikzlibrary{backgrounds}
\usetikzlibrary{math}
\usepackage{subcaption}
\usepackage{lscape}

\usepackage[all,cmtip]{xy}

\usepackage{verbatim}

\newdir^{ (}{{}*!/-3pt/\dir^{(}}
\newdir_{ (}{{}*!/-3pt/\dir_{(}}

\entrymodifiers={+!!<0pt,\fontdimen22\textfont2>}

\def\geqone{g\kern1pt{=}\kern1pt1}
\def\geqtwo{g\kern1pt{=}\kern1pt2}
\def\geqnine{g\kern1pt{=}\kern1pt9}

\def\FigThreeSomeA{\begin{tikzpicture}[scale=0.7]
		\draw [thick,black] plot  coordinates {(4/8+8/8,-1/8-8/32) (-20/8,5/8) };
		\begin{scope}[shift={(4/8,-1/8)}]
			\draw [thick,black] plot  coordinates {(-1/32,-1/8) (1/32,1/8) };
		\end{scope}
			\begin{scope}[shift={(8/8,-8/32)}]
			\draw [thick,black] plot  coordinates {(-1/32,-1/8) (1/32,1/8) };
		\end{scope}
	\begin{scope}[shift={(0,0)}]
		\draw [thick,black] plot  coordinates {(2/24,8/24) (-5/16,-20/16) };
	\fill [red] (0,0) circle (0.10);%enlarged because of grounded double point 
		\begin{scope}[shift={(-3/32,-12/32)}]
			\draw [thick,black] plot  coordinates {(-1/8,1/32) (1/8,-1/32) };
		\end{scope}
		\begin{scope}[shift={(-6/32,-24/32)}]
			\draw [thick,black] plot  coordinates {(-1/8,1/32) (1/8,-1/32) };
		\end{scope}
	\end{scope}
	\begin{scope}[shift={(-4/4,1/4)}]
		\draw [thick,black] plot  coordinates {(2/24,8/24) (-5/16,-20/16) };
	\fill [red] (0,0) circle (0.10);
		\begin{scope}[shift={(-3/32,-12/32)}]
			\draw [thick,black] plot  coordinates {(-1/8,1/32) (1/8,-1/32) };
		\end{scope}
		\begin{scope}[shift={(-6/32,-24/32)}]
			\draw [thick,black] plot  coordinates {(-1/8,1/32) (1/8,-1/32) };
		\end{scope}
	\end{scope}
	\begin{scope}[shift={(-8/4,2/4)}]
		\draw [thick,black] plot  coordinates {(2/24,8/24) (-5/16,-20/16) };
		\fill [red] (0,0) circle (0.10);
		\begin{scope}[shift={(-3/32,-12/32)}]
			\draw [thick,black] plot  coordinates {(-1/8,1/32) (1/8,-1/32) };
		\end{scope}
		\begin{scope}[shift={(-6/32,-24/32)}]
			\draw [thick,black] plot  coordinates {(-1/8,1/32) (1/8,-1/32) };
		\end{scope}
	\end{scope}
\end{tikzpicture}}

\def\FigThreeSomeC{\begin{tikzpicture}[scale=0.7]
		\draw [thick,black] plot  coordinates {(4/8+4/8,-1/8-4/32) (-20/8,5/8) };
		\begin{scope}[shift={(4/8,-1/8)}]
		\draw [thick,black] plot  coordinates {(-1/32,-1/8) (1/32,1/8) };
	\end{scope}
		\begin{scope}[yscale=1, xscale=1, shift={(0,0)}]
			\draw [thick,black] plot  coordinates {(5/24,20/24) (-3/16,-12/16) };
					\fill [white] (0,0) circle (0.1);
		\draw [thick] (0,0) circle (0.1);
			\begin{scope}[shift={(-3/32,-12/32)}]
				\draw [thick,black] plot  coordinates {(-1/8,1/32) (1/8,-1/32) };
			\end{scope}
			\begin{scope}[xscale=-1, yscale=-1, shift={(-4/32,-16/32)}]
				\draw [thick,black] plot  coordinates {(-12/8,12/32) (4/8,-4/32) };
				\fill [red] (0,0) circle (0.10);
				
				\begin{scope}[shift={(-4/8,1/8)}]
					\draw [thick,black] plot  coordinates {(-1/32,-1/8) (1/32,1/8) };
				\end{scope}
				\begin{scope}[shift={(-8/8,8/32)}]
					\draw [thick,black] plot  coordinates {(-1/32,-1/8) (1/32,1/8) };
				\end{scope}
				
			\end{scope}
			
		\end{scope}
		\begin{scope}[shift={(-4/4,1/4)}]
			\draw [thick,black] plot  coordinates {(2/24,8/24) (-5/16,-20/16) };
			\fill [red] (0,0) circle (0.10);
			\begin{scope}[shift={(-3/32,-12/32)}]
				\draw [thick,black] plot  coordinates {(-1/8,1/32) (1/8,-1/32) };
			\end{scope}
			\begin{scope}[shift={(-6/32,-24/32)}]
				\draw [thick,black] plot  coordinates {(-1/8,1/32) (1/8,-1/32) };
			\end{scope}
		\end{scope}
		\begin{scope}[shift={(-8/4,2/4)}]
			\draw [thick,black] plot  coordinates {(2/24,8/24) (-5/16,-20/16) };
			\fill [red] (0,0) circle (0.10);
			\begin{scope}[shift={(-3/32,-12/32)}]
				\draw [thick,black] plot  coordinates {(-1/8,1/32) (1/8,-1/32) };
			\end{scope}
			\begin{scope}[shift={(-6/32,-24/32)}]
				\draw [thick,black] plot  coordinates {(-1/8,1/32) (1/8,-1/32) };
			\end{scope}
		\end{scope}
\end{tikzpicture}}

\def\FigThreeSomeE{\begin{tikzpicture}[scale=0.7]
		\draw [thick,black] plot  coordinates {(4/8+0/8,-1/8-0/32) (-20/8,5/8) };
		\begin{scope}[yscale=1, xscale=1, shift={(0,0)}]
			\draw [thick,black] plot  coordinates {(5/24,20/24) (-3/16,-12/16) };
				\begin{scope}[shift={(0,0)}]
			\fill [white] (0,0) circle (0.1);
			\draw [thick] (0,0) circle (0.1);
		\end{scope}
			\begin{scope}[shift={(-3/32,-12/32)}]
				\draw [thick,black] plot  coordinates {(-1/8,1/32) (1/8,-1/32) };
			\end{scope}
			\begin{scope}[xscale=-1, yscale=-1, shift={(-4/32,-16/32)}]
				\draw [thick,black] plot  coordinates {(-12/8,12/32) (4/8,-4/32) };
				\fill [red] (0,0) circle (0.10);
				
				\begin{scope}[shift={(-4/8,1/8)}]
					\draw [thick,black] plot  coordinates {(-1/32,-1/8) (1/32,1/8) };
				\end{scope}
				\begin{scope}[shift={(-8/8,8/32)}]
					\draw [thick,black] plot  coordinates {(-1/32,-1/8) (1/32,1/8) };
				\end{scope}
				
			\end{scope}
			
		\end{scope}
		\begin{scope}[shift={(-4/4,1/4)}]
			\draw [thick,black] plot  coordinates {(2/24,8/24) (-5/16,-20/16) };
			\fill [red] (0,0) circle (0.10);
			\begin{scope}[shift={(-3/32,-12/32)}]
				\draw [thick,black] plot  coordinates {(-1/8,1/32) (1/8,-1/32) };
			\end{scope}
			\begin{scope}[shift={(-6/32,-24/32)}]
				\draw [thick,black] plot  coordinates {(-1/8,1/32) (1/8,-1/32) };
			\end{scope}
		\end{scope}

\begin{scope}[shift={(-8/4,2/4)}]
	\draw [thick,black] plot  coordinates {(5/24,20/24) (-3/16,-12/16) };
		\begin{scope}[shift={(0,0)}]
	\fill [white] (0,0) circle (0.1);
	\draw [thick] (0,0) circle (0.1);
\end{scope}
	\begin{scope}[shift={(-3/32,-12/32)}]
		\draw [thick,black] plot  coordinates {(-1/8,1/32) (1/8,-1/32) };
	\end{scope}
	\begin{scope}[xscale=1, yscale=1, shift={(4/32,16/32)}]
		\draw [thick,black] plot  coordinates {(-12/8,12/32) (4/8,-4/32) };
		\fill [red] (0,0) circle (0.10);
		
		\begin{scope}[shift={(-4/8,1/8)}]
			\draw [thick,black] plot  coordinates {(-1/32,-1/8) (1/32,1/8) };
		\end{scope}
		\begin{scope}[shift={(-8/8,8/32)}]
			\draw [thick,black] plot  coordinates {(-1/32,-1/8) (1/32,1/8) };
		\end{scope}
		
	\end{scope}
	
\end{scope}
\end{tikzpicture}}

\def\FigThreeSomeD{\begin{tikzpicture}[scale=0.7]
		\draw [thick,black] plot  coordinates {(-10/8,-10/8) (1/4,1/4) };
		\draw [thick,black] plot  coordinates {(11/8,-11/8) (-1/4,1/4) };
		\begin{scope}[shift={(-1/2,-1/2)}]
			\draw [thick,black] plot  coordinates {(-1/12,1/12) (1/12,-1/12) };
		\end{scope}
		\begin{scope}[shift={(-5/6,-5/6)}]
			\draw [thick,black] plot  coordinates {(-1/12,1/12) (1/12,-1/12) };
		\end{scope}
		\begin{scope}[shift={(3/8,-3/8)}]
			\draw [thick,black] plot  coordinates {(-1/12,-1/12) (1/12,1/12) };
		\end{scope}
				\begin{scope}[shift={(3/4,-3/4)}]
			\draw [thick,black] plot  coordinates {(-1/12,-1/12) (1/12,1/12) };
		\end{scope}
		
				\fill [red] (0,0) circle (0.1);
				
				\begin{scope}[xscale=-1,shift={(-18/8,0)}]
					\draw [thick,black] plot  coordinates {(-10/8,-10/8) (1/4,1/4) };
					\draw [thick,black] plot  coordinates {(11/8,-11/8) (-1/4,1/4) };
					\begin{scope}[shift={(-1/2,-1/2)}]
						\draw [thick,black] plot  coordinates {(-1/12,1/12) (1/12,-1/12) };
					\end{scope}
					\begin{scope}[shift={(-5/6,-5/6)}]
						\draw [thick,black] plot  coordinates {(-1/12,1/12) (1/12,-1/12) };
					\end{scope}
					\begin{scope}[shift={(3/8,-3/8)}]
						\draw [thick,black] plot  coordinates {(-1/12,-1/12) (1/12,1/12) };
					\end{scope}
					\begin{scope}[shift={(3/4,-3/4)}]
						\draw [thick,black] plot  coordinates {(-1/12,-1/12) (1/12,1/12) };
					\end{scope}
					
					\fill [red] (0,0) circle (0.1);
				\end{scope}

				\fill [red] (9/8,-9/8) circle (0.1);
\end{tikzpicture}}

\def\FigThreeSomeB{\begin{tikzpicture}[scale=0.7]
		\draw [thick,black] plot  coordinates {(-10/8,-10/8) (1/4,1/4) };
		\draw [thick,black] plot  coordinates {(8/8,-8/8) (-1/4,1/4) };
		\begin{scope}[shift={(-1/2,-1/2)}]
			\draw [thick,black] plot  coordinates {(-1/12,1/12) (1/12,-1/12) };
		\end{scope}
		\begin{scope}[shift={(-5/6,-5/6)}]
			\draw [thick,black] plot  coordinates {(-1/12,1/12) (1/12,-1/12) };
		\end{scope}
		\begin{scope}[shift={(3/8,-3/8)}]
			\draw [thick,black] plot  coordinates {(-1/12,-1/12) (1/12,1/12) };
		\end{scope}
	\draw [thick,black] plot  coordinates {(4/8,-8/8) (14/8,2/8) };
		\begin{scope}[shift={(3/2-3/8,-3/8)}]
		\draw [thick,black] plot  coordinates {(-1/12,1/12) (1/12,-1/12) };
	\end{scope}
		
		\fill [red] (0,0) circle (0.1);
		\begin{scope}[shift={(3/4,-3/4)}]
			\fill [white] (0,0) circle (0.1);
			\draw [thick] (0,0) circle (0.1);
		\end{scope}
		
		\begin{scope}[xscale=-1,shift={(-24/8,0)}]
		\draw [thick,black] plot  coordinates {(-10/8,-10/8) (1/4,1/4) };
		\draw [thick,black] plot  coordinates {(8/8,-8/8) (-1/4,1/4) };
		\begin{scope}[shift={(-1/2,-1/2)}]
			\draw [thick,black] plot  coordinates {(-1/12,1/12) (1/12,-1/12) };
		\end{scope}
		\begin{scope}[shift={(-5/6,-5/6)}]
			\draw [thick,black] plot  coordinates {(-1/12,1/12) (1/12,-1/12) };
		\end{scope}
		\begin{scope}[shift={(3/8,-3/8)}]
			\draw [thick,black] plot  coordinates {(-1/12,-1/12) (1/12,1/12) };
		\end{scope}
		\draw [thick,black] plot  coordinates {(4/8,-8/8) (14/8,2/8) };
		\begin{scope}[shift={(3/2-3/8,-3/8)}]
			\draw [thick,black] plot  coordinates {(-1/12,1/12) (1/12,-1/12) };
		\end{scope}
		
		\fill [red] (0,0) circle (0.1);
		\begin{scope}[shift={(3/4,-3/4)}]
			\fill [white] (0,0) circle (0.1);
			\draw [thick] (0,0) circle (0.1);
		\end{scope}
		\end{scope}
		
			\fill [red] (12/8,0) circle (0.1);

\end{tikzpicture}}

\def\FigThreeSomeF{\begin{tikzpicture}[scale=0.7]
		\draw [thick,black] plot  coordinates {(-10/8,-10/8) (1/4,1/4) };
		\draw [thick,black] plot  coordinates {(8/8,-8/8) (-1/4,1/4) };
		\begin{scope}[shift={(-1/2,-1/2)}]
			\draw [thick,black] plot  coordinates {(-1/12,1/12) (1/12,-1/12) };
		\end{scope}
		\begin{scope}[shift={(-5/6,-5/6)}]
			\draw [thick,black] plot  coordinates {(-1/12,1/12) (1/12,-1/12) };
		\end{scope}
		\begin{scope}[shift={(3/8,-3/8)}]
			\draw [thick,black] plot  coordinates {(-1/12,-1/12) (1/12,1/12) };
		\end{scope}
		\draw [thick,black] plot  coordinates {(4/8,-8/8) (14/8,2/8) };
		\begin{scope}[shift={(3/2-3/8,-3/8)}]
			\draw [thick,black] plot  coordinates {(-1/12,1/12) (1/12,-1/12) };
		\end{scope}
		
		\fill [red] (0,0) circle (0.1);
		\begin{scope}[shift={(3/4,-3/4)}]
			\fill [white] (0,0) circle (0.1);
			\draw [thick] (0,0) circle (0.1);
		\end{scope}
		
		\begin{scope}[xscale=-1,shift={(-24/8,0)}]
			\draw [thick,black] plot  coordinates {(-10/8,-10/8) (1/4,1/4) };
						\begin{scope}[shift={(-1/2,-1/2)}]
				\draw [thick,black] plot  coordinates {(-1/12,1/12) (1/12,-1/12) };
			\end{scope}
			\begin{scope}[shift={(-5/6,-5/6)}]
				\draw [thick,black] plot  coordinates {(-1/12,1/12) (1/12,-1/12) };
			\end{scope}
			\draw [thick,black] plot  coordinates {(14/8,0) (-1/4,0) };
			\begin{scope}[shift={(4/8,0)}]
				\draw [thick,black] plot  coordinates {(0,1/8) (0,-1/8) };
			\end{scope}
				\begin{scope}[shift={(8/8,0)}]
				\draw [thick,black] plot  coordinates {(0,1/8) (0,-1/8) };
			\end{scope}
			
			\fill [red] (0,0) circle (0.1);
		\end{scope}
		
		\fill [red] (12/8,0) circle (0.1);

\end{tikzpicture}}

\def\includefig#1#2{
	\vcenter{\hbox{\makebox(120pt,55pt){
				\fbox{\makebox(0pt,45pt)[lt]{\makebox(20pt,20pt){(\small #1)}}
					\makebox(120pt,60pt){#2}}}}}}

\def\FigSomeThree{
	$\xymatrix@C+0pt@R+0pt{
		\includefig{A}{ \FigThreeSomeA} & 	\includefig{B}{\FigThreeSomeB} \\
		\includefig{C}{\FigThreeSomeC} & \includefig{D}{\FigThreeSomeD}  \\
		\includefig{E}{\FigThreeSomeE}  & \includefig{F}{\FigThreeSomeF}  }$}
\def\FigGenusThreeWithOnlyTwoGroundedDoublePoints{\begin{tikzpicture}[scale=0.7]
		\draw [thick,black] plot  coordinates {(4/8+6/8,-1/8-6/32) (-20/8,5/8) };
		\begin{scope}[shift={(4/8,-1/8)}]
			\draw [thick,black] plot  coordinates {(-1/32,-1/8) (1/32,1/8) };
		\end{scope}
		\begin{scope}[shift={(0,0)}]
			\draw [thick,black] plot  coordinates {(2/24,8/24) (-5/16,-20/16) };
			\fill [red] (0,0) circle (0.10);%enlarged because of grounded double point 
			\begin{scope}[shift={(-3/32,-12/32)}]
				\draw [thick,black] plot  coordinates {(-1/8,1/32) (1/8,-1/32) };
			\end{scope}
			\begin{scope}[shift={(-6/32,-24/32)}]
				\draw [thick,black] plot  coordinates {(-1/8,1/32) (1/8,-1/32) };
			\end{scope}
		\end{scope}
		\begin{scope}[shift={(-4/4,1/4)}]
			\draw [thick,black] plot  coordinates {(2/24,8/24) (-5/16,-20/16) };
			\fill [red] (0,0) circle (0.10);
			\begin{scope}[shift={(-3/32,-12/32)}]
				\draw [thick,black] plot  coordinates {(-1/8,1/32) (1/8,-1/32) };
			\end{scope}
			\begin{scope}[shift={(-6/32,-24/32)}]
				\draw [thick,black] plot  coordinates {(-1/8,1/32) (1/8,-1/32) };
			\end{scope}
		\end{scope}
		\begin{scope}[shift={(-8/4,2/4)}]
			\draw [thick,black] plot  coordinates {(2/24,8/24) (-7/16,-28/16) };
			\fill [black] (0,0) circle (0.08);
			\begin{scope}[shift={(-3/32,-12/32)}]
				\draw [thick,black] plot  coordinates {(-1/8,1/32) (1/8,-1/32) };
			\end{scope}
			\begin{scope}[shift={(-6/32,-24/32)}]
				\draw [thick,black] plot  coordinates {(-1/8,1/32) (1/8,-1/32) };
			\end{scope}
			\begin{scope}[shift={(-9/32,-36/32)}]
				\draw [thick,black] plot  coordinates {(-1/8,1/32) (1/8,-1/32) };
			\end{scope}
		\end{scope}
\end{tikzpicture}}

\def\FigGenusThreeOnlyTwoGroundedDP{$\vcenter{\hbox{\makebox(120pt,55pt){
				\fbox{\makebox(0pt,45pt)[lt]{\makebox(20pt,20pt){}}
					\makebox(120pt,60pt){\FigGenusThreeWithOnlyTwoGroundedDoublePoints}}}}}$}

\def\FigA{
	\def\geqone{g\kern1pt{=}\kern1pt1}
	\def\arraystretch{0.75}
	\begin{tabular}{ | c | M | M | M |}
		\hline
		\rule{0pt}{2.5ex}    $\alpha$    & $\bar C_0$ & $\bbar{C}_0$  & $C_0$ \\
		\hline
		
		$\alpha=0$ & 
		
		\begin{tikzpicture}[baseline=0, scale=0.7]
			\draw [thick,black] plot  coordinates {(0,0) (4,1) };
			\begin{scope}[shift={(1/8,1/32)}]
				\begin{scope}[shift={(3/4,3/16)}]
					\draw [thick,black] plot  coordinates {(-1/32,1/8) (1/32,-1/8) };
					\draw (-2/32,1/8) node [above] {\hbox{\footnotesize$0$}};
				\end{scope}   
				\begin{scope}[shift={(6/4,6/16)}]
					\draw [thick,black] plot  coordinates {(-1/32,1/8) (1/32,-1/8) };
					\draw (-2/32,1/8) node [above] {\hbox{\footnotesize$a$}};
				\end{scope}
				\begin{scope}[shift={(9/4,9/16)}]
					\draw [thick,black] plot  coordinates {(-1/32,1/8) (1/32,-1/8) };
					\draw (-2/32,1/8) node [above] {\hbox{\footnotesize$1$}};
				\end{scope}
				\begin{scope}[shift={(12/4,12/16)}]
					\draw [thick,black] plot  coordinates {(-1/32,1/8) (1/32,-1/8) };
					\draw (-2/32,1/8) node [above] {\hbox{\footnotesize$\infty$}};
				\end{scope}
			\end{scope}
		\end{tikzpicture} &
		
		\begin{center}\begin{tikzpicture}[baseline=0, scale=0.7]
				\draw [thick,black] plot  coordinates {(0,0) (4,1) };
				\draw [thick,blue] plot  coordinates {(0,24/12) (8/12,-8/12) };
				\fill [black] (8/17,2/17) circle (0.08);
				\begin{scope}[shift={(2/8,2/32)}]
					\begin{scope}[shift={(3/4,3/16)}]
						\draw [thick,black] plot  coordinates {(-1/32,1/8) (1/32,-1/8) };
						\draw (-2/32,1/8) node [above] {\hbox{\footnotesize$0$}};
					\end{scope}   
					\begin{scope}[shift={(6/4,6/16)}]
						\draw [thick,black] plot  coordinates {(-1/32,1/8) (1/32,-1/8) };
						\draw (-2/32,1/8) node [above] {\hbox{\footnotesize$a$}};
					\end{scope}
					\begin{scope}[shift={(9/4,9/16)}]
						\draw [thick,black] plot  coordinates {(-1/32,1/8) (1/32,-1/8) };
						\draw (-2/32,1/8) node [above] {\hbox{\footnotesize$1$}};
					\end{scope}
					\begin{scope}[shift={(12/4,12/16)}]
						\draw [thick,black] plot  coordinates {(-1/32,1/8) (1/32,-1/8) };
						\draw (-2/32,1/8) node [above] {\hbox{\footnotesize$\infty$}};
					\end{scope}
				\end{scope}
		\end{tikzpicture}\end{center}&
		
		\begin{center}\begin{tikzpicture}[baseline=0, scale=0.7]
				\draw [thick,black] plot  coordinates {(0,0) (4,1) };
				\draw [thick,blue] plot  coordinates {(0,24/12) (8/12,-8/12) };
				\draw (1/12,20/12) node [right] {\hbox{\footnotesize\,$\geqone$\;}};
				\fill [black] (8/17,2/17) circle (0.08);
				\begin{scope}[shift={(2/8,2/32)}]
					\begin{scope}[shift={(3/4,3/16)}]
						\draw [thick,black] plot  coordinates {(-1/32,1/8) (1/32,-1/8) };
						\draw (-2/32,1/8) node [above] {\hbox{\footnotesize$0$}};
					\end{scope}   
					\begin{scope}[shift={(6/4,6/16)}]
						\draw [thick,black] plot  coordinates {(-1/32,1/8) (1/32,-1/8) };
						\draw (-2/32,1/8) node [above] {\hbox{\footnotesize$a$}};
					\end{scope}
					\begin{scope}[shift={(9/4,9/16)}]
						\draw [thick,black] plot  coordinates {(-1/32,1/8) (1/32,-1/8) };
						\draw (-2/32,1/8) node [above] {\hbox{\footnotesize$1$}};
					\end{scope}
					\begin{scope}[shift={(12/4,12/16)}]
						\draw [thick,black] plot  coordinates {(-1/32,1/8) (1/32,-1/8) };
						\draw (-2/32,1/8) node [above] {\hbox{\footnotesize$\infty$}};
					\end{scope}
				\end{scope}
		\end{tikzpicture}\end{center} \\
		
		\hline
		$0 < \alpha <4$ & 
		
		\begin{center}\begin{tikzpicture}[baseline=0, scale=0.7]		
				\draw [thick,black] plot  coordinates {(-2,-2) (1/2,1/2 ) };
				\draw [thick,black] plot  coordinates {(2,-2) (-1/2,1/2 ) };
				\draw [thick,black] plot  coordinates {(3/2-1/8   +1/8,-3/2+1/8   +1/8) (3/2-1/8   -1/8,-3/2+1/8   -1/8) };
				\draw [thick,black] plot  coordinates {(1/2 +1/8    +1/8,-1/2 -1/8  +1/8) (1/2 +1/8  -1/8,-1/2-1/8   -1/8) };
				\draw [thick,black] plot  coordinates {(-1/2 -1/8    -1/8,-1/2 -1/8  +1/8) (-1/2 -1/8  +1/8,-1/2-1/8   -1/8) };
				\draw [thick,black] plot  coordinates {(-3/2+1/8   -1/8,-3/2+1/8   +1/8) (-3/2+1/8   +1/8,-3/2+1/8   -1/8) };
				\fill [black] (0,0) circle (0.08);
				\draw (-1/2 -1/8 ,-1/4-1/16 ) node [left] {\hbox{\footnotesize$0$}};
				\draw (1/2 +1/8 ,-1/4-1/16 ) node [right] {\hbox{\footnotesize$1$}};
				\draw (-3/2 + 1/8 ,-5/4+1/8 +1/16 ) node [left] {\hbox{\footnotesize$a$}};
				\draw (+3/2 - 1/8 ,-5/4+1/8 +1/16 ) node [right] {\hbox{\footnotesize$\infty$}};
		\end{tikzpicture}\end{center} &
		
		\begin{center}\begin{tikzpicture}[baseline=0, scale=0.7]		
				\draw [thick,black] plot  coordinates {(-15/8,-14/8) (-3/4,1/2 ) };
				\draw [thick,black] plot  coordinates {(15/8,-14/8) (3/4,1/2 ) };
				\draw [thick,orange] plot coordinates {(-7/4,0) (7/4,0) };
				\draw [thick,blue] plot coordinates {(0,7/4-1/4) (0,-1) };
				\fill [black] (0,0) circle (0.08);
				\fill [black] (-1,0) circle (0.08);
				\fill [black] (1,0) circle (0.08);
				\begin{scope}[shift={(-1-2/6,-2/3)}]
					\draw [thick,black] plot  coordinates {(-1/8,1/16) (1/8,-1/16) };
					\draw (-0.05,0.2) node [left] {\hbox{\footnotesize$0$}};
				\end{scope}
				\begin{scope}[shift={(-1-4/6,-4/3)}]
					\draw [thick,black] plot  coordinates {(-1/8,1/16) (1/8,-1/16) };
					\draw (-0.05,0.2) node [left] {\hbox{\footnotesize$a$}};
				\end{scope}
				\begin{scope}[shift={(1+2/6,-2/3)}]
					\draw [thick,black] plot  coordinates {(1/8,1/16) (-1/8,-1/16) };
					\draw (0.05,0.15) node [right] {\hbox{\footnotesize$1$}};
				\end{scope}
				\begin{scope}[shift={(1+4/6,-4/3)}]
					\draw [thick,black] plot  coordinates {(1/8,1/16) (-1/8,-1/16) };
					\draw (0.05,0.15) node [right] {\hbox{\footnotesize$\infty$}};
				\end{scope}
		\end{tikzpicture}\end{center} &
		
		\begin{center}\begin{tikzpicture}[baseline=0, scale=0.7]		
				\draw [thick,black] plot  coordinates {(-15/8,-14/8) (-3/4,1/2 ) };
				\draw [thick,black] plot  coordinates {(15/8,-14/8) (3/4,1/2 ) };
				\draw [thick,orange] plot coordinates {(-7/4,0) (7/4,0) };
				\draw [thick,blue] plot coordinates {(0,7/4-1/4) (0,-1) };
				\fill [black] (0,0) circle (0.08);
				\fill [black] (-1,0) circle (0.08);
				\fill [black] (1,0) circle (0.08);
				\draw (0,6/4-3/8) node [right] {\hbox{\footnotesize$\geqone$}};   
				\begin{scope}[shift={(-1-2/6,-2/3)}]
					\draw [thick,black] plot  coordinates {(-1/8,1/16) (1/8,-1/16) };
					\draw (-0.05,0.2) node [left] {\hbox{\footnotesize$0$}};
				\end{scope}
				\begin{scope}[shift={(-1-4/6,-4/3)}]
					\draw [thick,black] plot  coordinates {(-1/8,1/16) (1/8,-1/16) };
					\draw (-0.05,0.2) node [left] {\hbox{\footnotesize$a$}};
				\end{scope}
				\begin{scope}[shift={(1+2/6,-2/3)}]
					\draw [thick,black] plot  coordinates {(1/8,1/16) (-1/8,-1/16) };
					\draw (0.05,0.15) node [right] {\hbox{\footnotesize$1$}};
				\end{scope}
				\begin{scope}[shift={(1+4/6,-4/3)}]
					\draw [thick,black] plot  coordinates {(1/8,1/16) (-1/8,-1/16) };
					\draw (0.05,0.15) node [right] {\hbox{\footnotesize$\infty$}};
				\end{scope}
		\end{tikzpicture}\end{center} \\
		
		\hline
		$\alpha =4$ & 
		
		\begin{center}\begin{tikzpicture}[baseline=0, scale=0.7]
				\draw [thick,black] plot  coordinates {(-2,-2) (1/2,1/2 ) };
				\draw [thick,black] plot  coordinates {(2,-2) (-1/2,1/2 ) };
				\draw [thick,black] plot  coordinates {(3/2-1/8   +1/8,-3/2+1/8   +1/8) (3/2-1/8   -1/8,-3/2+1/8   -1/8) };
				\draw [thick,black] plot  coordinates {(1/2 +1/8    +1/8,-1/2 -1/8  +1/8) (1/2 +1/8  -1/8,-1/2-1/8   -1/8) };
				\draw [thick,black] plot  coordinates {(-1/2 -1/8    -1/8,-1/2 -1/8  +1/8) (-1/2 -1/8  +1/8,-1/2-1/8   -1/8) };
				\draw [thick,black] plot  coordinates {(-3/2+1/8   -1/8,-3/2+1/8   +1/8) (-3/2+1/8   +1/8,-3/2+1/8   -1/8) };
				\fill [black] (0,0) circle (0.08);
				\draw (-1/2 -1/8 ,-1/4-1/16 ) node [left] {\hbox{\footnotesize$0$}};
				\draw (1/2 +1/8 ,-1/4-1/16 ) node [right] {\hbox{\footnotesize$1$}};
				\draw (-3/2 + 1/8 ,-5/4+1/8 +1/16 ) node [left] {\hbox{\footnotesize$a$}};
				\draw (+3/2 - 1/8 ,-5/4+1/8 +1/16 ) node [right] {\hbox{\footnotesize$\infty$}};
		\end{tikzpicture}\end{center} &
		
		\begin{center}\begin{tikzpicture}[baseline=0, scale=0.7]		
				\draw [thick,black] plot  coordinates {(-15/8,-14/8) (-3/4,1/2 ) };
				\draw [thick,black] plot  coordinates {(15/8,-14/8) (3/4,1/2 ) };
				\draw [thick,orange] plot coordinates {(-7/4,0) (7/4,0) };
				\fill [black] (-1,0) circle (0.08);
				\fill [black] (1,0) circle (0.08);
				\begin{scope}[shift={(-1-2/6,-2/3)}]
					\draw [thick,black] plot  coordinates {(-1/8,1/16) (1/8,-1/16) };
					\draw (-0.05,0.2) node [left] {\hbox{\footnotesize$0$}};
				\end{scope}
				\begin{scope}[shift={(-1-4/6,-4/3)}]
					\draw [thick,black] plot  coordinates {(-1/8,1/16) (1/8,-1/16) };
					\draw (-0.05,0.2) node [left] {\hbox{\footnotesize$a$}};
				\end{scope}
				\begin{scope}[shift={(1+2/6,-2/3)}]
					\draw [thick,black] plot  coordinates {(1/8,1/16) (-1/8,-1/16) };
					\draw (0.05,0.15) node [right] {\hbox{\footnotesize$1$}};
				\end{scope}
				\begin{scope}[shift={(1+4/6,-4/3)}]
					\draw [thick,black] plot  coordinates {(1/8,1/16) (-1/8,-1/16) };
					\draw (0.05,0.15) node [right] {\hbox{\footnotesize$\infty$}};
				\end{scope}
		\end{tikzpicture}\end{center} &
		
		\begin{center}\begin{tikzpicture}[baseline=0, scale=0.7]		
				\draw [thick,black] plot  coordinates {(-15/8,-14/8) (-3/4,1/2 ) };
				\draw [thick,black] plot  coordinates {(15/8,-14/8) (3/4,1/2 ) };
				\draw [thick,orange] plot coordinates {(-7/4,0) (7/4,0) };
				\fill [black] (-1,0) circle (0.08);
				\fill [black] (1,0) circle (0.08);
				\draw (0,0) node [above] {\hbox{\footnotesize$\geqone$}};   
				\begin{scope}[shift={(-1-2/6,-2/3)}]
					\draw [thick,black] plot  coordinates {(-1/8,1/16) (1/8,-1/16) };
					\draw (-0.05,0.2) node [left] {\hbox{\footnotesize$0$}};
				\end{scope}
				\begin{scope}[shift={(-1-4/6,-4/3)}]
					\draw [thick,black] plot  coordinates {(-1/8,1/16) (1/8,-1/16) };
					\draw (-0.05,0.2) node [left] {\hbox{\footnotesize$a$}};
				\end{scope}
				\begin{scope}[shift={(1+2/6,-2/3)}]
					\draw [thick,black] plot  coordinates {(1/8,1/16) (-1/8,-1/16) };
					\draw (0.05,0.15) node [right] {\hbox{\footnotesize$1$}};
				\end{scope}
				\begin{scope}[shift={(1+4/6,-4/3)}]
					\draw [thick,black] plot  coordinates {(1/8,1/16) (-1/8,-1/16) };
					\draw (0.05,0.15) node [right] {\hbox{\footnotesize$\infty$}};
				\end{scope}
		\end{tikzpicture}\end{center} \\
		
		\hline 
		$\alpha >4$ & 
		
		\begin{center}\begin{tikzpicture}[baseline=0, scale=0.7]
				\draw [thick,black] plot  coordinates {(-2,-2) (1/2,1/2 ) };
				\draw [thick,black] plot  coordinates {(2,-2) (-1/2,1/2 ) };
				\draw [thick,black] plot  coordinates {(3/2-1/8   +1/8,-3/2+1/8   +1/8) (3/2-1/8   -1/8,-3/2+1/8   -1/8) };
				\draw [thick,black] plot  coordinates {(1/2 +1/8    +1/8,-1/2 -1/8  +1/8) (1/2 +1/8  -1/8,-1/2-1/8   -1/8) };
				\draw [thick,black] plot  coordinates {(-1/2 -1/8    -1/8,-1/2 -1/8  +1/8) (-1/2 -1/8  +1/8,-1/2-1/8   -1/8) };
				\draw [thick,black] plot  coordinates {(-3/2+1/8   -1/8,-3/2+1/8   +1/8) (-3/2+1/8   +1/8,-3/2+1/8   -1/8) };
				\fill [black] (0,0) circle (0.08);
				\draw (-1/2 -1/8 ,-1/4-1/16 ) node [left] {\hbox{\footnotesize$0$}};
				\draw (1/2 +1/8 ,-1/4-1/16 ) node [right] {\hbox{\footnotesize$1$}};
				\draw (-3/2 + 1/8 ,-5/4+1/8 +1/16 ) node [left] {\hbox{\footnotesize$a$}};
				\draw (+3/2 - 1/8 ,-5/4+1/8 +1/16 ) node [right] {\hbox{\footnotesize$\infty$}};
		\end{tikzpicture}\end{center} &
		
		\begin{center}\begin{tikzpicture}[baseline=0, scale=0.7]		
				\draw [thick,black] plot  coordinates {(-15/8,-14/8) (-3/4,1/2 ) };
				\draw [thick,black] plot  coordinates {(15/8,-14/8) (3/4,1/2 ) };
				\draw [thick,orange] plot coordinates {(-18/12,4/12) (6/12,-12/12) };
				\draw [thick,orange] plot coordinates {(18/12,4/12) (-6/12,-12/12) };
				\fill [black] (-1,0) circle (0.08);
				\fill [black] (1,0) circle (0.08);
				\fill [black] (0,-2/3) circle (0.08);
				\begin{scope}[shift={(-1-2/6,-2/3)}]
					\draw [thick,black] plot  coordinates {(-1/8,1/16) (1/8,-1/16) };
					\draw (-0.05,0.2) node [left] {\hbox{\footnotesize$0$}};
				\end{scope}
				\begin{scope}[shift={(-1-4/6,-4/3)}]
					\draw [thick,black] plot  coordinates {(-1/8,1/16) (1/8,-1/16) };
					\draw (-0.05,0.2) node [left] {\hbox{\footnotesize$a$}};
				\end{scope}
				\begin{scope}[shift={(1+2/6,-2/3)}]
					\draw [thick,black] plot  coordinates {(1/8,1/16) (-1/8,-1/16) };
					\draw (0.05,0.15) node [right] {\hbox{\footnotesize$1$}};
				\end{scope}
				\begin{scope}[shift={(1+4/6,-4/3)}]
					\draw [thick,black] plot  coordinates {(1/8,1/16) (-1/8,-1/16) };
					\draw (0.05,0.15) node [right] {\hbox{\footnotesize$\infty$}};
				\end{scope}
		\end{tikzpicture}\end{center} &
		
		\begin{center}\begin{tikzpicture}[baseline=0, scale=0.7]		
				\draw [thick,black] plot  coordinates {(-15/8,-14/8) (-6/8,4/8 ) };
				\draw [thick,black] plot  coordinates {(15/8,-14/8) (3/4,1/2 ) };
				\draw [thick,orange] plot [smooth, tension=1] coordinates {(-3/2,2/8) (1/4,-4/8) (-1/3,-10/8) };
				\draw [thick,orange] plot [smooth, tension=1] coordinates {(3/2,2/8) (-1/4,-4/8) (1/3,-10/8) };
				\fill [black] (0.095/2-1,0.095) circle (0.08);
				\fill [black] (-0.095/2+1,0.095) circle (0.08);
				\fill [black] (0,-0.30) circle (0.08);
				\fill [black] (0,-1.085) circle (0.08);
				\begin{scope}[shift={(-1-2/6,-2/3)}]
					\draw [thick,black] plot  coordinates {(-1/8,1/16) (1/8,-1/16) };
					\draw (-0.05,0.2) node [left] {\hbox{\footnotesize$0$}};
				\end{scope}
				\begin{scope}[shift={(-1-4/6,-4/3)}]
					\draw [thick,black] plot  coordinates {(-1/8,1/16) (1/8,-1/16) };
					\draw (-0.05,0.2) node [left] {\hbox{\footnotesize$a$}};
				\end{scope}
				\begin{scope}[shift={(1+2/6,-2/3)}]
					\draw [thick,black] plot  coordinates {(1/8,1/16) (-1/8,-1/16) };
					\draw (0.05,0.15) node [right] {\hbox{\footnotesize$1$}};
				\end{scope}
				\begin{scope}[shift={(1+4/6,-4/3)}]
					\draw [thick,black] plot  coordinates {(1/8,1/16) (-1/8,-1/16) };
					\draw (0.05,0.15) node [right] {\hbox{\footnotesize$\infty$}};
				\end{scope}
		\end{tikzpicture}\end{center} \\
		
		\hline
\end{tabular}}

	\def\FigXten{
		\begin{tikzpicture}[
			outer frame sep=0ex,%
			/tikz/inner frame xsep=2ex,
			/tikz/inner frame ysep=2ex,
			show background top,%
			show background bottom,%
			show background left,%
			show background right]
			\draw[thick] plot (0,1/4) -- (0,-5-1/4);

			\begin{scope}[xscale=1, shift={(0,0)}]
			\draw[thick,orange] plot (-1/4,-1/16) -- (2+1/4,1/2+1/16);
			\fill [black]   (0,0) circle (0.05);%
			\begin{scope}[xscale=-1, shift={(-1,1/4)}]
				\draw[thick,blue] plot (-1/16,-1/4) -- (1/4+1/8,1+1/2);
				\fill [black]   (0,0) circle (0.05);%
				\draw (3/8-1/8,3/2-1/4) node [right] {\hbox{\footnotesize$\geqtwo$}};
			\end{scope}
			
			\begin{scope}[xscale=-1, shift={(-4,0)}]
				\draw[thick,orange] plot (-1/4,-1/16) -- (2+1/4,1/2+1/16);
				\fill [black]   (0,0) circle (0.05);%
				\begin{scope}[xscale=-1, shift={(-1,1/4)}]
					\draw[thick,blue] plot (-1/16,-1/4) -- (1/4+1/8,1+1/2);
					\fill [black]   (0,0) circle (0.05);%
					\draw (3/8,3/2-1/4) node [right] {\hbox{\footnotesize$\geqone$}};
				\end{scope}
				
			\end{scope}
			\fill [black]   (2,1/2) circle (0.05);
			\begin{scope}[shift={(4,0)}]
				\draw[thick] plot (0,1/4) -- (0,-3/2+1/4);
				
				\begin{scope}[shift={(0,-1/2)}]
					\fill [black]   (0,0) circle (0.05);
					\draw[thick] plot (-2*3/4,-1/2*3/4) -- (1/4,1/16);
					\begin{scope}[shift={(-1/2,-1/8)}]
						\draw[thick] plot (-1/32,1/8) -- (1/32,-1/8);
					\end{scope}
					\begin{scope}[shift={(-1,-1/4)}]
						\draw[thick] plot (-1/32,1/8) -- (1/32,-1/8);
					\end{scope}
				\end{scope}
				
				\begin{scope}[shift={(0,-1)}]
					\fill [black]   (0,0) circle (0.05);
					\draw[thick] plot (-2*3/4,-1/2*3/4) -- (1/4,1/16);
					\begin{scope}[shift={(-1/2,-1/8)}]
						\draw[thick] plot (-1/32,1/8) -- (1/32,-1/8);
					\end{scope}
					\begin{scope}[shift={(-1,-1/4)}]
						\draw[thick] plot (-1/32,1/8) -- (1/32,-1/8);
					\end{scope}
				\end{scope}
				
			\end{scope}
			\begin{scope}

			\end{scope}

		\end{scope}
						
						\begin{scope}[xscale=-1, shift={(0,-3/2-1/4)}]
							\draw[thick,orange] plot (-1/4,-1/16) -- (2+1/4,1/2+1/16);
							\fill [black]   (0,0) circle (0.05);%
							\begin{scope}[xscale=-1, shift={(-1,1/4)}]
								\draw[thick,blue] plot (-1/16,-1/4) -- (1/4+1/8,1+1/2);
								\fill [black]   (0,0) circle (0.05);%
								\draw (3/8,3/2-1/4) node [left] {\hbox{\footnotesize$\geqtwo$}};
							\end{scope}
							
							\begin{scope}[xscale=-1, shift={(-4,0)}]
								\draw[thick,orange] plot (-1/4,-1/16) -- (2+1/4,1/2+1/16);
								\fill [black]   (0,0) circle (0.05);%
								\begin{scope}[xscale=-1, shift={(-1,1/4)}]
									\draw[thick,blue] plot (-1/16,-1/4) -- (1/4+1/8,1+1/2);
									\fill [black]   (0,0) circle (0.05);%
									\draw (3/8,3/2-1/4) node [left] {\hbox{\footnotesize$\geqone$}};
								\end{scope}
								
							\end{scope}
							\fill [black]   (2,1/2) circle (0.05);
							\begin{scope}[shift={(4,0)}]
								\draw[thick] plot (0,1/4) -- (0,-3/2+1/4);
								
								\begin{scope}[shift={(0,-1/2)}]
									\fill [black]   (0,0) circle (0.05);
									\draw[thick] plot (-2*3/4,-1/2*3/4) -- (1/4,1/16);
									\begin{scope}[shift={(-1/2,-1/8)}]
										\draw[thick] plot (-1/32,1/8) -- (1/32,-1/8);
									\end{scope}
									\begin{scope}[shift={(-1,-1/4)}]
										\draw[thick] plot (-1/32,1/8) -- (1/32,-1/8);
									\end{scope}
								\end{scope}
								
								\begin{scope}[shift={(0,-1)}]
									\fill [black]   (0,0) circle (0.05);
									\draw[thick] plot (-2*3/4,-1/2*3/4) -- (1/4,1/16);
									\begin{scope}[shift={(-1/2,-1/8)}]
										\draw[thick] plot (-1/32,1/8) -- (1/32,-1/8);
									\end{scope}
									\begin{scope}[shift={(-1,-1/4)}]
										\draw[thick] plot (-1/32,1/8) -- (1/32,-1/8);
									\end{scope}
								\end{scope}
								
							\end{scope}
							\begin{scope}

							\end{scope}

						\end{scope}
						
						\begin{scope}[xscale=1,yscale=1, shift={(0,-3-1/2)}]
							\draw[thick,orange] plot (-1/4,-1/16) -- (2+1/4,1/2+1/16);
							\fill [black]   (0,0) circle (0.05);%
							\begin{scope}[xscale=-1, shift={(-1,1/4)}]
								\draw[thick,blue] plot (-1/16,-1/4) -- (1/4+1/8,1+1/2);
								\fill [black]   (0,0) circle (0.05);%
								\draw (3/8-1/8,3/2-1/4) node [right] {\hbox{\footnotesize$\geqtwo$}};
							\end{scope}
							
							\begin{scope}[xscale=-1, shift={(-4,0)}]
								\draw[thick,orange] plot (-1/4,-1/16) -- (2+1/4,1/2+1/16);
								\fill [black]   (0,0) circle (0.05);%
								\begin{scope}[xscale=-1, shift={(-1,1/4)}]
									\draw[thick,blue] plot (-1/16,-1/4) -- (1/4+1/8,1+1/2);
									\fill [black]   (0,0) circle (0.05);%
									\draw (3/8,3/2-1/4) node [right] {\hbox{\footnotesize$\geqone$}};
								\end{scope}
								
							\end{scope}
							\fill [black]   (2,1/2) circle (0.05);
							\begin{scope}[shift={(4,0)}]
								\draw[thick] plot (0,1/4) -- (0,-3/2+1/4);
								
								\begin{scope}[shift={(0,-1/2)}]
									\fill [black]   (0,0) circle (0.05);
									\draw[thick] plot (-2*3/4,-1/2*3/4) -- (1/4,1/16);
									\begin{scope}[shift={(-1/2,-1/8)}]
										\draw[thick] plot (-1/32,1/8) -- (1/32,-1/8);
									\end{scope}
									\begin{scope}[shift={(-1,-1/4)}]
										\draw[thick] plot (-1/32,1/8) -- (1/32,-1/8);
									\end{scope}
								\end{scope}
								
								\begin{scope}[shift={(0,-1)}]
									\fill [black]   (0,0) circle (0.05);
									\draw[thick] plot (-2*3/4,-1/2*3/4) -- (1/4,1/16);
									\begin{scope}[shift={(-1/2,-1/8)}]
										\draw[thick] plot (-1/32,1/8) -- (1/32,-1/8);
									\end{scope}
									\begin{scope}[shift={(-1,-1/4)}]
										\draw[thick] plot (-1/32,1/8) -- (1/32,-1/8);
									\end{scope}
								\end{scope}
								
							\end{scope}
							\begin{scope}

							\end{scope}
						\end{scope}
						
						\begin{scope}[shift={(0,-3)}]
								\begin{scope}[shift={(0,-1)}]
								\fill [black]   (0,0) circle (0.05);
								\draw[thick] plot (-5/2,-1/2-1/8) -- (1/4,1/16);
								\begin{scope}[shift={(-1/2,-1/8)}]
									\draw[thick] plot (-1/32,1/8) -- (1/32,-1/8);
								\end{scope}
								\begin{scope}[shift={(-1,-1/4)}]
									\draw[thick] plot (-1/32,1/8) -- (1/32,-1/8);
								\end{scope}
									\begin{scope}[shift={(-3/2,-3/8)}]
									\draw[thick] plot (-1/32,1/8) -- (1/32,-1/8);
								\end{scope}
								\begin{scope}[shift={(-4/2,-4/8)}]
									\draw[thick] plot (-1/32,1/8) -- (1/32,-1/8);
								\end{scope}
							\end{scope}
						\end{scope}
						\begin{scope}[shift={(0,-3-1/2)}]
							\begin{scope}[shift={(0,-1)}]
								\fill [black]   (0,0) circle (0.05);
								\draw[thick] plot (-5/2,-1/2-1/8) -- (1/4,1/16);
								\begin{scope}[shift={(-1/2,-1/8)}]
									\draw[thick] plot (-1/32,1/8) -- (1/32,-1/8);
								\end{scope}
								\begin{scope}[shift={(-1,-1/4)}]
									\draw[thick] plot (-1/32,1/8) -- (1/32,-1/8);
								\end{scope}
								\begin{scope}[shift={(-3/2,-3/8)}]
									\draw[thick] plot (-1/32,1/8) -- (1/32,-1/8);
								\end{scope}
								\begin{scope}[shift={(-4/2,-4/8)}]
									\draw[thick] plot (-1/32,1/8) -- (1/32,-1/8);
								\end{scope}
							\end{scope}
						\end{scope}
		
		\end{tikzpicture}
}

\def\FigXeleven{
	\begin{tikzpicture}[
		outer frame sep=0ex,%
		/tikz/inner frame xsep=2ex,
		/tikz/inner frame ysep=2ex,
		show background top,%
		show background bottom,%
		show background left,%
		show background right]
		\draw[thick] plot (0,1/4) -- (0,-7-1/4);

		\begin{scope}[xscale=1, shift={(0,0)}]
			\draw[thick] plot (-1/4,-1/16) -- (2+1/4,1/2+1/16);
			\fill [black]   (0,0) circle (0.05);%
			\begin{scope}[xscale=-1, shift={(-1,1/4)}]
				\draw[thick,black] plot (-1/16,-1/4) -- (1/4+1/8,1+1/2);
				\fill [black]   (0,0) circle (0.05);%
			\begin{scope}[shift={(1/8,1/2)}]
					\draw[thick] plot (1/8,-1/32) --(-1/8,1/32);
			\end{scope}
				\begin{scope}[shift={(2/8,2/2)}]
				\draw[thick] plot (1/8,-1/32) --(-1/8,1/32);
			\end{scope}
			\end{scope}

			\fill [black]   (2,1/2) circle (0.05);
			\begin{scope}[shift={(2,2/4)}]
				\draw[thick] plot (0,1/4) -- (0,-5/2+1/4);
				
				\begin{scope}[xscale=-1,shift={(0,-1/2)}]
					\fill [black]   (0,0) circle (0.05);
					\draw[thick] plot (-2*3/4,-1/2*3/4) -- (1/4,1/16);
					\begin{scope}[shift={(-1/2,-1/8)}]
						\draw[thick] plot (-1/32,1/8) -- (1/32,-1/8);
					\end{scope}
					\begin{scope}[shift={(-1,-1/4)}]
						\draw[thick] plot (-1/32,1/8) -- (1/32,-1/8);
					\end{scope}
				\end{scope}
				
				\begin{scope}[xscale=-1,shift={(0,-1)}]
					\fill [black]   (0,0) circle (0.05);
					\draw[thick] plot (-2*3/4,-1/2*3/4) -- (1/4,1/16);
					\begin{scope}[shift={(-1/2,-1/8)}]
						\draw[thick] plot (-1/32,1/8) -- (1/32,-1/8);
					\end{scope}
					\begin{scope}[shift={(-1,-1/4)}]
						\draw[thick] plot (-1/32,1/8) -- (1/32,-1/8);
					\end{scope}
				\end{scope}
				
				\begin{scope}[xscale=-1,shift={(0,-3/2)}]
					\draw[thick, blue] plot (-2*3/4,-1/2*3/4) -- (1/4,1/16);
					\fill [black]   (0,0) circle (0.05);
					\draw (-3/2,-3/8) node [right] {\hbox{\footnotesize$\geqone$}};
				\end{scope}
				
				\begin{scope}[xscale=1, yscale=-1,shift={(0,4/2)}]
					\draw[thick, blue] plot (-2*3/4,-1/2*3/4) -- (1/4,1/16);
					\fill [black]   (0,0) circle (0.05);
					\draw (-3/2,-3/8) node [above] {\hbox{\footnotesize$\geqone$}};
				\end{scope}

			\end{scope}
			\begin{scope}

			\end{scope}

		\end{scope}
		
		\begin{scope}[xscale=-1, shift={(0,-2)}]
			\draw[thick] plot (-1/4,-1/16) -- (2+1/4,1/2+1/16);
			\fill [black]   (0,0) circle (0.05);%
			\begin{scope}[xscale=-1, shift={(-1,1/4)}]
				\draw[thick,black] plot (-1/16,-1/4) -- (1/4+1/8,1+1/2);
				\fill [black]   (0,0) circle (0.05);%
				\begin{scope}[shift={(1/8,1/2)}]
					\draw[thick] plot (1/8,-1/32) --(-1/8,1/32);
				\end{scope}
				\begin{scope}[shift={(2/8,2/2)}]
					\draw[thick] plot (1/8,-1/32) --(-1/8,1/32);
				\end{scope}
			\end{scope}

			\fill [black]   (2,1/2) circle (0.05);
			\begin{scope}[shift={(2,2/4)}]
				\draw[thick] plot (0,1/4) -- (0,-5/2+1/4);
				
				\begin{scope}[xscale=-1,shift={(0,-1/2)}]
					\fill [black]   (0,0) circle (0.05);
					\draw[thick] plot (-2*3/4,-1/2*3/4) -- (1/4,1/16);
					\begin{scope}[shift={(-1/2,-1/8)}]
						\draw[thick] plot (-1/32,1/8) -- (1/32,-1/8);
					\end{scope}
					\begin{scope}[shift={(-1,-1/4)}]
						\draw[thick] plot (-1/32,1/8) -- (1/32,-1/8);
					\end{scope}
				\end{scope}
				
				\begin{scope}[xscale=-1,shift={(0,-1)}]
					\fill [black]   (0,0) circle (0.05);
					\draw[thick] plot (-2*3/4,-1/2*3/4) -- (1/4,1/16);
					\begin{scope}[shift={(-1/2,-1/8)}]
						\draw[thick] plot (-1/32,1/8) -- (1/32,-1/8);
					\end{scope}
					\begin{scope}[shift={(-1,-1/4)}]
						\draw[thick] plot (-1/32,1/8) -- (1/32,-1/8);
					\end{scope}
				\end{scope}
				
				\begin{scope}[xscale=-1,shift={(0,-3/2)}]
					\draw[thick, blue] plot (-2*3/4,-1/2*3/4) -- (1/4,1/16);
					\fill [black]   (0,0) circle (0.05);
					\draw (-3/2,-3/8) node [below] {\hbox{\footnotesize$\geqone$}};
				\end{scope}
				
				\begin{scope}[xscale=1, yscale=-1,shift={(0,4/2)}]
					\draw[thick, blue] plot (-2*3/4,-1/2*3/4) -- (1/4,1/16);
					\fill [black]   (0,0) circle (0.05);
					\draw (-3/2,-3/8) node [above] {\hbox{\footnotesize$\geqone$}};
				\end{scope}

			\end{scope}
			\begin{scope}

			\end{scope}

		\end{scope}
		
		\begin{scope}[xscale=1, shift={(0,-4)}]
			\draw[thick] plot (-1/4,-1/16) -- (2+1/4,1/2+1/16);
			\fill [black]   (0,0) circle (0.05);%
			\begin{scope}[xscale=-1, shift={(-1,1/4)}]
				\draw[thick,black] plot (-1/16,-1/4) -- (1/4+1/8,1+1/2);
				\fill [black]   (0,0) circle (0.05);%
				\begin{scope}[shift={(1/8,1/2)}]
					\draw[thick] plot (1/8,-1/32) --(-1/8,1/32);
				\end{scope}
				\begin{scope}[shift={(2/8,2/2)}]
					\draw[thick] plot (1/8,-1/32) --(-1/8,1/32);
				\end{scope}
			\end{scope}

			\fill [black]   (2,1/2) circle (0.05);
			\begin{scope}[shift={(2,2/4)}]
				\draw[thick] plot (0,1/4) -- (0,-5/2+1/4);
				
				\begin{scope}[xscale=-1,shift={(0,-1/2)}]
					\fill [black]   (0,0) circle (0.05);
					\draw[thick] plot (-2*3/4,-1/2*3/4) -- (1/4,1/16);
					\begin{scope}[shift={(-1/2,-1/8)}]
						\draw[thick] plot (-1/32,1/8) -- (1/32,-1/8);
					\end{scope}
					\begin{scope}[shift={(-1,-1/4)}]
						\draw[thick] plot (-1/32,1/8) -- (1/32,-1/8);
					\end{scope}
				\end{scope}
				
				\begin{scope}[xscale=-1,shift={(0,-1)}]
					\fill [black]   (0,0) circle (0.05);
					\draw[thick] plot (-2*3/4,-1/2*3/4) -- (1/4,1/16);
					\begin{scope}[shift={(-1/2,-1/8)}]
						\draw[thick] plot (-1/32,1/8) -- (1/32,-1/8);
					\end{scope}
					\begin{scope}[shift={(-1,-1/4)}]
						\draw[thick] plot (-1/32,1/8) -- (1/32,-1/8);
					\end{scope}
				\end{scope}
				
				\begin{scope}[xscale=-1,shift={(0,-3/2)}]
					\draw[thick, blue] plot (-2*3/4,-1/2*3/4) -- (1/4,1/16);
					\fill [black]   (0,0) circle (0.05);
					\draw (-3/2,-3/8) node [right] {\hbox{\footnotesize$\geqone$}};
				\end{scope}
				
				\begin{scope}[xscale=1, yscale=-1,shift={(0,4/2)}]
					\draw[thick, blue] plot (-2*3/4,-1/2*3/4) -- (1/4,1/16);
					\fill [black]   (0,0) circle (0.05);
					\draw (-3/2,-3/8) node [above] {\hbox{\footnotesize$\geqone$}};
				\end{scope}

			\end{scope}
			\begin{scope}

			\end{scope}

		\end{scope}

		\begin{scope}[shift={(0,-6)}]
						\draw[thick, orange] plot (-8/2,-8/8) -- (1/4,1/16);
									\fill [black]   (0,0) circle (0.05);
	\begin{scope}[shift={(-1/2,-1/8)}]
		\draw[thick, blue] plot (-3/8,3/2) -- (1/16,-1/4);
			\fill [black]   (0,0) circle (0.05);
				\draw (-3/8,3/2) node [above] {\hbox{\footnotesize$\geqone$}};
		\end{scope}
		\begin{scope}[shift={(-10/4,-10/16)}]
			\draw[thick, blue] plot (-3/8,3/2) -- (1/16,-1/4);
			\fill [black]   (0,0) circle (0.05);
						\draw (-3/8,3/2) node [above] {\hbox{\footnotesize$\geqone$}};
		\end{scope}
		\begin{scope}[shift={(-6/4,-6/16)}]
			\draw[thick, blue] plot (-3/8,3/2) -- (1/16,-1/4);
			\fill [black]   (0,0) circle (0.05);
						\draw (-3/8,3/2) node [above] {\hbox{\footnotesize$\geqone$}};
		\end{scope}
			\begin{scope}[shift={(-14/4,-14/16)}]
			\draw[thick] plot (-5/8,5/2) -- (1/16,-1/4);
			\fill [black]   (0,0) circle (0.05);
			\begin{scope}[shift={(-1/8,1/2)}]
				\draw[thick] plot (-1/8,-1/32) -- (1/8,1/32);
			\end{scope}
			\begin{scope}[shift={(-2/8,2/2)}]
				\draw[thick] plot (-1/8,-1/32) -- (1/8,1/32);
			\end{scope}
			\begin{scope}[shift={(-3/8,3/2)}]
				\draw[thick] plot (-1/8,-1/32) -- (1/8,1/32);
			\end{scope}
			\begin{scope}[shift={(-4/8,4/2)}]
				\draw[thick] plot (-1/8,-1/32) -- (1/8,1/32);
			\end{scope}
		\end{scope}
	\end{scope}
	
	\begin{scope}[xscale=-1, yscale=-1,shift={(0,7)}]
		\draw[thick, orange] plot (-8/2,-8/8) -- (1/4,1/16);
		\fill [black]   (0,0) circle (0.05);
		\begin{scope}[shift={(-1/2,-1/8)}]
			\draw[thick, blue] plot (-3/8,3/2) -- (1/16,-1/4);
			\fill [black]   (0,0) circle (0.05);
			\draw (-3/8,3/2) node [below] {\hbox{\footnotesize$\geqone$}};
		\end{scope}
		\begin{scope}[shift={(-10/4,-10/16)}]
			\draw[thick, blue] plot (-3/8,3/2) -- (1/16,-1/4);
			\fill [black]   (0,0) circle (0.05);
			\draw (-3/8,3/2) node [below] {\hbox{\footnotesize$\geqone$}};
		\end{scope}
		\begin{scope}[shift={(-6/4,-6/16)}]
			\draw[thick, blue] plot (-3/8,3/2) -- (1/16,-1/4);
			\fill [black]   (0,0) circle (0.05);
			\draw (-3/8,3/2) node [below] {\hbox{\footnotesize$\geqone$}};
		\end{scope}
		\begin{scope}[shift={(-14/4,-14/16)}]
			\draw[thick] plot (-5/8,5/2) -- (1/16,-1/4);
			\fill [black]   (0,0) circle (0.05);
			\begin{scope}[shift={(-1/8,1/2)}]
				\draw[thick] plot (-1/8,-1/32) -- (1/8,1/32);
			\end{scope}
			\begin{scope}[shift={(-2/8,2/2)}]
				\draw[thick] plot (-1/8,-1/32) -- (1/8,1/32);
			\end{scope}
			\begin{scope}[shift={(-3/8,3/2)}]
				\draw[thick] plot (-1/8,-1/32) -- (1/8,1/32);
			\end{scope}
			\begin{scope}[shift={(-4/8,4/2)}]
				\draw[thick] plot (-1/8,-1/32) -- (1/8,1/32);
			\end{scope}
		\end{scope}
	\end{scope}
	
	\end{tikzpicture}
}

\def\FigXsixteen{
	\begin{tikzpicture}[
		outer frame sep=0ex,%
		/tikz/inner frame xsep=2ex,
		/tikz/inner frame ysep=2ex,
		show background top,%
		show background bottom,%
		show background left,%
		show background right]
		\draw[thick] plot (0,1/4) -- (0,-6-1/2);
	\draw [thick, black] (0,1/4) arc (0:90:0.25);
	\draw [thick] plot (-1/4,1/2) --(-3-1/2,1/2);
	\begin{scope}[shift={(-1/2,1/2)}]
	\draw[thick] (0,-1/8) -- (0,1/8);
	\end{scope}
	\begin{scope}[shift={(-2/2,1/2)}]
		\draw[thick] (0,-1/8) -- (0,1/8);
	\end{scope}
	\begin{scope}[shift={(-3/2,1/2)}]
		\draw[thick] (0,-1/8) -- (0,1/8);
	\end{scope}
	\begin{scope}[shift={(-4/2,1/2)}]
		\draw[thick] (0,-1/8) -- (0,1/8);
	\end{scope}
	\begin{scope}[shift={(-5/2,1/2)}]
		\draw[thick] (0,-1/8) -- (0,1/8);
	\end{scope}
	\begin{scope}[shift={(-6/2,1/2)}]
		\draw[thick] (0,-1/8) -- (0,1/8);
	\end{scope}
	
	\begin{scope}[xscale=1, yscale=-1, shift={(0,6+1/4)}]
	\draw [thick, black] (0,1/4) arc (0:90:0.25);
\draw [thick] plot (-1/4,1/2) --(-3-1/2,1/2);
\begin{scope}[shift={(-1/2,1/2)}]
	\draw[thick] (0,-1/8) -- (0,1/8);
\end{scope}
\begin{scope}[shift={(-2/2,1/2)}]
	\draw[thick] (0,-1/8) -- (0,1/8);
\end{scope}
\begin{scope}[shift={(-3/2,1/2)}]
	\draw[thick] (0,-1/8) -- (0,1/8);
\end{scope}
\begin{scope}[shift={(-4/2,1/2)}]
	\draw[thick] (0,-1/8) -- (0,1/8);
\end{scope}
\begin{scope}[shift={(-5/2,1/2)}]
	\draw[thick] (0,-1/8) -- (0,1/8);
\end{scope}
\begin{scope}[shift={(-6/2,1/2)}]
	\draw[thick] (0,-1/8) -- (0,1/8);
\end{scope}
	\end{scope}

	\begin{scope}[xscale=-1, yscale=-1, shift={(0,0)}]
		\draw[thick, orange] plot (-8/2,-8/8) -- (1/4,1/16);
		\fill [black]   (0,0) circle (0.05);
		\begin{scope}[shift={(-1/2,-1/8)}]
			\draw[thick, blue] plot (-3/8,3/2) -- (1/16,-1/4);
			\fill [black]   (0,0) circle (0.05);
			\draw (-3/8,3/2) node [below] {\hbox{\footnotesize$\geqone$}};
		\end{scope}
		\begin{scope}[shift={(-10/4,-10/16)}]
			\draw[thick, blue] plot (-3/8,3/2) -- (1/16,-1/4);
			\fill [black]   (0,0) circle (0.05);
			\draw (-3/8,3/2) node [below] {\hbox{\footnotesize$\geqone$}};
		\end{scope}
		\begin{scope}[shift={(-6/4,-6/16)}]
			\draw[thick, blue] plot (-3/8,3/2) -- (1/16,-1/4);
			\fill [black]   (0,0) circle (0.05);
			\draw (-3/8,3/2) node [below] {\hbox{\footnotesize$\geqone$}};
		\end{scope}
		\begin{scope}[shift={(-14/4,-14/16)}]
			\draw[thick] plot (-5/8+1/16,5/2-1/4) -- (1/16,-1/4);
			\fill [black]   (0,0) circle (0.05);
			\begin{scope}[shift={(-1/8,1/2)}]
				\draw[thick] plot (-12/8,-12/32) -- (2/8,2/32);
					\fill [black]   (0,0) circle (0.05);
					\begin{scope}[shift={(-1/2,-1/8)}]
							\draw[thick] plot (1/32,-1/8) -- (-1/32,1/8);
					\end{scope}
						\begin{scope}[shift={(-2/2,-2/8)}]
						\draw[thick] plot (1/32,-1/8) -- (-1/32,1/8);
					\end{scope}
			\end{scope}
			\begin{scope}[shift={(-2/8,2/2)}]
				\draw[thick] plot (-12/8,-12/32) -- (2/8,2/32);
				\fill [black]   (0,0) circle (0.05);
				\begin{scope}[shift={(-1/2,-1/8)}]
					\draw[thick] plot (1/32,-1/8) -- (-1/32,1/8);
				\end{scope}
				\begin{scope}[shift={(-2/2,-2/8)}]
					\draw[thick] plot (1/32,-1/8) -- (-1/32,1/8);
				\end{scope}
			\end{scope}
			\begin{scope}[shift={(-3/8,3/2)}]
					\draw[thick] plot (-12/8,-12/32) -- (2/8,2/32);
				\fill [black]   (0,0) circle (0.05);
				\begin{scope}[shift={(-1/2,-1/8)}]
					\draw[thick] plot (1/32,-1/8) -- (-1/32,1/8);
				\end{scope}
				\begin{scope}[shift={(-2/2,-2/8)}]
					\draw[thick] plot (1/32,-1/8) -- (-1/32,1/8);
				\end{scope}
			\end{scope}
			\begin{scope}[shift={(-4/8,4/2)}]
			\draw[thick, blue] plot (-12/8,-12/32) -- (2/8,2/32);
			\fill [black]   (0,0) circle (0.05);
	\draw (-8/8,-8/32) node [below] {\hbox{\footnotesize$\geqone$}};
			\end{scope}
		\end{scope}
	\end{scope}

		\begin{scope}[xscale=1, yscale=-1, shift={(0,3/2)}]
		\draw[thick, orange] plot (-8/2,-8/8) -- (1/4,1/16);
		\fill [black]   (0,0) circle (0.05);
		\begin{scope}[shift={(-1/2,-1/8)}]
			\draw[thick, blue] plot (-3/8,3/2) -- (1/16,-1/4);
			\fill [black]   (0,0) circle (0.05);
			\draw (-3/8,3/2) node [below] {\hbox{\footnotesize$\geqone$}};
		\end{scope}
		\begin{scope}[shift={(-10/4,-10/16)}]
			\draw[thick, blue] plot (-3/8,3/2) -- (1/16,-1/4);
			\fill [black]   (0,0) circle (0.05);
			\draw (-3/8,3/2) node [below] {\hbox{\footnotesize$\geqone$}};
		\end{scope}
		\begin{scope}[shift={(-6/4,-6/16)}]
			\draw[thick, blue] plot (-3/8,3/2) -- (1/16,-1/4);
			\fill [black]   (0,0) circle (0.05);
			\draw (-3/8,3/2) node [below] {\hbox{\footnotesize$\geqone$}};
		\end{scope}
		\begin{scope}[shift={(-14/4,-14/16)}]
			\draw[thick] plot (-5/8+1/16,5/2-1/4) -- (1/16,-1/4);
			\fill [black]   (0,0) circle (0.05);
			\begin{scope}[shift={(-1/8,1/2)}]
				\draw[thick] plot (-12/8,-12/32) -- (2/8,2/32);
				\fill [black]   (0,0) circle (0.05);
				\begin{scope}[shift={(-1/2,-1/8)}]
					\draw[thick] plot (1/32,-1/8) -- (-1/32,1/8);
				\end{scope}
				\begin{scope}[shift={(-2/2,-2/8)}]
					\draw[thick] plot (1/32,-1/8) -- (-1/32,1/8);
				\end{scope}
			\end{scope}
			\begin{scope}[shift={(-2/8,2/2)}]
				\draw[thick] plot (-12/8,-12/32) -- (2/8,2/32);
				\fill [black]   (0,0) circle (0.05);
				\begin{scope}[shift={(-1/2,-1/8)}]
					\draw[thick] plot (1/32,-1/8) -- (-1/32,1/8);
				\end{scope}
				\begin{scope}[shift={(-2/2,-2/8)}]
					\draw[thick] plot (1/32,-1/8) -- (-1/32,1/8);
				\end{scope}
			\end{scope}
			\begin{scope}[shift={(-3/8,3/2)}]
				\draw[thick] plot (-12/8,-12/32) -- (2/8,2/32);
				\fill [black]   (0,0) circle (0.05);
				\begin{scope}[shift={(-1/2,-1/8)}]
					\draw[thick] plot (1/32,-1/8) -- (-1/32,1/8);
				\end{scope}
				\begin{scope}[shift={(-2/2,-2/8)}]
					\draw[thick] plot (1/32,-1/8) -- (-1/32,1/8);
				\end{scope}
			\end{scope}
			\begin{scope}[shift={(-4/8,4/2)}]
				\draw[thick, blue] plot (-12/8,-12/32) -- (2/8,2/32);
				\fill [black]   (0,0) circle (0.05);
				\draw (-8/8,-8/32) node [below] {\hbox{\footnotesize$\geqone$}};
				
			\end{scope}
		\end{scope}
	\end{scope}

		\begin{scope}[xscale=-1, yscale=-1, shift={(0,3)}]
		\draw[thick, orange] plot (-8/2,-8/8) -- (1/4,1/16);
		\fill [black]   (0,0) circle (0.05);
		\begin{scope}[shift={(-1/2,-1/8)}]
			\draw[thick, blue] plot (-3/8,3/2) -- (1/16,-1/4);
			\fill [black]   (0,0) circle (0.05);
			\draw (-3/8,3/2) node [below] {\hbox{\footnotesize$\geqone$}};
		\end{scope}
		\begin{scope}[shift={(-10/4,-10/16)}]
			\draw[thick, blue] plot (-3/8,3/2) -- (1/16,-1/4);
			\fill [black]   (0,0) circle (0.05);
			\draw (-3/8,3/2) node [below] {\hbox{\footnotesize$\geqone$}};
		\end{scope}
		\begin{scope}[shift={(-6/4,-6/16)}]
			\draw[thick, blue] plot (-3/8,3/2) -- (1/16,-1/4);
			\fill [black]   (0,0) circle (0.05);
			\draw (-3/8,3/2) node [below] {\hbox{\footnotesize$\geqone$}};
		\end{scope}
		\begin{scope}[shift={(-14/4,-14/16)}]
			\draw[thick] plot (-5/8+1/16,5/2-1/4) -- (1/16,-1/4);
			\fill [black]   (0,0) circle (0.05);
			\begin{scope}[shift={(-1/8,1/2)}]
				\draw[thick] plot (-12/8,-12/32) -- (2/8,2/32);
				\fill [black]   (0,0) circle (0.05);
				\begin{scope}[shift={(-1/2,-1/8)}]
					\draw[thick] plot (1/32,-1/8) -- (-1/32,1/8);
				\end{scope}
				\begin{scope}[shift={(-2/2,-2/8)}]
					\draw[thick] plot (1/32,-1/8) -- (-1/32,1/8);
				\end{scope}
			\end{scope}
			\begin{scope}[shift={(-2/8,2/2)}]
				\draw[thick] plot (-12/8,-12/32) -- (2/8,2/32);
				\fill [black]   (0,0) circle (0.05);
				\begin{scope}[shift={(-1/2,-1/8)}]
					\draw[thick] plot (1/32,-1/8) -- (-1/32,1/8);
				\end{scope}
				\begin{scope}[shift={(-2/2,-2/8)}]
					\draw[thick] plot (1/32,-1/8) -- (-1/32,1/8);
				\end{scope}
			\end{scope}
			\begin{scope}[shift={(-3/8,3/2)}]
				\draw[thick] plot (-12/8,-12/32) -- (2/8,2/32);
				\fill [black]   (0,0) circle (0.05);
				\begin{scope}[shift={(-1/2,-1/8)}]
					\draw[thick] plot (1/32,-1/8) -- (-1/32,1/8);
				\end{scope}
				\begin{scope}[shift={(-2/2,-2/8)}]
					\draw[thick] plot (1/32,-1/8) -- (-1/32,1/8);
				\end{scope}
			\end{scope}
			\begin{scope}[shift={(-4/8,4/2)}]
				\draw[thick, blue] plot (-12/8,-12/32) -- (2/8,2/32);
				\fill [black]   (0,0) circle (0.05);
				\draw (-8/8,-8/32) node [below] {\hbox{\footnotesize$\geqone$}};
				
			\end{scope}
		\end{scope}
	\end{scope}
	
		\begin{scope}[xscale=1, yscale=-1, shift={(0,4.5)}]
		\draw[thick, orange] plot (-8/2,-8/8) -- (1/4,1/16);
		\fill [black]   (0,0) circle (0.05);
		\begin{scope}[shift={(-1/2,-1/8)}]
			\draw[thick, blue] plot (-3/8,3/2) -- (1/16,-1/4);
			\fill [black]   (0,0) circle (0.05);
			\draw (-3/8,3/2) node [below] {\hbox{\footnotesize$\geqone$}};
		\end{scope}
		\begin{scope}[shift={(-10/4,-10/16)}]
			\draw[thick, blue] plot (-3/8,3/2) -- (1/16,-1/4);
			\fill [black]   (0,0) circle (0.05);
			\draw (-3/8,3/2) node [below] {\hbox{\footnotesize$\geqone$}};
		\end{scope}
		\begin{scope}[shift={(-6/4,-6/16)}]
			\draw[thick, blue] plot (-3/8,3/2) -- (1/16,-1/4);
			\fill [black]   (0,0) circle (0.05);
			\draw (-3/8,3/2) node [below] {\hbox{\footnotesize$\geqone$}};
		\end{scope}
		\begin{scope}[shift={(-14/4,-14/16)}]
			\draw[thick] plot (-5/8+1/16,5/2-1/4) -- (1/16,-1/4);
			\fill [black]   (0,0) circle (0.05);
			\begin{scope}[shift={(-1/8,1/2)}]
				\draw[thick] plot (-12/8,-12/32) -- (2/8,2/32);
				\fill [black]   (0,0) circle (0.05);
				\begin{scope}[shift={(-1/2,-1/8)}]
					\draw[thick] plot (1/32,-1/8) -- (-1/32,1/8);
				\end{scope}
				\begin{scope}[shift={(-2/2,-2/8)}]
					\draw[thick] plot (1/32,-1/8) -- (-1/32,1/8);
				\end{scope}
			\end{scope}
			\begin{scope}[shift={(-2/8,2/2)}]
				\draw[thick] plot (-12/8,-12/32) -- (2/8,2/32);
				\fill [black]   (0,0) circle (0.05);
				\begin{scope}[shift={(-1/2,-1/8)}]
					\draw[thick] plot (1/32,-1/8) -- (-1/32,1/8);
				\end{scope}
				\begin{scope}[shift={(-2/2,-2/8)}]
					\draw[thick] plot (1/32,-1/8) -- (-1/32,1/8);
				\end{scope}
			\end{scope}
			\begin{scope}[shift={(-3/8,3/2)}]
				\draw[thick] plot (-12/8,-12/32) -- (2/8,2/32);
				\fill [black]   (0,0) circle (0.05);
				\begin{scope}[shift={(-1/2,-1/8)}]
					\draw[thick] plot (1/32,-1/8) -- (-1/32,1/8);
				\end{scope}
				\begin{scope}[shift={(-2/2,-2/8)}]
					\draw[thick] plot (1/32,-1/8) -- (-1/32,1/8);
				\end{scope}
			\end{scope}
			\begin{scope}[shift={(-4/8,4/2)}]
				\draw[thick, blue] plot (-12/8,-12/32) -- (2/8,2/32);
				\fill [black]   (0,0) circle (0.05);
				\draw (-8/8,-8/32) node [below] {\hbox{\footnotesize$\geqone$}};
				
			\end{scope}
		\end{scope}
	\end{scope}
	
		\begin{scope}[xscale=-1, yscale=-1, shift={(0,6)}]
		\draw[thick, orange] plot (-8/2,-8/8) -- (1/4,1/16);
		\fill [black]   (0,0) circle (0.05);
		\begin{scope}[shift={(-1/2,-1/8)}]
			\draw[thick, blue] plot (-3/8,3/2) -- (1/16,-1/4);
			\fill [black]   (0,0) circle (0.05);
			\draw (-3/8,3/2) node [below] {\hbox{\footnotesize$\geqone$}};
		\end{scope}
		\begin{scope}[shift={(-10/4,-10/16)}]
			\draw[thick, blue] plot (-3/8,3/2) -- (1/16,-1/4);
			\fill [black]   (0,0) circle (0.05);
			\draw (-3/8,3/2) node [below] {\hbox{\footnotesize$\geqone$}};
		\end{scope}
		\begin{scope}[shift={(-6/4,-6/16)}]
			\draw[thick, blue] plot (-3/8,3/2) -- (1/16,-1/4);
			\fill [black]   (0,0) circle (0.05);
			\draw (-3/8,3/2) node [below] {\hbox{\footnotesize$\geqone$}};
		\end{scope}
		\begin{scope}[shift={(-14/4,-14/16)}]
			\draw[thick] plot (-5/8+1/16,5/2-1/4) -- (1/16,-1/4);
			\fill [black]   (0,0) circle (0.05);
			\begin{scope}[shift={(-1/8,1/2)}]
				\draw[thick] plot (-12/8,-12/32) -- (2/8,2/32);
				\fill [black]   (0,0) circle (0.05);
				\begin{scope}[shift={(-1/2,-1/8)}]
					\draw[thick] plot (1/32,-1/8) -- (-1/32,1/8);
				\end{scope}
				\begin{scope}[shift={(-2/2,-2/8)}]
					\draw[thick] plot (1/32,-1/8) -- (-1/32,1/8);
				\end{scope}
			\end{scope}
			\begin{scope}[shift={(-2/8,2/2)}]
				\draw[thick] plot (-12/8,-12/32) -- (2/8,2/32);
				\fill [black]   (0,0) circle (0.05);
				\begin{scope}[shift={(-1/2,-1/8)}]
					\draw[thick] plot (1/32,-1/8) -- (-1/32,1/8);
				\end{scope}
				\begin{scope}[shift={(-2/2,-2/8)}]
					\draw[thick] plot (1/32,-1/8) -- (-1/32,1/8);
				\end{scope}
			\end{scope}
			\begin{scope}[shift={(-3/8,3/2)}]
				\draw[thick] plot (-12/8,-12/32) -- (2/8,2/32);
				\fill [black]   (0,0) circle (0.05);
				\begin{scope}[shift={(-1/2,-1/8)}]
					\draw[thick] plot (1/32,-1/8) -- (-1/32,1/8);
				\end{scope}
				\begin{scope}[shift={(-2/2,-2/8)}]
					\draw[thick] plot (1/32,-1/8) -- (-1/32,1/8);
				\end{scope}
			\end{scope}
			\begin{scope}[shift={(-4/8,4/2)}]
				\draw[thick, blue] plot (-12/8,-12/32) -- (2/8,2/32);
				\fill [black]   (0,0) circle (0.05);
				\draw (-8/8,-8/32) node [below] {\hbox{\footnotesize$\geqone$}};
				
			\end{scope}
		\end{scope}
	\end{scope}

	\end{tikzpicture}
}

\def\FigXseventeen{
	\begin{tikzpicture}[
		outer frame sep=0ex,%
		/tikz/inner frame xsep=2ex,
		/tikz/inner frame ysep=2ex,
		show background top,%
		show background bottom,%
		show background left,%
		show background right]
		\draw[thick] plot (0,1/4) -- (0,-6-1/2);
	\draw (-1,-6-2/2+1/4) node [below] {\hbox{\footnotesize$\geqnine$}};
	\begin{scope}[xscale=1, yscale=-1, shift={(0,6+1/4)}]
		\draw [thick, black] (0,1/4) arc (0:90:0.25);
		\draw [thick] plot (-1/4,1/2) --(-1-1/2,1/2);
		\end{scope}
		
		\begin{scope}[xscale=-1, yscale=-1, shift={(0,0)}]
			\draw[thick] plot (-5/2,-5/8) -- (1/4,1/16);
			\fill [black]   (0,0) circle (0.05);
			\begin{scope}[shift={(-1/2,-1/8)}]
				\draw[thick] plot (-3/8,3/2) -- (1/16,-1/4);
				\fill [black]   (0,0) circle (0.05);
				\begin{scope}[shift={(-1/8,1/2)}]
				\draw[thick] plot (-1/8,-1/32) -- (1/8,1/32);
				\end{scope}
					\begin{scope}[shift={(-2/8,2/2)}]
					\draw[thick] plot (-1/8,-1/32) -- (1/8,1/32);
				\end{scope}
			\end{scope}
			
			\begin{scope}[shift={(-2/2,-2/8)}]
				\draw[thick] plot (-3/8,3/2) -- (1/16,-1/4);
				\fill [black]   (0,0) circle (0.05);
				\begin{scope}[shift={(-1/8,1/2)}]
					\draw[thick] plot (-1/8,-1/32) -- (1/8,1/32);
				\end{scope}
				\begin{scope}[shift={(-2/8,2/2)}]
					\draw[thick] plot (-1/8,-1/32) -- (1/8,1/32);
				\end{scope}
			\end{scope}
			
			\begin{scope}[shift={(-3/2,-3/8)}]
				\draw[thick] plot (-3/8,3/2) -- (1/16,-1/4);
				\fill [black]   (0,0) circle (0.05);
				\begin{scope}[shift={(-1/8,1/2)}]
					\draw[thick] plot (-1/8,-1/32) -- (1/8,1/32);
				\end{scope}
				\begin{scope}[shift={(-2/8,2/2)}]
					\draw[thick] plot (-1/8,-1/32) -- (1/8,1/32);
				\end{scope}
			\end{scope}

			\begin{scope}[shift={(-9/4,-9/16)}]
				\draw[thick, orange] plot (-5/8+1/16,5/2-1/4) -- (1/16,-1/4);
				\fill [black]   (0,0) circle (0.05);
				\begin{scope}[shift={(-1/8,1/2)}]
					\draw[thick, blue] plot (-12/8,-12/32) -- (2/8,2/32);
					\fill [black]   (0,0) circle (0.05);
					\draw(-16/8,-12/32) node [] {\hbox{\footnotesize$\geqone$}};
				\end{scope}
				\begin{scope}[shift={(-2/8,2/2)}]
						\draw[thick, blue] plot (-12/8,-12/32) -- (2/8,2/32);
					\fill [black]   (0,0) circle (0.05);
					\draw(-16/8,-12/32) node [] {\hbox{\footnotesize$\geqone$}};
				\end{scope}
				\begin{scope}[shift={(-3/8,3/2)}]
					\draw[thick, blue] plot (-12/8,-12/32) -- (2/8,2/32);
				\fill [black]   (0,0) circle (0.05);
				\draw(-16/8,-12/32) node [] {\hbox{\footnotesize$\geqone$}};
				\end{scope}
				\begin{scope}[shift={(-4/8,4/2)}]
						\draw[thick] plot (-5/2,-5/8) -- (1/4,1/16);
					\fill [black]   (0,0) circle (0.05);
						\begin{scope}[shift={(-1/2,-1/8)}]
						\draw[thick] plot (1/32,-1/8) -- (-1/32,1/8);
					\end{scope}
						\begin{scope}[shift={(-2/2,-2/8)}]
						\draw[thick] plot (1/32,-1/8) -- (-1/32,1/8);
					\end{scope}
						\begin{scope}[shift={(-3/2,-3/8)}]
						\draw[thick] plot (1/32,-1/8) -- (-1/32,1/8);
					\end{scope}
						\begin{scope}[shift={(-4/2,-4/8)}]
						\draw[thick] plot (1/32,-1/8) -- (-1/32,1/8);
					\end{scope}
				\end{scope}
		\end{scope}
		\end{scope}

		\begin{scope}[xscale=1, yscale=-1, shift={(0,3/2)}]
			\draw[thick] plot (-5/2,-5/8) -- (1/4,1/16);
			\fill [black]   (0,0) circle (0.05);
			\begin{scope}[shift={(-1/2,-1/8)}]
				\draw[thick] plot (-3/8,3/2) -- (1/16,-1/4);
				\fill [black]   (0,0) circle (0.05);
				\begin{scope}[shift={(-1/8,1/2)}]
					\draw[thick] plot (-1/8,-1/32) -- (1/8,1/32);
				\end{scope}
				\begin{scope}[shift={(-2/8,2/2)}]
					\draw[thick] plot (-1/8,-1/32) -- (1/8,1/32);
				\end{scope}
			\end{scope}
			
			\begin{scope}[shift={(-2/2,-2/8)}]
				\draw[thick] plot (-3/8,3/2) -- (1/16,-1/4);
				\fill [black]   (0,0) circle (0.05);
				\begin{scope}[shift={(-1/8,1/2)}]
					\draw[thick] plot (-1/8,-1/32) -- (1/8,1/32);
				\end{scope}
				\begin{scope}[shift={(-2/8,2/2)}]
					\draw[thick] plot (-1/8,-1/32) -- (1/8,1/32);
				\end{scope}
			\end{scope}
			
			\begin{scope}[shift={(-3/2,-3/8)}]
				\draw[thick] plot (-3/8,3/2) -- (1/16,-1/4);
				\fill [black]   (0,0) circle (0.05);
				\begin{scope}[shift={(-1/8,1/2)}]
					\draw[thick] plot (-1/8,-1/32) -- (1/8,1/32);
				\end{scope}
				\begin{scope}[shift={(-2/8,2/2)}]
					\draw[thick] plot (-1/8,-1/32) -- (1/8,1/32);
				\end{scope}
			\end{scope}

			\begin{scope}[shift={(-9/4,-9/16)}]
				\draw[thick, orange] plot (-5/8+1/16,5/2-1/4) -- (1/16,-1/4);
				\fill [black]   (0,0) circle (0.05);
				\begin{scope}[shift={(-1/8,1/2)}]
					\draw[thick, blue] plot (-12/8,-12/32) -- (2/8,2/32);
					\fill [black]   (0,0) circle (0.05);
					\draw(-16/8,-12/32) node [] {\hbox{\footnotesize$\geqone$}};
				\end{scope}
				\begin{scope}[shift={(-2/8,2/2)}]
					\draw[thick, blue] plot (-12/8,-12/32) -- (2/8,2/32);
					\fill [black]   (0,0) circle (0.05);
					\draw(-16/8,-12/32) node [] {\hbox{\footnotesize$\geqone$}};
				\end{scope}
				\begin{scope}[shift={(-3/8,3/2)}]
					\draw[thick, blue] plot (-12/8,-12/32) -- (2/8,2/32);
					\fill [black]   (0,0) circle (0.05);
					\draw(-16/8,-12/32) node [] {\hbox{\footnotesize$\geqone$}};
				\end{scope}
				\begin{scope}[shift={(-4/8,4/2)}]
					\draw[thick] plot (-5/2,-5/8) -- (1/4,1/16);
					\fill [black]   (0,0) circle (0.05);
					\begin{scope}[shift={(-1/2,-1/8)}]
						\draw[thick] plot (1/32,-1/8) -- (-1/32,1/8);
					\end{scope}
					\begin{scope}[shift={(-2/2,-2/8)}]
						\draw[thick] plot (1/32,-1/8) -- (-1/32,1/8);
					\end{scope}
					\begin{scope}[shift={(-3/2,-3/8)}]
						\draw[thick] plot (1/32,-1/8) -- (-1/32,1/8);
					\end{scope}
					\begin{scope}[shift={(-4/2,-4/8)}]
						\draw[thick] plot (1/32,-1/8) -- (-1/32,1/8);
					\end{scope}
				\end{scope}
			\end{scope}
		\end{scope}

		\begin{scope}[xscale=-1, yscale=-1, shift={(0,3)}]
			\draw[thick] plot (-5/2,-5/8) -- (1/4,1/16);
			\fill [black]   (0,0) circle (0.05);
			\begin{scope}[shift={(-1/2,-1/8)}]
				\draw[thick] plot (-3/8,3/2) -- (1/16,-1/4);
				\fill [black]   (0,0) circle (0.05);
				\begin{scope}[shift={(-1/8,1/2)}]
					\draw[thick] plot (-1/8,-1/32) -- (1/8,1/32);
				\end{scope}
				\begin{scope}[shift={(-2/8,2/2)}]
					\draw[thick] plot (-1/8,-1/32) -- (1/8,1/32);
				\end{scope}
			\end{scope}
			
			\begin{scope}[shift={(-2/2,-2/8)}]
				\draw[thick] plot (-3/8,3/2) -- (1/16,-1/4);
				\fill [black]   (0,0) circle (0.05);
				\begin{scope}[shift={(-1/8,1/2)}]
					\draw[thick] plot (-1/8,-1/32) -- (1/8,1/32);
				\end{scope}
				\begin{scope}[shift={(-2/8,2/2)}]
					\draw[thick] plot (-1/8,-1/32) -- (1/8,1/32);
				\end{scope}
			\end{scope}
			
			\begin{scope}[shift={(-3/2,-3/8)}]
				\draw[thick] plot (-3/8,3/2) -- (1/16,-1/4);
				\fill [black]   (0,0) circle (0.05);
				\begin{scope}[shift={(-1/8,1/2)}]
					\draw[thick] plot (-1/8,-1/32) -- (1/8,1/32);
				\end{scope}
				\begin{scope}[shift={(-2/8,2/2)}]
					\draw[thick] plot (-1/8,-1/32) -- (1/8,1/32);
				\end{scope}
			\end{scope}

			\begin{scope}[shift={(-9/4,-9/16)}]
				\draw[thick, orange] plot (-5/8+1/16,5/2-1/4) -- (1/16,-1/4);
				\fill [black]   (0,0) circle (0.05);
				\begin{scope}[shift={(-1/8,1/2)}]
					\draw[thick, blue] plot (-12/8,-12/32) -- (2/8,2/32);
					\fill [black]   (0,0) circle (0.05);
					\draw(-16/8,-12/32) node [] {\hbox{\footnotesize$\geqone$}};
				\end{scope}
				\begin{scope}[shift={(-2/8,2/2)}]
					\draw[thick, blue] plot (-12/8,-12/32) -- (2/8,2/32);
					\fill [black]   (0,0) circle (0.05);
					\draw(-16/8,-12/32) node [] {\hbox{\footnotesize$\geqone$}};
				\end{scope}
				\begin{scope}[shift={(-3/8,3/2)}]
					\draw[thick, blue] plot (-12/8,-12/32) -- (2/8,2/32);
					\fill [black]   (0,0) circle (0.05);
					\draw(-16/8,-12/32) node [] {\hbox{\footnotesize$\geqone$}};
				\end{scope}
				\begin{scope}[shift={(-4/8,4/2)}]
					\draw[thick] plot (-5/2,-5/8) -- (1/4,1/16);
					\fill [black]   (0,0) circle (0.05);
					\begin{scope}[shift={(-1/2,-1/8)}]
						\draw[thick] plot (1/32,-1/8) -- (-1/32,1/8);
					\end{scope}
					\begin{scope}[shift={(-2/2,-2/8)}]
						\draw[thick] plot (1/32,-1/8) -- (-1/32,1/8);
					\end{scope}
					\begin{scope}[shift={(-3/2,-3/8)}]
						\draw[thick] plot (1/32,-1/8) -- (-1/32,1/8);
					\end{scope}
					\begin{scope}[shift={(-4/2,-4/8)}]
						\draw[thick] plot (1/32,-1/8) -- (-1/32,1/8);
					\end{scope}
				\end{scope}
			\end{scope}
		\end{scope}
		
		\begin{scope}[xscale=1, yscale=-1, shift={(0,4.5)}]
			\draw[thick] plot (-5/2,-5/8) -- (1/4,1/16);
			\fill [black]   (0,0) circle (0.05);
			\begin{scope}[shift={(-1/2,-1/8)}]
				\draw[thick] plot (-3/8,3/2) -- (1/16,-1/4);
				\fill [black]   (0,0) circle (0.05);
				\begin{scope}[shift={(-1/8,1/2)}]
					\draw[thick] plot (-1/8,-1/32) -- (1/8,1/32);
				\end{scope}
				\begin{scope}[shift={(-2/8,2/2)}]
					\draw[thick] plot (-1/8,-1/32) -- (1/8,1/32);
				\end{scope}
			\end{scope}
			
			\begin{scope}[shift={(-2/2,-2/8)}]
				\draw[thick] plot (-3/8,3/2) -- (1/16,-1/4);
				\fill [black]   (0,0) circle (0.05);
				\begin{scope}[shift={(-1/8,1/2)}]
					\draw[thick] plot (-1/8,-1/32) -- (1/8,1/32);
				\end{scope}
				\begin{scope}[shift={(-2/8,2/2)}]
					\draw[thick] plot (-1/8,-1/32) -- (1/8,1/32);
				\end{scope}
			\end{scope}
			
			\begin{scope}[shift={(-3/2,-3/8)}]
				\draw[thick] plot (-3/8,3/2) -- (1/16,-1/4);
				\fill [black]   (0,0) circle (0.05);
				\begin{scope}[shift={(-1/8,1/2)}]
					\draw[thick] plot (-1/8,-1/32) -- (1/8,1/32);
				\end{scope}
				\begin{scope}[shift={(-2/8,2/2)}]
					\draw[thick] plot (-1/8,-1/32) -- (1/8,1/32);
				\end{scope}
			\end{scope}

			\begin{scope}[shift={(-9/4,-9/16)}]
				\draw[thick, orange] plot (-5/8+1/16,5/2-1/4) -- (1/16,-1/4);
				\fill [black]   (0,0) circle (0.05);
				\begin{scope}[shift={(-1/8,1/2)}]
					\draw[thick, blue] plot (-12/8,-12/32) -- (2/8,2/32);
					\fill [black]   (0,0) circle (0.05);
					\draw(-16/8,-12/32) node [] {\hbox{\footnotesize$\geqone$}};
				\end{scope}
				\begin{scope}[shift={(-2/8,2/2)}]
					\draw[thick, blue] plot (-12/8,-12/32) -- (2/8,2/32);
					\fill [black]   (0,0) circle (0.05);
					\draw(-16/8,-12/32) node [] {\hbox{\footnotesize$\geqone$}};
				\end{scope}
				\begin{scope}[shift={(-3/8,3/2)}]
					\draw[thick, blue] plot (-12/8,-12/32) -- (2/8,2/32);
					\fill [black]   (0,0) circle (0.05);
					\draw(-16/8,-12/32) node [] {\hbox{\footnotesize$\geqone$}};
				\end{scope}
				\begin{scope}[shift={(-4/8,4/2)}]
					\draw[thick] plot (-5/2,-5/8) -- (1/4,1/16);
					\fill [black]   (0,0) circle (0.05);
					\begin{scope}[shift={(-1/2,-1/8)}]
						\draw[thick] plot (1/32,-1/8) -- (-1/32,1/8);
					\end{scope}
					\begin{scope}[shift={(-2/2,-2/8)}]
						\draw[thick] plot (1/32,-1/8) -- (-1/32,1/8);
					\end{scope}
					\begin{scope}[shift={(-3/2,-3/8)}]
						\draw[thick] plot (1/32,-1/8) -- (-1/32,1/8);
					\end{scope}
					\begin{scope}[shift={(-4/2,-4/8)}]
						\draw[thick] plot (1/32,-1/8) -- (-1/32,1/8);
					\end{scope}
				\end{scope}
			\end{scope}
		\end{scope}
		
		\begin{scope}[xscale=-1, yscale=-1, shift={(0,6)}]
			\draw[thick] plot (-5/2,-5/8) -- (1/4,1/16);
			\fill [black]   (0,0) circle (0.05);
			\begin{scope}[shift={(-1/2,-1/8)}]
				\draw[thick] plot (-3/8,3/2) -- (1/16,-1/4);
				\fill [black]   (0,0) circle (0.05);
				\begin{scope}[shift={(-1/8,1/2)}]
					\draw[thick] plot (-1/8,-1/32) -- (1/8,1/32);
				\end{scope}
				\begin{scope}[shift={(-2/8,2/2)}]
					\draw[thick] plot (-1/8,-1/32) -- (1/8,1/32);
				\end{scope}
			\end{scope}
			
			\begin{scope}[shift={(-2/2,-2/8)}]
				\draw[thick] plot (-3/8,3/2) -- (1/16,-1/4);
				\fill [black]   (0,0) circle (0.05);
				\begin{scope}[shift={(-1/8,1/2)}]
					\draw[thick] plot (-1/8,-1/32) -- (1/8,1/32);
				\end{scope}
				\begin{scope}[shift={(-2/8,2/2)}]
					\draw[thick] plot (-1/8,-1/32) -- (1/8,1/32);
				\end{scope}
			\end{scope}
			
			\begin{scope}[shift={(-3/2,-3/8)}]
				\draw[thick] plot (-3/8,3/2) -- (1/16,-1/4);
				\fill [black]   (0,0) circle (0.05);
				\begin{scope}[shift={(-1/8,1/2)}]
					\draw[thick] plot (-1/8,-1/32) -- (1/8,1/32);
				\end{scope}
				\begin{scope}[shift={(-2/8,2/2)}]
					\draw[thick] plot (-1/8,-1/32) -- (1/8,1/32);
				\end{scope}
			\end{scope}

			\begin{scope}[shift={(-9/4,-9/16)}]
				\draw[thick, orange] plot (-5/8+1/16,5/2-1/4) -- (1/16,-1/4);
				\fill [black]   (0,0) circle (0.05);
				\begin{scope}[shift={(-1/8,1/2)}]
					\draw[thick, blue] plot (-12/8,-12/32) -- (2/8,2/32);
					\fill [black]   (0,0) circle (0.05);
					\draw(-16/8,-12/32) node [] {\hbox{\footnotesize$\geqone$}};
				\end{scope}
				\begin{scope}[shift={(-2/8,2/2)}]
					\draw[thick, blue] plot (-12/8,-12/32) -- (2/8,2/32);
					\fill [black]   (0,0) circle (0.05);
					\draw(-16/8,-12/32) node [] {\hbox{\footnotesize$\geqone$}};
				\end{scope}
				\begin{scope}[shift={(-3/8,3/2)}]
					\draw[thick, blue] plot (-12/8,-12/32) -- (2/8,2/32);
					\fill [black]   (0,0) circle (0.05);
					\draw(-16/8,-12/32) node [] {\hbox{\footnotesize$\geqone$}};
				\end{scope}
				\begin{scope}[shift={(-4/8,4/2)}]
					\draw[thick] plot (-5/2,-5/8) -- (1/4,1/16);
					\fill [black]   (0,0) circle (0.05);
					\begin{scope}[shift={(-1/2,-1/8)}]
						\draw[thick] plot (1/32,-1/8) -- (-1/32,1/8);
					\end{scope}
					\begin{scope}[shift={(-2/2,-2/8)}]
						\draw[thick] plot (1/32,-1/8) -- (-1/32,1/8);
					\end{scope}
					\begin{scope}[shift={(-3/2,-3/8)}]
						\draw[thick] plot (1/32,-1/8) -- (-1/32,1/8);
					\end{scope}
					\begin{scope}[shift={(-4/2,-4/8)}]
						\draw[thick] plot (1/32,-1/8) -- (-1/32,1/8);
					\end{scope}
				\end{scope}
			\end{scope}
		\end{scope}

	\end{tikzpicture}
}

\def\FigXeightteen{
	\begin{tikzpicture}[
		outer frame sep=0ex,%
		/tikz/inner frame xsep=2ex,
		/tikz/inner frame ysep=2ex,
		show background top,%
		show background bottom,%
		show background left,%
		show background right]
		\draw[thick] plot (0,1/4) -- (0,-8-1/2);
		\draw [thick, black] (0,1/4) arc (0:90:0.25);
		\draw [thick] plot (-1/4,1/2) --(-3-1/2,1/2);
		\begin{scope}[shift={(-1/2,1/2)}]
			\draw[thick] (0,-1/8) -- (0,1/8);
		\end{scope}
		\begin{scope}[shift={(-2/2,1/2)}]
			\draw[thick] (0,-1/8) -- (0,1/8);
		\end{scope}
		\begin{scope}[shift={(-3/2,1/2)}]
			\draw[thick] (0,-1/8) -- (0,1/8);
		\end{scope}
		\begin{scope}[shift={(-4/2,1/2)}]
			\draw[thick] (0,-1/8) -- (0,1/8);
		\end{scope}
		\begin{scope}[shift={(-5/2,1/2)}]
			\draw[thick] (0,-1/8) -- (0,1/8);
		\end{scope}
		\begin{scope}[shift={(-6/2,1/2)}]
			\draw[thick] (0,-1/8) -- (0,1/8);
		\end{scope}
		
		\begin{scope}[xscale=1, yscale=-1, shift={(0,8+1/4)}]
			\draw [thick, black] (0,1/4) arc (0:90:0.25);
			\draw [thick] plot (-1/4,1/2) --(-3-1/2,1/2);
			\begin{scope}[shift={(-1/2,1/2)}]
				\draw[thick] (0,-1/8) -- (0,1/8);
			\end{scope}
			\begin{scope}[shift={(-2/2,1/2)}]
				\draw[thick] (0,-1/8) -- (0,1/8);
			\end{scope}
			\begin{scope}[shift={(-3/2,1/2)}]
				\draw[thick] (0,-1/8) -- (0,1/8);
			\end{scope}
			\begin{scope}[shift={(-4/2,1/2)}]
				\draw[thick] (0,-1/8) -- (0,1/8);
			\end{scope}
			\begin{scope}[shift={(-5/2,1/2)}]
				\draw[thick] (0,-1/8) -- (0,1/8);
			\end{scope}
			\begin{scope}[shift={(-6/2,1/2)}]
				\draw[thick] (0,-1/8) -- (0,1/8);
			\end{scope}
		\end{scope}

		\begin{scope}[xscale=-1, yscale=-1, shift={(0,0)}]
			\draw[thick, orange] plot (-8/2,-8/8) -- (1/4,1/16);
			\fill [black]   (0,0) circle (0.05);
			\begin{scope}[shift={(-1/2,-1/8)}]
				\draw[thick, blue] plot (-3/8,3/2) -- (1/16,-1/4);
				\fill [black]   (0,0) circle (0.05);
				\draw (-3/8,3/2) node [below] {\hbox{\footnotesize$\geqone$}};
			\end{scope}
			\begin{scope}[shift={(-10/4,-10/16)}]
				\draw[thick, blue] plot (-3/8,3/2) -- (1/16,-1/4);
				\fill [black]   (0,0) circle (0.05);
				\draw (-3/8,3/2) node [below] {\hbox{\footnotesize$\geqone$}};
			\end{scope}
			\begin{scope}[shift={(-6/4,-6/16)}]
				\draw[thick, blue] plot (-3/8,3/2) -- (1/16,-1/4);
				\fill [black]   (0,0) circle (0.05);
				\draw (-3/8,3/2) node [below] {\hbox{\footnotesize$\geqone$}};
			\end{scope}
			\begin{scope}[shift={(-14/4,-14/16)}]
				\draw[thick] plot (-5/8+1/16,5/2-1/4) -- (1/16,-1/4);
				\fill [black]   (0,0) circle (0.05);
				\begin{scope}[shift={(-1/8,1/2)}]
					\draw[thick] plot (-12/8,-12/32) -- (2/8,2/32);
					\fill [black]   (0,0) circle (0.05);
					\begin{scope}[shift={(-1/2,-1/8)}]
						\draw[thick] plot (1/32,-1/8) -- (-1/32,1/8);
					\end{scope}
					\begin{scope}[shift={(-2/2,-2/8)}]
						\draw[thick] plot (1/32,-1/8) -- (-1/32,1/8);
					\end{scope}
				\end{scope}
				\begin{scope}[shift={(-2/8,2/2)}]
					\draw[thick] plot (-12/8,-12/32) -- (2/8,2/32);
					\fill [black]   (0,0) circle (0.05);
					\begin{scope}[shift={(-1/2,-1/8)}]
						\draw[thick] plot (1/32,-1/8) -- (-1/32,1/8);
					\end{scope}
					\begin{scope}[shift={(-2/2,-2/8)}]
						\draw[thick] plot (1/32,-1/8) -- (-1/32,1/8);
					\end{scope}
				\end{scope}
				\begin{scope}[shift={(-3/8,3/2)}]
					\draw[thick] plot (-12/8,-12/32) -- (2/8,2/32);
					\fill [black]   (0,0) circle (0.05);
					\begin{scope}[shift={(-1/2,-1/8)}]
						\draw[thick] plot (1/32,-1/8) -- (-1/32,1/8);
					\end{scope}
					\begin{scope}[shift={(-2/2,-2/8)}]
						\draw[thick] plot (1/32,-1/8) -- (-1/32,1/8);
					\end{scope}
				\end{scope}
				\begin{scope}[shift={(-4/8,4/2)}]
				\draw[thick, orange] plot (-8/2,-8/8) -- (1/4,1/16);
				\fill [black]   (0,0) circle (0.05);
				\begin{scope}[shift={(-1/2,-1/8)}]
					\draw[thick, blue] plot (-3/8,3/2) -- (1/16,-1/4);
					\fill [black]   (0,0) circle (0.05);
					\draw (-3/8,3/2) node [below] {\hbox{\footnotesize$\geqone$}};
				\end{scope}
				\begin{scope}[shift={(-10/4,-10/16)}]
					\draw[thick, blue] plot (-3/8,3/2) -- (1/16,-1/4);
					\fill [black]   (0,0) circle (0.05);
					\draw (-3/8,3/2) node [below] {\hbox{\footnotesize$\geqone$}};
				\end{scope}
				\begin{scope}[shift={(-6/4,-6/16)}]
					\draw[thick, blue] plot (-3/8,3/2) -- (1/16,-1/4);
					\fill [black]   (0,0) circle (0.05);
					\draw (-3/8,3/2) node [below] {\hbox{\footnotesize$\geqone$}};
				\end{scope}
				\begin{scope}[shift={(-14/4,-14/16)}]
					\draw[thick] plot (-5/8,5/2) -- (1/16,-1/4);
					\fill [black]   (0,0) circle (0.05);
					\begin{scope}[shift={(-1/8,1/2)}]
						\draw[thick] plot (-1/8,-1/32) -- (1/8,1/32);
					\end{scope}
					\begin{scope}[shift={(-2/8,2/2)}]
						\draw[thick] plot (-1/8,-1/32) -- (1/8,1/32);
					\end{scope}
					\begin{scope}[shift={(-3/8,3/2)}]
						\draw[thick] plot (-1/8,-1/32) -- (1/8,1/32);
					\end{scope}
					\begin{scope}[shift={(-4/8,4/2)}]
						\draw[thick] plot (-1/8,-1/32) -- (1/8,1/32);
					\end{scope}
					
				\end{scope}
				\end{scope}
			\end{scope}
		\end{scope}

		\begin{scope}[xscale=1, yscale=-1, shift={(0,3/2)}]
		\draw[thick, orange] plot (-8/2,-8/8) -- (1/4,1/16);
		\fill [black]   (0,0) circle (0.05);
		\begin{scope}[shift={(-1/2,-1/8)}]
			\draw[thick, blue] plot (-3/8,3/2) -- (1/16,-1/4);
			\fill [black]   (0,0) circle (0.05);
			\draw (-3/8,3/2) node [below] {\hbox{\footnotesize$\geqone$}};
		\end{scope}
		\begin{scope}[shift={(-10/4,-10/16)}]
			\draw[thick, blue] plot (-3/8,3/2) -- (1/16,-1/4);
			\fill [black]   (0,0) circle (0.05);
			\draw (-3/8,3/2) node [below] {\hbox{\footnotesize$\geqone$}};
		\end{scope}
		\begin{scope}[shift={(-6/4,-6/16)}]
			\draw[thick, blue] plot (-3/8,3/2) -- (1/16,-1/4);
			\fill [black]   (0,0) circle (0.05);
			\draw (-3/8,3/2) node [below] {\hbox{\footnotesize$\geqone$}};
		\end{scope}
		\begin{scope}[shift={(-14/4,-14/16)}]
			\draw[thick] plot (-5/8+1/16,5/2-1/4) -- (1/16,-1/4);
			\fill [black]   (0,0) circle (0.05);
			\begin{scope}[shift={(-1/8,1/2)}]
				\draw[thick] plot (-12/8,-12/32) -- (2/8,2/32);
				\fill [black]   (0,0) circle (0.05);
				\begin{scope}[shift={(-1/2,-1/8)}]
					\draw[thick] plot (1/32,-1/8) -- (-1/32,1/8);
				\end{scope}
				\begin{scope}[shift={(-2/2,-2/8)}]
					\draw[thick] plot (1/32,-1/8) -- (-1/32,1/8);
				\end{scope}
			\end{scope}
			\begin{scope}[shift={(-2/8,2/2)}]
				\draw[thick] plot (-12/8,-12/32) -- (2/8,2/32);
				\fill [black]   (0,0) circle (0.05);
				\begin{scope}[shift={(-1/2,-1/8)}]
					\draw[thick] plot (1/32,-1/8) -- (-1/32,1/8);
				\end{scope}
				\begin{scope}[shift={(-2/2,-2/8)}]
					\draw[thick] plot (1/32,-1/8) -- (-1/32,1/8);
				\end{scope}
			\end{scope}
			\begin{scope}[shift={(-3/8,3/2)}]
				\draw[thick] plot (-12/8,-12/32) -- (2/8,2/32);
				\fill [black]   (0,0) circle (0.05);
				\begin{scope}[shift={(-1/2,-1/8)}]
					\draw[thick] plot (1/32,-1/8) -- (-1/32,1/8);
				\end{scope}
				\begin{scope}[shift={(-2/2,-2/8)}]
					\draw[thick] plot (1/32,-1/8) -- (-1/32,1/8);
				\end{scope}
			\end{scope}
			\begin{scope}[shift={(-4/8,4/2)}]
				\draw[thick, orange] plot (-8/2,-8/8) -- (1/4,1/16);
				\fill [black]   (0,0) circle (0.05);
				\begin{scope}[shift={(-1/2,-1/8)}]
					\draw[thick, blue] plot (-3/8,3/2) -- (1/16,-1/4);
					\fill [black]   (0,0) circle (0.05);
					\draw (-3/8,3/2) node [below] {\hbox{\footnotesize$\geqone$}};
				\end{scope}
				\begin{scope}[shift={(-10/4,-10/16)}]
					\draw[thick, blue] plot (-3/8,3/2) -- (1/16,-1/4);
					\fill [black]   (0,0) circle (0.05);
					\draw (-3/8,3/2) node [below] {\hbox{\footnotesize$\geqone$}};
				\end{scope}
				\begin{scope}[shift={(-6/4,-6/16)}]
					\draw[thick, blue] plot (-3/8,3/2) -- (1/16,-1/4);
					\fill [black]   (0,0) circle (0.05);
					\draw (-3/8,3/2) node [below] {\hbox{\footnotesize$\geqone$}};
				\end{scope}
				\begin{scope}[shift={(-14/4,-14/16)}]
					\draw[thick] plot (-5/8,5/2) -- (1/16,-1/4);
					\fill [black]   (0,0) circle (0.05);
					\begin{scope}[shift={(-1/8,1/2)}]
						\draw[thick] plot (-1/8,-1/32) -- (1/8,1/32);
					\end{scope}
					\begin{scope}[shift={(-2/8,2/2)}]
						\draw[thick] plot (-1/8,-1/32) -- (1/8,1/32);
					\end{scope}
					\begin{scope}[shift={(-3/8,3/2)}]
						\draw[thick] plot (-1/8,-1/32) -- (1/8,1/32);
					\end{scope}
					\begin{scope}[shift={(-4/8,4/2)}]
						\draw[thick] plot (-1/8,-1/32) -- (1/8,1/32);
					\end{scope}
					
				\end{scope}
			\end{scope}
		\end{scope}
		\end{scope}

		\begin{scope}[xscale=-1, yscale=-1, shift={(0,4)}]
		\draw[thick, orange] plot (-8/2,-8/8) -- (1/4,1/16);
		\fill [black]   (0,0) circle (0.05);
		\begin{scope}[shift={(-1/2,-1/8)}]
			\draw[thick, blue] plot (-3/8,3/2) -- (1/16,-1/4);
			\fill [black]   (0,0) circle (0.05);
			\draw (-3/8,3/2) node [below] {\hbox{\footnotesize$\geqone$}};
		\end{scope}
		\begin{scope}[shift={(-10/4,-10/16)}]
			\draw[thick, blue] plot (-3/8,3/2) -- (1/16,-1/4);
			\fill [black]   (0,0) circle (0.05);
			\draw (-3/8,3/2) node [below] {\hbox{\footnotesize$\geqone$}};
		\end{scope}
		\begin{scope}[shift={(-6/4,-6/16)}]
			\draw[thick, blue] plot (-3/8,3/2) -- (1/16,-1/4);
			\fill [black]   (0,0) circle (0.05);
			\draw (-3/8,3/2) node [below] {\hbox{\footnotesize$\geqone$}};
		\end{scope}
		\begin{scope}[shift={(-14/4,-14/16)}]
			\draw[thick] plot (-5/8+1/16,5/2-1/4) -- (1/16,-1/4);
			\fill [black]   (0,0) circle (0.05);
			\begin{scope}[shift={(-1/8,1/2)}]
				\draw[thick] plot (-12/8,-12/32) -- (2/8,2/32);
				\fill [black]   (0,0) circle (0.05);
				\begin{scope}[shift={(-1/2,-1/8)}]
					\draw[thick] plot (1/32,-1/8) -- (-1/32,1/8);
				\end{scope}
				\begin{scope}[shift={(-2/2,-2/8)}]
					\draw[thick] plot (1/32,-1/8) -- (-1/32,1/8);
				\end{scope}
			\end{scope}
			\begin{scope}[shift={(-2/8,2/2)}]
				\draw[thick] plot (-12/8,-12/32) -- (2/8,2/32);
				\fill [black]   (0,0) circle (0.05);
				\begin{scope}[shift={(-1/2,-1/8)}]
					\draw[thick] plot (1/32,-1/8) -- (-1/32,1/8);
				\end{scope}
				\begin{scope}[shift={(-2/2,-2/8)}]
					\draw[thick] plot (1/32,-1/8) -- (-1/32,1/8);
				\end{scope}
			\end{scope}
			\begin{scope}[shift={(-3/8,3/2)}]
				\draw[thick] plot (-12/8,-12/32) -- (2/8,2/32);
				\fill [black]   (0,0) circle (0.05);
				\begin{scope}[shift={(-1/2,-1/8)}]
					\draw[thick] plot (1/32,-1/8) -- (-1/32,1/8);
				\end{scope}
				\begin{scope}[shift={(-2/2,-2/8)}]
					\draw[thick] plot (1/32,-1/8) -- (-1/32,1/8);
				\end{scope}
			\end{scope}
			\begin{scope}[shift={(-4/8,4/2)}]
				\draw[thick, orange] plot (-8/2,-8/8) -- (1/4,1/16);
				\fill [black]   (0,0) circle (0.05);
				\begin{scope}[shift={(-1/2,-1/8)}]
					\draw[thick, blue] plot (-3/8,3/2) -- (1/16,-1/4);
					\fill [black]   (0,0) circle (0.05);
					\draw (-3/8,3/2) node [below] {\hbox{\footnotesize$\geqone$}};
				\end{scope}
				\begin{scope}[shift={(-10/4,-10/16)}]
					\draw[thick, blue] plot (-3/8,3/2) -- (1/16,-1/4);
					\fill [black]   (0,0) circle (0.05);
					\draw (-3/8,3/2) node [below] {\hbox{\footnotesize$\geqone$}};
				\end{scope}
				\begin{scope}[shift={(-6/4,-6/16)}]
					\draw[thick, blue] plot (-3/8,3/2) -- (1/16,-1/4);
					\fill [black]   (0,0) circle (0.05);
					\draw (-3/8,3/2) node [below] {\hbox{\footnotesize$\geqone$}};
				\end{scope}
				\begin{scope}[shift={(-14/4,-14/16)}]
					\draw[thick] plot (-5/8,5/2) -- (1/16,-1/4);
					\fill [black]   (0,0) circle (0.05);
					\begin{scope}[shift={(-1/8,1/2)}]
						\draw[thick] plot (-1/8,-1/32) -- (1/8,1/32);
					\end{scope}
					\begin{scope}[shift={(-2/8,2/2)}]
						\draw[thick] plot (-1/8,-1/32) -- (1/8,1/32);
					\end{scope}
					\begin{scope}[shift={(-3/8,3/2)}]
						\draw[thick] plot (-1/8,-1/32) -- (1/8,1/32);
					\end{scope}
					\begin{scope}[shift={(-4/8,4/2)}]
						\draw[thick] plot (-1/8,-1/32) -- (1/8,1/32);
					\end{scope}
					
				\end{scope}
			\end{scope}
		\end{scope}
		\end{scope}
		
		\begin{scope}[xscale=1, yscale=-1, shift={(0,5.5)}]
		\draw[thick, orange] plot (-8/2,-8/8) -- (1/4,1/16);
		\fill [black]   (0,0) circle (0.05);
		\begin{scope}[shift={(-1/2,-1/8)}]
			\draw[thick, blue] plot (-3/8,3/2) -- (1/16,-1/4);
			\fill [black]   (0,0) circle (0.05);
			\draw (-3/8,3/2) node [below] {\hbox{\footnotesize$\geqone$}};
		\end{scope}
		\begin{scope}[shift={(-10/4,-10/16)}]
			\draw[thick, blue] plot (-3/8,3/2) -- (1/16,-1/4);
			\fill [black]   (0,0) circle (0.05);
			\draw (-3/8,3/2) node [below] {\hbox{\footnotesize$\geqone$}};
		\end{scope}
		\begin{scope}[shift={(-6/4,-6/16)}]
			\draw[thick, blue] plot (-3/8,3/2) -- (1/16,-1/4);
			\fill [black]   (0,0) circle (0.05);
			\draw (-3/8,3/2) node [below] {\hbox{\footnotesize$\geqone$}};
		\end{scope}
		\begin{scope}[shift={(-14/4,-14/16)}]
			\draw[thick] plot (-5/8+1/16,5/2-1/4) -- (1/16,-1/4);
			\fill [black]   (0,0) circle (0.05);
			\begin{scope}[shift={(-1/8,1/2)}]
				\draw[thick] plot (-12/8,-12/32) -- (2/8,2/32);
				\fill [black]   (0,0) circle (0.05);
				\begin{scope}[shift={(-1/2,-1/8)}]
					\draw[thick] plot (1/32,-1/8) -- (-1/32,1/8);
				\end{scope}
				\begin{scope}[shift={(-2/2,-2/8)}]
					\draw[thick] plot (1/32,-1/8) -- (-1/32,1/8);
				\end{scope}
			\end{scope}
			\begin{scope}[shift={(-2/8,2/2)}]
				\draw[thick] plot (-12/8,-12/32) -- (2/8,2/32);
				\fill [black]   (0,0) circle (0.05);
				\begin{scope}[shift={(-1/2,-1/8)}]
					\draw[thick] plot (1/32,-1/8) -- (-1/32,1/8);
				\end{scope}
				\begin{scope}[shift={(-2/2,-2/8)}]
					\draw[thick] plot (1/32,-1/8) -- (-1/32,1/8);
				\end{scope}
			\end{scope}
			\begin{scope}[shift={(-3/8,3/2)}]
				\draw[thick] plot (-12/8,-12/32) -- (2/8,2/32);
				\fill [black]   (0,0) circle (0.05);
				\begin{scope}[shift={(-1/2,-1/8)}]
					\draw[thick] plot (1/32,-1/8) -- (-1/32,1/8);
				\end{scope}
				\begin{scope}[shift={(-2/2,-2/8)}]
					\draw[thick] plot (1/32,-1/8) -- (-1/32,1/8);
				\end{scope}
			\end{scope}
			\begin{scope}[shift={(-4/8,4/2)}]
				\draw[thick, orange] plot (-8/2,-8/8) -- (1/4,1/16);
				\fill [black]   (0,0) circle (0.05);
				\begin{scope}[shift={(-1/2,-1/8)}]
					\draw[thick, blue] plot (-3/8,3/2) -- (1/16,-1/4);
					\fill [black]   (0,0) circle (0.05);
					\draw (-3/8,3/2) node [below] {\hbox{\footnotesize$\geqone$}};
				\end{scope}
				\begin{scope}[shift={(-10/4,-10/16)}]
					\draw[thick, blue] plot (-3/8,3/2) -- (1/16,-1/4);
					\fill [black]   (0,0) circle (0.05);
					\draw (-3/8,3/2) node [below] {\hbox{\footnotesize$\geqone$}};
				\end{scope}
				\begin{scope}[shift={(-6/4,-6/16)}]
					\draw[thick, blue] plot (-3/8,3/2) -- (1/16,-1/4);
					\fill [black]   (0,0) circle (0.05);
					\draw (-3/8,3/2) node [below] {\hbox{\footnotesize$\geqone$}};
				\end{scope}
				\begin{scope}[shift={(-14/4,-14/16)}]
					\draw[thick] plot (-5/8,5/2) -- (1/16,-1/4);
					\fill [black]   (0,0) circle (0.05);
					\begin{scope}[shift={(-1/8,1/2)}]
						\draw[thick] plot (-1/8,-1/32) -- (1/8,1/32);
					\end{scope}
					\begin{scope}[shift={(-2/8,2/2)}]
						\draw[thick] plot (-1/8,-1/32) -- (1/8,1/32);
					\end{scope}
					\begin{scope}[shift={(-3/8,3/2)}]
						\draw[thick] plot (-1/8,-1/32) -- (1/8,1/32);
					\end{scope}
					\begin{scope}[shift={(-4/8,4/2)}]
						\draw[thick] plot (-1/8,-1/32) -- (1/8,1/32);
					\end{scope}
					
				\end{scope}
			\end{scope}
		\end{scope}
		\end{scope}
		
		\begin{scope}[xscale=-1, yscale=-1, shift={(0,8)}]
		\draw[thick, orange] plot (-8/2,-8/8) -- (1/4,1/16);
		\fill [black]   (0,0) circle (0.05);
		\begin{scope}[shift={(-1/2,-1/8)}]
			\draw[thick, blue] plot (-3/8,3/2) -- (1/16,-1/4);
			\fill [black]   (0,0) circle (0.05);
			\draw (-3/8,3/2) node [below] {\hbox{\footnotesize$\geqone$}};
		\end{scope}
		\begin{scope}[shift={(-10/4,-10/16)}]
			\draw[thick, blue] plot (-3/8,3/2) -- (1/16,-1/4);
			\fill [black]   (0,0) circle (0.05);
			\draw (-3/8,3/2) node [below] {\hbox{\footnotesize$\geqone$}};
		\end{scope}
		\begin{scope}[shift={(-6/4,-6/16)}]
			\draw[thick, blue] plot (-3/8,3/2) -- (1/16,-1/4);
			\fill [black]   (0,0) circle (0.05);
			\draw (-3/8,3/2) node [below] {\hbox{\footnotesize$\geqone$}};
		\end{scope}
		\begin{scope}[shift={(-14/4,-14/16)}]
			\draw[thick] plot (-5/8+1/16,5/2-1/4) -- (1/16,-1/4);
			\fill [black]   (0,0) circle (0.05);
			\begin{scope}[shift={(-1/8,1/2)}]
				\draw[thick] plot (-12/8,-12/32) -- (2/8,2/32);
				\fill [black]   (0,0) circle (0.05);
				\begin{scope}[shift={(-1/2,-1/8)}]
					\draw[thick] plot (1/32,-1/8) -- (-1/32,1/8);
				\end{scope}
				\begin{scope}[shift={(-2/2,-2/8)}]
					\draw[thick] plot (1/32,-1/8) -- (-1/32,1/8);
				\end{scope}
			\end{scope}
			\begin{scope}[shift={(-2/8,2/2)}]
				\draw[thick] plot (-12/8,-12/32) -- (2/8,2/32);
				\fill [black]   (0,0) circle (0.05);
				\begin{scope}[shift={(-1/2,-1/8)}]
					\draw[thick] plot (1/32,-1/8) -- (-1/32,1/8);
				\end{scope}
				\begin{scope}[shift={(-2/2,-2/8)}]
					\draw[thick] plot (1/32,-1/8) -- (-1/32,1/8);
				\end{scope}
			\end{scope}
			\begin{scope}[shift={(-3/8,3/2)}]
				\draw[thick] plot (-12/8,-12/32) -- (2/8,2/32);
				\fill [black]   (0,0) circle (0.05);
				\begin{scope}[shift={(-1/2,-1/8)}]
					\draw[thick] plot (1/32,-1/8) -- (-1/32,1/8);
				\end{scope}
				\begin{scope}[shift={(-2/2,-2/8)}]
					\draw[thick] plot (1/32,-1/8) -- (-1/32,1/8);
				\end{scope}
			\end{scope}
			\begin{scope}[shift={(-4/8,4/2)}]
				\draw[thick, orange] plot (-8/2,-8/8) -- (1/4,1/16);
				\fill [black]   (0,0) circle (0.05);
				\begin{scope}[shift={(-1/2,-1/8)}]
					\draw[thick, blue] plot (-3/8,3/2) -- (1/16,-1/4);
					\fill [black]   (0,0) circle (0.05);
					\draw (-3/8,3/2) node [below] {\hbox{\footnotesize$\geqone$}};
				\end{scope}
				\begin{scope}[shift={(-10/4,-10/16)}]
					\draw[thick, blue] plot (-3/8,3/2) -- (1/16,-1/4);
					\fill [black]   (0,0) circle (0.05);
					\draw (-3/8,3/2) node [below] {\hbox{\footnotesize$\geqone$}};
				\end{scope}
				\begin{scope}[shift={(-6/4,-6/16)}]
					\draw[thick, blue] plot (-3/8,3/2) -- (1/16,-1/4);
					\fill [black]   (0,0) circle (0.05);
					\draw (-3/8,3/2) node [below] {\hbox{\footnotesize$\geqone$}};
				\end{scope}
				\begin{scope}[shift={(-14/4,-14/16)}]
					\draw[thick] plot (-5/8,5/2) -- (1/16,-1/4);
					\fill [black]   (0,0) circle (0.05);
					\begin{scope}[shift={(-1/8,1/2)}]
						\draw[thick] plot (-1/8,-1/32) -- (1/8,1/32);
					\end{scope}
					\begin{scope}[shift={(-2/8,2/2)}]
						\draw[thick] plot (-1/8,-1/32) -- (1/8,1/32);
					\end{scope}
					\begin{scope}[shift={(-3/8,3/2)}]
						\draw[thick] plot (-1/8,-1/32) -- (1/8,1/32);
					\end{scope}
					\begin{scope}[shift={(-4/8,4/2)}]
						\draw[thick] plot (-1/8,-1/32) -- (1/8,1/32);
					\end{scope}
					
				\end{scope}
			\end{scope}
		\end{scope}
		\end{scope}

	\end{tikzpicture}
}

\def\FigGenusThreeExampleOneGroundedDPdownstairs{
		\begin{tikzpicture}[
		outer frame sep=0ex,%
		/tikz/inner frame xsep=2ex,
		/tikz/inner frame ysep=2ex,
		show background top,%
		show background bottom,%
		show background left,%
		show background right]
		\draw[thick] plot (0,1/4+0/4) -- (0,-3/2-1/4);

		\begin{scope}[xscale=1, shift={(0,0)}]
			\draw[thick,black] plot (-1/4,-1/16) -- (2+4/4,1/2+4/16);
				markings 
							\begin{scope}[shift={(3/2,3/8)}]
									\draw[thick] plot (-1/32,1/8) -- (1/32,-1/8);
								\end{scope}
								\begin{scope}[shift={(4/2,4/8)}]
								\draw[thick] plot (-1/32,1/8) -- (1/32,-1/8);
							\end{scope}
								\begin{scope}[shift={(5/2,5/8)}]
								\draw[thick] plot (-1/32,1/8) -- (1/32,-1/8);
							\end{scope}
			\fill [black]   (0,0) circle (0.05);%
			\begin{scope}[xscale=-1, shift={(-1,1/4)}]
				\draw[thick,blue] plot (-1/16,-1/4) -- (1/4+1/8,1+1/2);
				\fill [black]   (0,0) circle (0.05);%
				\draw (3/8-1/8,3/2-1/4) node [right] {\hbox{\footnotesize$\geqone$}};
			\end{scope}

		\end{scope}
		
		\begin{scope}[xscale=-1, shift={(0,-1/2-1/4)}]
			\draw[thick,orange] plot (-1/4,-1/16) -- (2+1/4,1/2+1/16);
			\fill [black]   (0,0) circle (0.05);%
			\begin{scope}[xscale=1,yscale=-1, shift={(1,-1/4)}]
				\draw[thick,blue] plot (-1/16,-1/4) -- (1/4+1/8,1+1/2);
				\fill [black]   (0,0) circle (0.05);%
				\draw (3/8-0/8,3/2-1/4) node [left] {\hbox{\footnotesize$\geqone$}};
			\end{scope}

			component of type (a) with two marked points
			\begin{scope}[xscale=-1,shift={(-2,1/2)}]
				\fill [black]   (0,0) circle (0.05);
				\draw[thick] plot (-2*3/4,-1/2*3/4) -- (1/4,1/16);
				\begin{scope}[shift={(-1/2,-1/8)}]
					\draw[thick] plot (-1/32,1/8) -- (1/32,-1/8);
				\end{scope}
				\begin{scope}[shift={(-1,-1/4)}]
					\draw[thick] plot (-1/32,1/8) -- (1/32,-1/8);
				\end{scope}
			\end{scope}

		\end{scope}
		
		\begin{scope}[xscale=1,yscale=1, shift={(0,-1-1/2)}]
		\draw[thick,black] plot (-1/4,-1/16) -- (2+4/4,1/2+4/16);
		markings 
		\begin{scope}[shift={(3/2,3/8)}]
			\draw[thick] plot (-1/32,1/8) -- (1/32,-1/8);
		\end{scope}
		\begin{scope}[shift={(4/2,4/8)}]
			\draw[thick] plot (-1/32,1/8) -- (1/32,-1/8);
		\end{scope}
		\begin{scope}[shift={(5/2,5/8)}]
			\draw[thick] plot (-1/32,1/8) -- (1/32,-1/8);
		\end{scope}
		\fill [black]   (0,0) circle (0.05);%
		\begin{scope}[yscale=-1, shift={(1,-1/4)}]
			\draw[thick,blue] plot (-1/16,-1/4) -- (1/4+1/8,1+1/2);
			\fill [black]   (0,0) circle (0.05);%
			\draw (3/8-1/8,3/2-1/4) node [left] {\hbox{\footnotesize$\geqone$}};
		\end{scope}

		\end{scope}

	\end{tikzpicture}
}

\def\FigGenusThreeExampleWithThreeGroundedDoublePoints{
	\begin{tikzpicture}[
		outer frame sep=0ex,%
		/tikz/inner frame xsep=2ex,
		/tikz/inner frame ysep=2ex,
		show background top,%
		show background bottom,%
		show background left,%
		show background right]
		\draw[thick] plot (0,1/4+4/4) -- (0,-3/2-1/4);
		\draw[thick] plot (3/16,2/4) -- (-3/16,2/4);
			\draw[thick] plot (3/16,4/4) -- (-3/16,4/4);
		
		\begin{scope}[xscale=1, shift={(0,0)}]
			\draw[thick,orange] plot (-1/4,-1/16) -- (2+1/4,1/2+1/16);
			\fill [black]   (0,0) circle (0.05);%
			\begin{scope}[xscale=-1, shift={(-1,1/4)}]
				\draw[thick,blue] plot (-1/16,-1/4) -- (1/4+1/8,1+1/2);
				\fill [black]   (0,0) circle (0.05);%
				\draw (3/8-1/8,3/2-1/4) node [right] {\hbox{\footnotesize$\geqone$}};
			\end{scope}

			component of type (a) with two marked points
				\begin{scope}[xscale=-1,shift={(-2,1/2)}]
				\fill [black]   (0,0) circle (0.05);
				\draw[thick] plot (-2*3/4,-1/2*3/4) -- (1/4,1/16);
				\begin{scope}[shift={(-1/2,-1/8)}]
					\draw[thick] plot (-1/32,1/8) -- (1/32,-1/8);
				\end{scope}
				\begin{scope}[shift={(-1,-1/4)}]
					\draw[thick] plot (-1/32,1/8) -- (1/32,-1/8);
				\end{scope}
			\end{scope}

		\end{scope}
		
		\begin{scope}[xscale=-1, shift={(0,-1/2-1/4)}]
		\draw[thick,orange] plot (-1/4,-1/16) -- (2+1/4,1/2+1/16);
		\fill [black]   (0,0) circle (0.05);%
		\begin{scope}[xscale=1,yscale=-1, shift={(1,-1/4)}]
			\draw[thick,blue] plot (-1/16,-1/4) -- (1/4+1/8,1+1/2);
			\fill [black]   (0,0) circle (0.05);%
			\draw (3/8-0/8,3/2-1/4) node [left] {\hbox{\footnotesize$\geqone$}};
		\end{scope}

		component of type (a) with two marked points
		\begin{scope}[xscale=-1,shift={(-2,1/2)}]
			\fill [black]   (0,0) circle (0.05);
			\draw[thick] plot (-2*3/4,-1/2*3/4) -- (1/4,1/16);
			\begin{scope}[shift={(-1/2,-1/8)}]
				\draw[thick] plot (-1/32,1/8) -- (1/32,-1/8);
			\end{scope}
			\begin{scope}[shift={(-1,-1/4)}]
				\draw[thick] plot (-1/32,1/8) -- (1/32,-1/8);
			\end{scope}
		\end{scope}

		\end{scope}
		
		\begin{scope}[xscale=1,yscale=1, shift={(0,-1-1/2)}]
		\draw[thick,orange] plot (-1/4,-1/16) -- (2+1/4,1/2+1/16);
		\fill [black]   (0,0) circle (0.05);%
		\begin{scope}[xscale=1,yscale=-1, shift={(1,-1/4)}]
			\draw[thick,blue] plot (-1/16,-1/4) -- (1/4+1/8,1+1/2);
			\fill [black]   (0,0) circle (0.05);%
			\draw (3/8-1/8,3/2-1/4) node [left] {\hbox{\footnotesize$\geqone$}};
		\end{scope}

		component of type (a) with two marked points
		\begin{scope}[xscale=-1,shift={(-2,1/2)}]
			\fill [black]   (0,0) circle (0.05);
			\draw[thick] plot (-2*3/4,-1/2*3/4) -- (1/4,1/16);
			\begin{scope}[shift={(-1/2,-1/8)}]
				\draw[thick] plot (-1/32,1/8) -- (1/32,-1/8);
			\end{scope}
			\begin{scope}[shift={(-1,-1/4)}]
				\draw[thick] plot (-1/32,1/8) -- (1/32,-1/8);
			\end{scope}
		\end{scope}

		\end{scope}

	\end{tikzpicture}
}

\def\FigGenusThreeExampleNotingOverDP{
	\begin{tikzpicture}[
		outer frame sep=0ex,%
		/tikz/inner frame xsep=2ex,
		/tikz/inner frame ysep=2ex,
		show background top,%
		show background bottom,%
		show background left,%
		show background right]
		\draw[thick] plot (0,1/4+9/4) -- (0,-3/2-1/4);
		\begin{scope}[yscale=-1, shift={(0,0)}]
			\fill [black]   (0,0) circle (0.05);
			\draw[thick] plot (-2*3/4,-1/2*3/4) -- (1/4,1/16);
			\begin{scope}[shift={(-1/2,-1/8)}]
				\draw[thick] plot (-1/32,1/8) -- (1/32,-1/8);
			\end{scope}
			\begin{scope}[shift={(-1,-1/4)}]
				\draw[thick] plot (-1/32,1/8) -- (1/32,-1/8);
			\end{scope}

		\end{scope}
		
				\begin{scope}[yscale=-1, shift={(0,3/2)}]
			\fill [black]   (0,0) circle (0.05);
			\draw[thick] plot (-2*3/4,-1/2*3/4) -- (1/4,1/16);
			\begin{scope}[shift={(-1/2,-1/8)}]
				\draw[thick] plot (-1/32,1/8) -- (1/32,-1/8);
			\end{scope}
			\begin{scope}[shift={(-1,-1/4)}]
				\draw[thick] plot (-1/32,1/8) -- (1/32,-1/8);
			\end{scope}

		\end{scope}
		
		\begin{scope}[xscale=1, shift={(0,-1/2-1/4)}]
			\draw[thick,black] plot (-1/4,-1/16) -- (2+2/4,1/2+2/16);
			\fill [black]   (0,0) circle (0.05);%
					\begin{scope}[shift={(1/2,1/8)}]
						\draw[thick] plot (-1/32,1/8) -- (1/32,-1/8);
					\end{scope}
					\begin{scope}[shift={(1,1/4)}]
						\draw[thick] plot (-1/32,1/8) -- (1/32,-1/8);
					\end{scope}
				\begin{scope}[shift={(3/2,3/8)}]
				\draw[thick] plot (-1/32,1/8) -- (1/32,-1/8);
			\end{scope}
				\begin{scope}[shift={(4/2,4/8)}]
				\draw[thick] plot (-1/32,1/8) -- (1/32,-1/8);
			\end{scope}
		\end{scope}
		
			\begin{scope}[xscale=1, shift={(0,3/4)}]
				\draw[thick,blue] plot (-1/4,-1/16) -- (2+1/4,1/2+1/16);
				\fill [black]   (0,0) circle (0.05);%
				\draw (2+0/4,1/2+0/16) node [below] {\hbox{\footnotesize$\geqone$}};
			\end{scope}
			
						\begin{scope}[xscale=-1, shift={(0,6/4)}]
				\draw[thick,blue] plot (-1/4,-1/16) -- (2+1/4,1/2+1/16);
				\fill [black]   (0,0) circle (0.05);%
				\draw (2+0/4,1/2+0/16) node [below] {\hbox{\footnotesize$\geqone$}};
			\end{scope}
			
					\begin{scope}[xscale=1, shift={(0,9/4)}]
				\draw[thick,blue] plot (-1/4,-1/16) -- (2+1/4,1/2+1/16);
				\fill [black]   (0,0) circle (0.05);%
				\draw (2+0/4,1/2+0/16) node [below] {\hbox{\footnotesize$\geqone$}};
			\end{scope}

	\end{tikzpicture}
}

\def\FigMotivatingExampleFigure{
		\begin{tikzpicture}[scale=0.7,
		outer frame sep=0ex,%
		/tikz/inner frame xsep=2ex,
		/tikz/inner frame ysep=2ex,
		show background top,%
		show background bottom,%
		show background left,%
		show background right]
		\draw [thick,black] plot  coordinates {(4/8+4/8,-1/8-4/32) (-20/8,5/8) };
		\begin{scope}[shift={(4/8,-1/8)}]
			\draw [thick,black] plot  coordinates {(-1/32,-1/8) (1/32,1/8) };
		\end{scope}
		\begin{scope}[yscale=1, xscale=1, shift={(0,0)}]
			\draw [thick,black] plot  coordinates {(5/24,20/24) (-3/16,-12/16) };
			\fill [white] (0,0) circle (0.1);
			\draw [thick] (0,0) circle (0.1);
			\begin{scope}[shift={(-3/32,-12/32)}]
				\draw [thick,black] plot  coordinates {(-1/8,1/32) (1/8,-1/32) };
			\end{scope}
			\begin{scope}[xscale=-1, yscale=-1, shift={(-4/32,-16/32)}]
				\draw [thick,black] plot  coordinates {(-12/8,12/32) (4/8,-4/32) };
				\fill [red] (0,0) circle (0.10);
				
				\begin{scope}[shift={(-4/8,1/8)}]
					\draw [thick,black] plot  coordinates {(-1/32,-1/8) (1/32,1/8) };
				\end{scope}
				\begin{scope}[shift={(-8/8,8/32)}]
					\draw [thick,black] plot  coordinates {(-1/32,-1/8) (1/32,1/8) };
				\end{scope}
				
			\end{scope}
			
		\end{scope}
		\begin{scope}[shift={(-4/4,1/4)}]
			\draw [thick,black] plot  coordinates {(2/24,8/24) (-5/16,-20/16) };
			\fill [red] (0,0) circle (0.10);
			\begin{scope}[shift={(-3/32,-12/32)}]
				\draw [thick,black] plot  coordinates {(-1/8,1/32) (1/8,-1/32) };
			\end{scope}
			\begin{scope}[shift={(-6/32,-24/32)}]
				\draw [thick,black] plot  coordinates {(-1/8,1/32) (1/8,-1/32) };
			\end{scope}
		\end{scope}
		\begin{scope}[shift={(-8/4,2/4)}]
			\draw [thick,black] plot  coordinates {(2/24,8/24) (-5/16,-20/16) };
			\fill [red] (0,0) circle (0.10);
			\begin{scope}[shift={(-3/32,-12/32)}]
				\draw [thick,black] plot  coordinates {(-1/8,1/32) (1/8,-1/32) };
			\end{scope}
			\begin{scope}[shift={(-6/32,-24/32)}]
				\draw [thick,black] plot  coordinates {(-1/8,1/32) (1/8,-1/32) };
			\end{scope}
		\end{scope}
\end{tikzpicture}}
\def\bigtimes{\mathop{\raise-2pt\hbox{\huge$\times$}}}

\newbox\circbulletbox
\setbox\circbulletbox=\hbox{\,{\large$\circledcirc$}\kern-8.7pt\raise .6pt\hbox{$\bullet$}\,}

\let\le\leqslant
\let\ge\geqslant
\let\leq\leqslant
\let\geq\geqslant

\def\circVbig{\hbox{\text{\it\r{V}}}}
\def\circVscript{\hbox{\scriptsize\text{\it\r{V}}}}
\def\circVscriptscript{\mbox{\tiny\text{\it\r{V}}}}

\def\circVlimits_#1^#2{{\mathchoice%
{\circVbig{}^{\kern2pt #2}_{\kern-2pt #1}}%
{\circVbig{}^{\kern2pt #2}_{\kern-2pt #1}}%
{\scriptstyle\circVscript{}^{\kern1.7pt #2}_{\kern-1pt #1}}%
{\scriptscriptstyle\circVscriptscript{}^{\kern1.5pt #2}_{\kern-1pt #1}}%
}}
\def\circVr_#1{\circVlimits_#1^r}
\def\circVs_#1{\circVlimits_#1^s}

\def\circWbig{\hbox{\text{\it\r{W}}}}
\def\circWscript{\hbox{\scriptsize\text{\it\r{W}}}}
\def\circWscriptscript{\mbox{\tiny\text{\it\r{W}}}}

\def\circWlimits_#1^#2{{\mathchoice%
{\circWbig{}^{\kern2pt #2}_{\kern-2pt #1}}%
{\circWbig{}^{\kern2pt #2}_{\kern-2pt #1}}%
{\scriptstyle\circWscript{}^{\kern1.7pt #2}_{\kern-1pt #1}}%
{\scriptscriptstyle\circWscriptscript{}^{\kern1.5pt #2}_{\kern-1pt #1}}%
}}

\def\OM{\mathchoice
{\rlap{\kern3.2pt$\overline{\phantom{L}}$}M}
{\rlap{\kern3.2pt$\overline{\phantom{L}}$}M}
{\rlap{\kern2.4pt$\scriptstyle\overline{\phantom{L}}$}M}
{\rlap{\kern1.8pt$\scriptscriptstyle\overline{\phantom{L}}$}M}}

\def\mycirc{{\kern1pt\circ\kern2pt}}

\let\bbar\hat

\def\Spec{\mathop{\rm Spec}\nolimits}

\def\deg{\mathop{\rm deg}\nolimits}

\def\PGL{\mathop{\rm PGL}\nolimits}

\let\phi\varphi
\let\theta\vartheta
\let\epsilon\varepsilon
\let\setminus\smallsetminus

\newcommand{\BP}{{\mathbb{P}}}
\newcommand{\BQ}{{\mathbb{Q}}}
\newcommand{\BR}{{\mathbb{R}}}

\newcommand{\BZ}{{\mathbb{Z}}}

\newcommand{\Fm}{{\mathfrak{m}}}

\newcommand{\CC}{{\cal C}}

\newcommand{\CP}{{\cal P}}

\newcommand{\CU}{{\cal U}}

\newcommand{\CX}{{\cal X}}
\newcommand{\CY}{{\cal Y}}

\newcommand{\op}{\operatorname}

\newcommand{\wbar}{{\overline{w}}}
\newbox\mybox
\def\arrover#1{\mathrel{
\setbox\mybox=\hbox spread 1.4em
{\hfil$\scriptstyle#1$\hfil}
\vbox{\offinterlineskip\copy\mybox
\hbox to\wd\mybox{\rightarrowfill}}}}

\def\larrover#1{\mathrel{
\setbox\mybox=\hbox spread 1.4em
{\hfil$\scriptstyle#1\vphantom{g}$\hfil}
\vbox{\offinterlineskip\copy\mybox
\hbox to\wd\mybox{\leftarrowfill}}}}

\def\ontoover#1{\mathrel{
\setbox\mybox=\hbox spread 1.4em
{\hfil$\scriptstyle#1\vphantom{g}$\hfil}
\vbox{\offinterlineskip\copy\mybox
\hbox to\wd\mybox{\rightarrowfill\hskip-2.8mm
$\rightarrow$}}}}
\def\leftontoover#1{\mathrel{
\setbox\mybox=\hbox spread 1.4em
{\hfil$\scriptstyle#1\vphantom{g}$\hfil}
\vbox{\offinterlineskip\copy\mybox
\hbox to\wd\mybox{$\leftarrow$\hskip-2.8mm
\leftarrowfill}}}}
\let\longto\longrightarrow

\let\onto\twoheadrightarrow

\def\isoto{\mathrel{
\setbox\mybox=\hbox spread 0.9em
{\hfil$\scriptstyle\sim$\hfil}
\vbox{\offinterlineskip\copy\mybox
\hbox to\wd\mybox{\rightarrowfill}}}}

\newtheorem{Thm}{Theorem}[section]
\newtheorem{Prop}[Thm]{Proposition}
\newtheorem{Proposition}[Thm]{Proposition}
\newtheorem{proposition}[Thm]{Proposition}
\newtheorem{Lem}[Thm]{Lemma}
\newtheorem{Lemma}[Thm]{Lemma}

\newtheorem{Cor}[Thm]{Corollary}
\newtheorem{Conj}[Thm]{Conjecture}
\newtheorem{Def}[Thm]{Definition}

\newtheorem{ThmDef}[Thm]{Theorem-Definition}
\newtheorem{Rem}[Thm]{Remark}
\newtheorem{Ex}[Thm]{Example}

\newtheorem{thmx}{Theorem}

\newtheorem{propx}[thmx]{Proposition}

\numberwithin{Thm}{subsection}
\def\UseTheoremCounterForNextEquation{\setcounter{equation}{\value{Thm}}\addtocounter{Thm}{1}}

\def\qed{{\hskip0pt\unskip\unskip\nobreak\hfil\penalty50
\hskip1em\hbox{}\nobreak\hfil
{$\square$}
\parfillskip=0pt\finalhyphendemerits=0
\par}\medskip}

\newenvironment{Proof}
{\noindent{\bf Proof.}\ }
{\qed}

{\noindent{\bf Proof of #1.}\ }
{\qed}

\begin{document}

\title{\strut
	\vskip-80pt
	Computing the Stable Reduction of \\
	 Hyperelliptic Curves in Residue Characteristic $2$\\
}
\author{
	\begin{minipage}{.3\hsize}
		Tim Gehrunger\\[12pt]
		\small Department of Mathematics \\
		ETH Z\"urich\\
		8092 Z\"urich\\
		Switzerland \\
		tim.gehrunger@math.ethz.ch\\[9pt]
	\end{minipage}
%	\qquad
%	\begin{minipage}{.3\hsize}
%		Richard Pink\\[12pt]
%		\small Department of Mathematics \\
%		ETH Z\"urich\\
%		8092 Z\"urich\\
%		Switzerland \\
%		pink@math.ethz.ch\\[9pt]
%	\end{minipage}
}
\date{\today}
%\date{October 11, 2161}

\maketitle

%\Bigskip
\centerline{In memory of  Hjalti Þór Ísleifsson}
\bigskip\bigskip
\begin{abstract}
	%We construct the stable marked model of a hyperelliptic curve over a discrete valuation ring of residue characteristic $2$ by direct computations. For genus $2$ cuves, we work out relatively simple conditions for their stable unmarked reduction.
%	Consider a hyperelliptic curve of genus $g$ over a field $K$ of characteristic zero. After extending $K$ we can view it as a marked curve with its $2g+2$ Weierstrass points. We present an explicit algorithm to compute the stable reduction of this marked curve for a valuation of residue characteristic~$2$ over a finite extension of~$K$. In the cases $g\le2$ we work out relatively simple conditions for the structure of this reduction.
Consider a hyperelliptic curve of genus $g$ over a field $K$ of characteristic zero. After extending $K$ we can view it as a marked curve with its $2g+2$ Weierstrass points. We prove some general properties of the stable reduction of this marked curve for a valuation of residue characteristic~$2$ and refine an existing algorithm for its computation. We work out explicit examples up to genus $g=30$. 
% an explicit algorithm to compute the stable reduction of this marked curve for a valuation of residue characteristic~$2$ over a finite extension of~$K$. In the cases $g\le2$ we work out relatively simple conditions for the structure of this reduction.
\end{abstract}

{\renewcommand{\thefootnote}{}
	%11F52 Modular forms associated to Drinfeld modules
	%11G09 Drinfeld modules; higher-dimensional motives, etc.
	%14D20 Algebraic moduli problems, moduli of vector bundles
	%14D22 Fine and coarse moduli spaces
	%14M27 Compactifications; symmetric and spherical varieties
}

\newpage
\tableofcontents
\newpage
%\newpage
%%\renewcommand{\baselinestretch}{0.6}\normalsize
%\tableofcontents
%%\renewcommand{\baselinestretch}{1.0}\normalsize
%\newpage
%\htbegin 
%All of this is in a pretty rough state edit-wise. Some of the statements may be stated in a sub-optimal way and require rephrasing. For Example the smooth unmarked points section could also use the Lemmas of the double points section, ToDo. 
%\htend \medskip 

\section{Introduction}
\label{Intro}

{\bf 1.1 Motivation and strategy:}
Let $K$ be a valued field of characteristic $0$ and residue characteristic $2$. Moreover, let $C$ be a hyperelliptic curve over~$K$, that is, a curve defined by a Weierstrass equation $z^2=F(x)$ for some polynomial~$F$. After extending $K$ if necessary, the ramification points of the canonical double covering $\pi\colon C\onto\bar C \cong \mathbb{P}^1_K$ are defined over $K$ and there exists a stable marked model $\CC$ of $C$ with these points marked. The special fiber of this model encodes many arithmetic properties of $C$. 

In \cite{GP24}, Richard Pink and the author of this article provide an algorithm that describes $\CC$ explicitly. However, the algorithm becomes infeasible to perform for large $g$, even when using computer algebra systems. This is mainly due to the requirement to approximate roots of a polynomial of large degree, which we call the stability polynomial.
In this article, we study properties of $\CC$ and use these to reduce the computational complexity by allowing us to avoid the stability polynomial or to replace it with a polynomial of smaller degree in many cases. We then use our methods to compute some explicit examples up to genus $g=30$.

%, allowing for an implicit analysis of the stable reduction of $C$.
%%%%%%%%%%%%%%%%%%%%%%%%%%%%%%%%%%%%%
\medskip
\textbf{1.2 Overview:} 
We now explain the content of this article in greater detail. First, we reduce the general case throughout to the case that $K$ is algebraically closed. Let $R$ denote its valuation ring with maximal ideal $\Fm$ and $k=R/\Fm$ its residue field. 

The advantage of working with marked models is that one can use the stable model of the rational curve $\bar C$ marked by the branch points of~$\pi$. These points are the zeros of $F$ and possibly the point $x=\infty$, and the stable model $\bar\CC$ of $\bar C$ as a marked curve can by computed easily.
It turns out that the model $\CC$ dominates $\bar\CC$ and that some of its complexity is already present in~$\bar\CC$. 
%Let $\hat \CC$  be the minimal model of $\bar C$ dominating  $\bar \CC$ such that its normalization in $K(C)$ is semistable. 
The first step in the algorithm from \cite{GP24} is to compute $\bar \CC$. Next, one determines the minimal model $\hat \CC$ of $\bar C$ dominating  $\bar \CC$ such that its normalization in $K(C)$ is semistable.
%Then $\CC$ dominates $\hat \CC$. 
  We denote the special fibers of $\CC, \hat \CC$ and $\bar \CC$ by $C_0, \hat C_0$ and $\bar C_0$ respectively and collect the schemes we have introduced in the following diagram:
\begin{equation*}\label{AllCPDiagram}
	\vcenter{\xymatrix{
			\ C\  \ar@{^{ (}->}[r] \ar@{->>}[d]_-\pi & 
			\ \CC\ \ar@{->>}[d]_-{} & 
			\ C_0\ \ar@{_{ (}->}[l] \ar@{->>}[d]_-{}\\
			\ \bbar C\ \ar@{^{ (}->}[r] \ar@{=}[d] & 
			\ \bbar\CC\ \ar@{->>}[d] & 
			\ \bbar C_0\ \ar@{->>}[d] \ar@{_{ (}->}[l]  \\
			\ \bar C\ \ar@{^{ (}->}[r] \ar@{->>}[d] &
			\ \bar\CC\ \ar@{->>}[d] & 
			\ \bar C_0\ \ar@{_{ (}->}[l] \ar@{->>}[d] \\
			\ \Spec K\ \ar@{^{ (}->}[r] & 
			\ \Spec R\  & 
			\ \Spec k\ \ar@{_{ (}->}[l] \\}}
\end{equation*}

\textbf{Local genus} 
For any closed point $\hat p \in \hat C_0$, we define the local genus $g_{\hat p}$ of $C_0$ over $\hat p$ as $1$ if its preimage in $C_0$ contains two double points and as $0$ otherwise. Moreover,  if $\hat p \in \hat C_0$ is the generic point of an irreducible component $X$, we set $g_{\hat p}=g'$ if the preimage of $X$ in $C_0$ is irreducible of geometric genus $g'$ and $g_{\hat p}=-1$ when the preimage is a disjoint union of two rational curves. For any point $\bar p \in \bar C_0$, we define $g_{\bar p}$ as the sum of the local genera of the elements in the preimage of $\bar p$ in $\hat C_0$. In the case that $\bar p\in \bar C_0$ is a closed point, the local genus $g_{\bar p}$ is the sum of the geometric genera of the irreducible components of $C_0$ which are fully contained in the preimage of $\bar p$.

In Proposition \ref{LocalGenusSumIsGFormula}, we establish
$$g \ = \ \sum_{\bar p\in \bar C_0}g_{\bar p} .$$
Moreover, we express the local genus in elementary formulas involving a Weierstrass equation of $C$, which allows us to derive properties of $C_0$ without computing $\CC$ first.  Let $X$ be an irreducible component of $\bar C_0$ which contains exactly one double point $\bar p$ and $n$ marked points. Then by Propositions \ref{EvenLeavGenusProp} and \ref{OddLeavGenusProp}, we have
$$\sum_{\bar q \in X\setminus\{\bar p\}} g_{\bar q}\leq \frac{n-1}{2}, $$ where equality holds if and only if $n$ is odd.

\medskip 

\textbf{Approximating Laurent polynomials by squares} 
%The first step in the algorithm from \cite{GP24} is to compute $\bar \CC$. Next, one determines the minimal model $\hat \CC$ of $\bar C$ dominating  $\bar \CC$ such that its normalization in $K(C)$ is semistable. 
It turns out that computing $\hat \CC$ is a local problem on the special fiber of $\bar \CC$, and is done by explicitly computing the normalization of finitely many open subsets of affine lines in $K(C)$.  To explain this computation, we introduce more notation.

Let $v$ denote the valuation on $K$ that is normalized to $v(2)=1$. For any rational number $\alpha$ we choose a suitable fractional power $2^\alpha\in K$ with $v(2^\alpha)=\alpha$. For any Laurent polynomial $F$ over $K$ we let $v(F)$ denote the minimum of the valuations of its coefficients.

In \cite{GP24}, for any Laurent polynomial $F\in R[x^{\pm1}]$ we defined
$$w(F)\ :=\ \sup \bigl\{ v(F-H^2) \bigm| H\in R[x^{\pm 1}] \bigr\}\ \in\ \BR\cup\{\infty\},$$
which measures how well $F$ can be approximated by squares. A decomposition of the form $F=G+H^2$ with $G,H\in R[x^{\pm 1}]$ is called \emph{optimal} if
$$v(G)=w(F)\quad\hbox{or}\quad v(G)>2.$$
Now we fix an optimal decomposition $F=H^2+G$ and let $\gamma:=\min\{2,w(F)\}.$
By \cite[Prop. 3.3.2]{GP24}, the normalization of $R[x^{\pm1}]$ in $K(x,z)$ with $z^2=F$ is isomorphic to
$$B\ :=\ R[x^{\pm1}][t]\!\bigm/\!\bigl(2^{1-\gamma/2}Ht+t^2-g/2^\gamma\bigr),$$
which we use in \cite{GP24} to compute the  normalizations of the open subsets of $\bar \CC$ in $K(C)$.

%In the present article, we show how the computation of such normalizations can be simplified by approximating polynomials $2$-adically. Moreover, we introduce the concept of the local genus over a point and show how to compute it without computing $\CC$.

\medskip
\textbf{Stable reduction over smooth closed points} 
Now assume that $\bar p$ is a smooth closed point contained in an irreducible component $X\subset \bar C_0$. Then there exist $F\in R[x^{\pm 1}]$ with $v(F)=0$ and a neighborhood $\bar \CU$ of $\bar p$ in $\bar \CC$ isomorphic to an open subscheme of $\op{Spec}(R[x])$ such that $z^2=F$ is a Weierstrass equation for $C$ and $\bar p$ is given by $x=0$.
%Computing the preimage of $\bar p$ in $\hat C_0$ using the algorithm in \cite{GP24} requires us to study the zeros of the stability polynomial with positive valuation, which gets 
Let $\tilde F\in R[x^{\pm 1}]$ with $v(F-\tilde F)>2$ and let $\tilde C$ be the curve defined by $z^2=\tilde F$. We show in Propositions \ref{BareNormalizationSpecialFiberIsIsomorphic} and \ref{stableReductionOverSmoothUnmarkedPoint} that the preimage $\tilde \CU$ of $\bar \CU$ in $\hat \CC$ is the minimal semistable model of $\bar \CU$ whose normalization in $K(\tilde C)$ is semistable. Moreover, for each component of the special fibers of  $\tilde \CU$, their preimages in the normalization in $K(C)$ and $K(\tilde C)$ are isomorphic. This means that when computing $C_0$ over $\bar p$ using the algorithm from \cite{GP24}, we can replace $F$ by $\tilde F$. This is fruitful in many cases, especially when the degree of $\tilde F$ is smaller than that of $F$. Among other things, this reduces the computational complexity of determining $C_0$ with the algorithm in \cite{GP24} due to the stability polynomial associated to $\tilde F$ having a smaller degree than that associated fo $F$.

Moreover, we define the square defect  $\wbar(X):=\min\{2, w(F)\}$ of $X$ which is independent of the choices of $\CU, \bar p$ and $F$. 
%We prove a bound for the thickness of double points above $\bar p$ in $\hat C_0$ in terms of $\wbar(X)$ and $g_{\bar p}$. 
 Fix an irreducible component $\bbar T$ of $\hat C_0$ above~$\bar p$ whose preimage $T$ in $C_0$ has positive geometric genus $g(T)$ for the genus of $T$. Furthermore, let $\CY$ be the semistable model of $\bar C$ obtained from blowing down all components above $\bar p$ other than $\bbar T$. Then $\CY$ has a unique double point $\bar q$ above $\bar p$.
\begin{propx}[= Proposition \ref{ComputeThickness}]
	Let $\epsilon$ be thickness of $\bar q$. Then $$\epsilon\leq \frac{2-\wbar(X)}{2g(T)+1}.$$ Moreover, equality holds if and only if $\bbar T$ is the only component of $\hat C_0$ above $\bar p$.  
\end{propx}

Now let  $F=H^2+G$ be an optimal decomposition.
\begin{thmx}[=Theorem \ref{LocalGenusUnmarkedPointProp}]
	Let $s$ be the multiplicity of the root of $[\frac{dG}{dx}/2^\gamma]$ at $\bar p$. Then $s$ is even and we have
$$\frac{s}{2}=g_{\bar p}.$$
\end{thmx}

%By Theorem \ref{LocalGenusUnmarkedPointProp}, we have  $$[\frac{dG}{dx}/2^{\wbar(X)}]\in k[x^2]$$ and the  local genus $g_{\bar p}$ of $C_0$ over $\bar p$ is  half the multiplicity of the root of $[\frac{dG}{dx}/2^\gamma]$ at $\bar p$. 

In many situations, Theorem \ref{LocalGenusUnmarkedPointProp} is enough to determine the combinatorial structure of the preimage of $\bar p$ in $C_0$. If for example $g_{\bar p}=1$, we know that there is exactly one component of $C_0$ over $\bar p$, and this component has genus $1$. 

\medskip

\textbf{Stable reduction over double points} 
Next, assume that $\bar p\in \bar C_0$ is a double point. Then there exist $\alpha>0$ and a neighborhood $\bar\CU \subset \bar \CC$ isomorphic to an open subscheme of $\Spec R[x,y]/(xy-2^\alpha)$ such that $\bar p$ is given by $x=y=0$. 
Identifying $y$ with $\frac{2^\alpha}{x}$, the ring in question is isomorphic to the subring $R[x,\tfrac{2^\alpha}{x}]$ of the ring of Laurent polynomials $K[x^{\pm1}]$. 
Moreover, there exists $F  \in R[x,\frac{2^\alpha}{x}]$ with constant coefficient $1$ such that function field of $C$ is $K(x,z)$ with $z^2=F (x)$. Let $X$ be the irreducible component of $\bar C_0$ for which $x$ is a coordinate. We say that $\bar p$ is  \emph{even} if each connected component of $\bar C_0 \setminus \{\bar p\}$ contains an even number of the points $\bar p_1,\ldots,\bar p_{2g+2}$. Otherwise, it is called \emph{odd}. If $\bar p$ is odd, we have $g_{\bar p}=0$ by Proposition \ref{LocalGenusTrivialForMarkedAndOddDP}.

Now assume that $\bar p$ is an even double point and let $u$ be a new variable. Let $\tilde F\in R[x^{\pm 1}]$ such that  for $\lambda \in [0,\alpha]\cap \BQ$, we have 
$$v(F(2^\lambda u)-\tilde F(2^\lambda u)) >2.$$ 
Let $\tilde C$ be the curve defined by $z^2=\tilde F$ and let $\tilde \CU$ be the minimal semistable model of $\bar \CU$ whose normalization in $K(\tilde C)$ is semistable. Then we show in Proposition \ref{ComponensofTypeCDApproxProp} that the preimage of $\bar \CU$ in $\hat \CC$ is equal to $\tilde \CU.$ As before, this means that when computing $C_0$ over $\bar \CU$ using the algorithm from \cite{GP24}, we can replace $F$ by $\tilde F$.

%	Moreover, the preimages of the components of the special fiber of $\tilde \CU$ of type (b), (c) and (d) above $\bar p$ in  the normalizations in $K(\tilde C)$ are isomorphic to their preimages in $\CC$.  
% that regarding the computation of the stable reduction of $C$ above $\bar \CU$, we can replace $F$ by $\tilde F$. 
We define the square defect function at $\bar p$ with respect to $X$ by  
\begin{equation*}
	\wbar_{\bar p, X}\colon\ \BQ\,\cap\, [0, \alpha] \longto\BR,\ \ \lambda\mapsto 
	\wbar_{\bar p, X}(\lambda) := \wbar(F(2^\lambda u)).
	% \min \{ 2, w(F(2^\lambda u)) \}.
\end{equation*}
By  \cite[Prop. 3.4.4]{GP24} and Theorem-Definition \ref{Square-defect}, this is a piecewise linear concave function. 
Together with \cite[Prop. 4.5.1]{GP24}, its break-points correspond to the irreducible components of $\hat C_0$ above $\bar p$ that lie between the proper transforms of the two irreducible components of $\bar C_0$ that meet at $\bar p$.  
\begin{thmx}[=Theorem \ref{RelativeGenusOverDoublePoint}]
Let $s_1$  and $s_2$ be the largest, respectively the smallest slope of $\wbar_{\bar p, X}$. Then the local genus of $C$ over $\bar p$ is
$$g_{\bar p}=\frac{s_1-s_2+\delta_{0s_1}+\delta_{0s_2}}{2},$$
where $\delta_{0s_i}$ is $0$ if $s_i\neq 0$ and $1$ otherwise.  	
\end{thmx}

An even double point $\bar p\in \bar C_0$ is called {\it grounded} if for the irreducible components $X,Y$ of $\bar C_0$ containing it, we have  $\wbar(X)=\wbar(Y)=0$. For such a double point, we have $g_{\bar p}>0$ by Proposition \ref{GroundedDPAtlestgenus1}. Moreover, if $g_{\bar p}=1$, we establish in Proposition \ref{Genus1BehaviouroverGroundedDP} that the reduction behavior above $\bar p$ only depends on the thickness of $\bar p$ and is as in the genus $g=1$ case.  

\medskip

\textbf{1.3 Examples:} 
Applying our results regarding the local genus of $C_0$ to the case of genus $3$ curves, we find 115 cases in which the stable marked reduction of $C$ only depends on the thickness of its double points. This count is not exhaustive, and we expect that there are many more cases with this property. 

Combining  local genus computations with our approximation results, we explicitly work out the stable marked reduction of five curves with many automorphisms in the sense of \cite{MuellerPink}. The highest of these curves has genus $30$. 

%%%%%%%%%%%%%%%%%%%%%%%%%%%%%%%%%%%%%
\medskip
{\bf 1.4 Structure of the paper:}
% From here on we assume that $R$ has residue characteristic~$2$.
Section \ref{SetUpChapter} contains preparatory material: In Subsection~\ref{SubsectionValuationPolynomials} we review valuations and optimal decompositions of Laurent polynomials. 
In the following two Subsections \ref{ModelsofCurvesSubSection} and \ref{SubSectionHyperellipticCurves} we review basic facts about semistable and stable marked models of hyperelliptic curves over~$R$. 
In the final Subsections \ref{SubSectionSDF} and \ref{SubSectionSDI} we define the square defect and square defect function and study their properties.

In Section \ref{ApproximationSection} we study how to approximate $F$ in the computation of $C_0$ from $\bar C_0$.
In Subsection \ref{SubSectionOptimalDecompositions}, we study how optimal decompositions of two Laurent polynomials whose difference is small relate to each other. In Subsection \ref{SubSectionSmoothPointApproximation}, we apply this to the computation of $C_0$ above smooth points and in Subsection \ref{SubSectionDoublePointApprox} we apply it to the computation of $C_0$ above double points. 

The local genus is the topic of Section \ref{SectionLocalGenus}. In Subsection \ref{SubSectionLocalGenus} we define the local genus, show that it is well-defined and trivial for marked points and odd double points. In Subsection \ref{SubSectionsLocalGenusUnmarkedPoints} we prove a formula for the local genus over a smooth unmarked point and in Subsection \ref{SubSectionEvenDoublePointsLocalGenus} we do the same for even double points. 

 In Section \ref{SectionComputationalCriterions} we provide tools for the computation of invariants of $C_0$ that do not require to compute $\CC$ first. In Subsection \ref{SubSectionThicknessBound} we prove a bound on the thickness of certain double points of $\hat C_0$ and $C_0$. In Subsection \ref{SubSectionGroundedDoublePoints} we introduce the notion of grounded double points and study their properties. We also briefly discuss the reduction behavior of curves containing as many such points as possible. In Subsections \ref{SubSectionLeafCompoenents} and \ref{SubSectionConjectureOddDoublePoints}, we study restrictions on the distribution of the local genera of points on the connected components of $C_0$ minus a double point.

In the final Section \ref{SectionExamples} we apply these methods to work out the reduction behavior of all genus $3$ curves containing  three grounded double points and for five curves with many automorphisms and of high genus.
 For the latter curves, we only list the final results, leaving the detailed computations to look up in the associated computer algebra worksheets~\cite{WorksheetsG24}.

 %%%%%%%%%%%%%%%%%%%%%%%%%%%%%%%%%%%%%
 \medskip
 {\bf 1.5 Relation with other work:}
 This article builds on \cite{GP24}, which was written by Richard Pink and the author of this article and in which an explicit algorithm for the computation of $\CC$ is given. Similar approaches can be found in the work of Arzdorf and Wewers in \cite{ArzdorfWewers} in the language of Berkovich analytic spaces and in the preprint \cite{FioreYelton} of Fiore and Yelton. In the earlier work \cite{LehrMatignon06} of Lehr and Matignon, the special fiber was described fully in the case that $\bar C_0$ is smooth. 
 
 All of these articles are based on the work of Raynaud, which in a series of articles \cite{Raynaud1970}, \cite{Raynaud1990} and \cite{Raynaud1999}, studied the case of mixed characteristic $(0, 2)$ extensively, also under the assumption that $\bar C_0$ is smooth. 
 
 In \cite{Henrio2000}, Henrio studied local lifting problems and introduced the notion of Hurwitz trees, which contains the information required to lift a curve over $k$ to a curve over $R$. The Hurwitz tree of a curve includes a value  $\delta\in [0,1]$ for every irreducible component $X$ of the curve, which is related to our square defect via $\delta=\frac{2-\wbar(X)}{2}$. 
In a series of articles,   \cite{Saidi1998} and \cite{Saidi2003}, Sa\"idi generalized some of this work to the global situation and introduced a notion very similar to our local genus over points.

 Historically, semistable reductions of hyperelliptic curves have mainly been studied and constructed in residue characteristic $\not =$ 2, see for example the construction of Bosch in \cite{Bosch1980}. A more recent approach is the article \cite{DDMM} by Dokchitser, Dokchitser, Maistret, and Morgan, which describes the special fiber in their notion of cluster pictures. Similarly, in \cite{GP21}, Richard Pink and the author of the present article have given a description of the stable marked reduction.

 Moreover, in  \cite{Liucriterionarticle}, Liu gives criteria for the type of the stable reduction of the unmarked curve $C$ in terms of Igusa invariants. This result is first proved in the setting of $\op{char}(k)\neq 2$ and carries over to the wild case by a moduli argument.

\section{Set-up}\label{SetUpChapter}

This article uses the notation and conventions of \cite{GP24}, which was written by Richard Pink and the author of this article. We briefly review these conventions. %in this section.%For more details and proofs, we refer to this article and the references therein. 
Throughout this article let $R_1$ be a complete discrete valuation ring with quotient
field $K_1$. At several places we will need to replace $R_1$ by its integral closure in a finite
extension of $K_1$. As in \cite{GP24}, we find it more convenient to
work over an algebraic closure instead. So we fix an algebraic closure $K$ of $K_1$ and let $R$
denote the integral closure of $R_1$ in $K$. Since $R_1$ is complete, the valuation on $K_1$ extends
to a unique valuation with values in $\BQ$ on $K$, whose associated valuation ring is $R.$
Let $v$ denote a corresponding valuation on $K$, and let $\Fm$ be the maximal ideal and $k := R/\Fm$ the
residue field of $R$.
By \cite[Thm. 2.1.5]{GP24} we can work over
the single field $K$ and avoid cumbersome changes of notation. 

\medskip

From now on we assume that $K$ has characteristic $0$ and $k$ has characteristic $2$. We normalize the valuation $v$ on $K$ in such a way that $v(2)=1$. For every integer $n\ge1$ we fix an $n$-th root $2^{1/n}\in K$ in a compatible way, such that for all $n,m\ge1$ we have $(2^{1/mn})^m=2^{1/n}$. For any rational number $\alpha = m/n$ we then set $2^\alpha := (2^{1/n})^m$. This defines a group homomorphism $\BQ \to K^\times$, which by the normalization of $v$ satisfies $v(2^\alpha)=\alpha$.

\subsection{Valuations of Polynomials}\label{SubsectionValuationPolynomials}
For any Laurent polynomial $F=\sum_i a_i x^i  \in K[x^{\pm 1}]$ we set
\UseTheoremCounterForNextEquation
\begin{equation}\label{vDef}
	v(F)\ :=\ \inf\,\{v(a_i) \,|\, i\in\BZ\}.
\end{equation}
This extends $v$ to a valuation on $K[x^{\pm 1}]$, which by the Gauss lemma satisfies the equation $v(FH)=v(F)+v(H)$ for all $F,H\in K[x^{\pm 1}]$.
The elements $F$ with $v(F)\ge0$ make up the subring $R[x^{\pm1}]$,
%The elements with $v(F)=0$ are those elements of $R[x^{\pm1}]$ that are not congruent to $0$ modulo~$\Fm$.
%$$\begin{array}{rl}
	%K[x^{\pm 1}]^{\leq d_1}_{\geq d_2} &:=\ 
	%\bigl\{ g \in K[x^{\pm1}] \bigm| \op{ord}(g)\geq d_2 \wedge \op{deg}(g)\leq d_1 \bigr\}, \\[3pt]
	%R[x^{\pm 1}]^{\leq d_1}_{\geq d_2} &:=\ 
	%\bigl\{ g \in R[x^{\pm1}] \bigm| \op{ord}(g)\geq d_2 \wedge \op{deg}(g)\leq d_1 \bigr\}.
	%\end{array}$$
	and for any such $F$ we let $[F]$ denote its residue class in $k[x^{\pm1}]$. We define 
	\UseTheoremCounterForNextEquation
	\begin{align}
		\begin{split}
	\op{deg}(F)\ &:= \ \sup \{i\in \BZ \ : \ a_i \neq 0\}; \\
	\op{ldeg}(F)\ &:= \ \inf \{i\in \BZ \ : \ a_i \neq 0\}.\label{ldegreeDefEq}
	\end{split}
	\end{align}
	%	\medskip
	%	ToDo: Add what ldegree means. 
	%%%%%%%%%%%%%%%%%%%%%%%%%%%%%%%%%%%%%%%%%%
	%	\subsection{Optimal decompositions}
	%	\label{OptDecomp}
	%	
	
	Now fix a Laurent polynomial $F\in R[x^{\pm 1}]$ with $F \not \equiv 0 \mod\Fm$, in other words with $v(F)=0$. To this we associate the values
	\UseTheoremCounterForNextEquation
	\begin{equation}\label{wDef}
		w(F)\ :=\ \sup \bigl\{ v(F-H^2) \bigm| H\in R[x^{\pm 1}] \bigr\}\ \in\ \BR\cup\{\infty\}
	\end{equation}
	and $\wbar(F):= \min\{w(F),2\}$, which measure how well $F$ can be approximated by squares.
	
	\begin{Def}\label{opOddtDef}
		A decomposition of the form $F=G+H^2$ with $G,H\in R[x^{\pm 1}]$ is called
		\begin{itemize}
			\item 
			\textbf{odd} if $G\in xR[x^{\pm 2}]$; 
			\item	
			\textbf{optimal} if it satisfies the condition
			\UseTheoremCounterForNextEquation
			$$v(G)=w(F)\quad\hbox{or}\quad v(G)>2.$$
		\end{itemize}
	\end{Def}
	By \cite[Prop. 3.2.4]{GP24}, an odd decomposition of $F$ exists. 
	Fix such a decomposition $F=G+H^2$ and let $\gamma:=v(G)$. By  \cite[Prop. 3.2.5]{GP24}, we have $\gamma=\min \{\wbar(F), 2\}$, which implies that the decomposition is optimal if $\gamma<2$. 
	%	By Proposition 3.2.5 of \cite{GP24}, an odd decomposition $F=H^2+G$ with $v(G)<2$ is optimal. 
	%	ToDo: Add odd decomposition. 

	%	In this section, let $F\in R[x, \frac{2^\alpha}{x}]$ for $0<\alpha\in \BQ$  and consider an odd decomposition $F=H^2+G$. Set $\gamma:=\wbar(F)$. ToDo: maybe make more clear that these things are defined? 
	
	\begin{Lem}\label{ExchangeOptimalForOddLem}
		Assume  $\gamma<2$. 
		For any optimal decomposition $F=\tilde G + \tilde H^2 $, we have
		$$\bigg[\frac{dG}{dx}/2^\gamma\bigg]=\bigg[\frac{d\tilde G}{dx}/2^\gamma\bigg].$$
	\end{Lem}
	\begin{Proof}
		As both decompositions are optimal, we have $v(G)=v(\tilde G)=\gamma$. This yields 
		$$v(H^2-\tilde H^2)=v(G-\tilde G)\geq \gamma. $$
		As $H^2- \tilde H^2= (H-\tilde H)(H+\tilde H)$ and $\gamma<2$, we get $v(H-\tilde H)=\frac{1}{2}v(H^2-\tilde H^2)\geq \frac{\gamma}{2}$. 
		Hence there exists   $L \in R[x^{\pm 1}]$ with $H=\tilde H + L 2^{\gamma/2}$. Then
		$$G=F+H^2=F+\tilde H^2 + 2^{1+\gamma/2} L +L^2 2^\gamma=\tilde G+ 2^{1+\gamma/2} L +L^2 2^\gamma.$$
		As $\gamma<2$, this yields 
		$$\bigg[\frac{dG}{dx}/2^\gamma\bigg]=\bigg[\frac{d\tilde G}{dx}/2^\gamma\bigg]+ [2^{1-\gamma/2}L'+2 L L']=\bigg[\frac{d\tilde G}{dx}/2^\gamma\bigg], $$
		as desired. 
	\end{Proof}
	
	\medskip
	
	%	ToDo: Maybe eliminate section? 
	%	\subsection{Valuations under Scaling}\label{SlopesSection}
	Writing $F=\sum_ia_ix^i$ with $a_i\in R$, we now assume that the constant coefficient $a_0$ is a unit and that we are given a positive number $\alpha\in\BQ$ satisfying 
	\UseTheoremCounterForNextEquation
	\begin{equation}\label{AlphaDef}
		v(a_i)\ge\alpha|i|\ \ \ \hbox{for all}\ \ \  i<0. 
	\end{equation}
	Since $a_0$ is a unit, property (\ref{AlphaDef}) implies that all negative slopes of the Newton polygon of $F$ are $\le-\alpha$. Thus for any $\lambda\in\BQ\cap[0,\alpha]$, the Laurent polynomial
	$$F(2^\lambda u)\ =\ \sum_i a_i2^{\lambda i}u^i$$
	in the new variable $u$ again has coefficients in~$R$. %Fix an odd decomposition$F=H^2+G$.
	Consider the function
	\UseTheoremCounterForNextEquation
	\begin{equation}\label{WBarDef}
		\wbar\colon\ \BQ\,\cap\, [0, \alpha] \longto\BR,\ \ \lambda\mapsto 
		\wbar(\lambda) := \wbar(F(2^\lambda u)).
		% \min \{ 2, w(F(2^\lambda u)) \}.
	\end{equation}
	By Proposition 3.4.2  in \cite{GP24} we have 
	\UseTheoremCounterForNextEquation
	\begin{equation}\wbar(\lambda)=\min\{2,v(G(2^\lambda u))\}\label{wbarIdentity}\end{equation}
	for all $\lambda \in \BQ \cap [0,\alpha]$. In particular, the function $\wbar$ is piecewise linear. 	
	\begin{Lem}\label{slopeLemLeftSide}
		Let $\beta \in (0,\alpha]\cap \BQ$. Then the line segment of $\wbar$ to the left of $\beta$ has slope 
	$$	\begin{cases}
			 \deg([G(2^\beta u)/2^{\wbar(\beta)}]) &\text{if } \wbar(\beta)< 2; \\
			 \max \{\deg([G(2^\beta u)/2^{\wbar(\beta)}]),0\} &\text{otherwise.}
		\end{cases}$$

%		 If $\wbar(\beta)\neq 2$, the line segment of $\wbar$ to the left of $\beta$ has slope $\deg([G(2^\beta x)/2^{\wbar(\beta)}]$. If otherwise $\wbar(\beta)=2$, it has slope $\min \{\deg([G(2^\beta x)/2^{\wbar(\beta)}],0\}$. 
	\end{Lem}
	\begin{Proof} Consider the function $w: \lambda \mapsto v(G(2^\lambda u)))$, which is also piecewise linear. From the definition of $v$ on $R[u^{\pm 1}]$, it follows that the slope of $w$ to the left of $\beta$ is 
		$\deg([G(2^\beta u)/2^{v(G(2^\beta u))}]).$
		
		If $\wbar(\beta)<2$, we have $\wbar(\beta)=w(\beta)$ by \eqref{wbarIdentity}. Furthermore, if $\deg([G(2^\beta x)/2^{\wbar(\beta))}])>0$, we have $\wbar(\beta)=w(\beta)$ as otherwise $[G(2^\beta u)/2^{\wbar(\beta))}]=0$.  In both cases, there exists $\epsilon>0$ such that $w$ is strictly smaller than $2$ on $I:=(\beta-\epsilon, \beta]$, which implies  $w_I=\wbar|_{I}$ by \eqref{wbarIdentity}. Hence the slope of $\wbar$ to the left of $\beta$ is equal to the slope of $w$ to the left of $\beta$, which is 
		$$\deg([G(2^\beta u)/2^{v(G(2^\beta u))}])\ =\ \deg([G(2^\beta u)/2^{\wbar(\beta)}]).$$
		
%		If $\wbar(\beta)=2$ and  $\deg([G(2^\beta x)/2^{\wbar(\beta)}]>0$, we have $w(\beta)=2$ and the slope of $w$ to the left of $\beta$ is $>0$, hence there exists $\epsilon>0$ such that $w$ is strictly smaller than $2$ on $I:=(\beta-\epsilon, \beta)$. Hence the slope of $\wbar$ to the left of $\beta$ is equal to the slope of $w$ to the left of $\beta$, which is $$\deg([G(2^\beta x)/2^{v(G(2^\beta x))}])\ =\ \deg([G(2^\beta x)/2^{\wbar(\beta)}]).$$
		
		If $\wbar(\beta)=2$ and  $\deg([G(2^\beta x)/2^{\wbar(\beta)}])<0$, we have $w(\beta)\geq 2$. Now either $w(\beta)>2$ or $w(\beta)=2$, in which case we have 
		$$\deg([G(2^\beta u)/2^{\wbar(\beta)}]=\deg([G(2^\beta u)/2^{w(\beta)}]<0$$
		and hence the slope of $w$ to the left of $\beta$ is $<0$. In either case, there exists $\epsilon>0$ such that the function $w$ is $\geq 2$ on $I:=(\beta-\epsilon, \beta]$. Thus by \eqref{wbarIdentity}, the restriction of $\wbar$ to $I$ is the constant function with value $2$. Hence the slope of $\wbar$ to the left of $\beta$ is $0$ in that case, as desired. 			
	\end{Proof}

	\begin{Lem}\label{slopeLemRightSide}
		Let $\beta \in [0,\alpha)\cap \BQ$. Then the line segment of $\wbar$ to the right of $\beta$ has slope 
			$$	\begin{cases}
			\operatorname{ldeg}([G(2^\beta x)/2^{\wbar(\beta)}]) &\text{if } \wbar(\beta)< 2; \\
			\min \{\operatorname{ldeg}([G(2^\beta x)/2^{\wbar(\beta)}]),0\} &\text{otherwise.}
		\end{cases}$$
%		
%		
%		 If $\wbar(\beta)\neq 2$, the line segment of $\wbar$ to the right of $\beta$ has slope $\op{ldeg}([G(2^\beta x)/2^{\wbar(\beta)}]$. If otherwise $\wbar(\beta)=2$, it has slope $\max \{\deg([G(2^\beta x)/2^{\wbar(\beta)}],0\}$. 
		%Old:
		%		If $\beta\neq \alpha$, there exists $\epsilon>0$ such that the function $\wbar_{f}|_{[\beta,\beta+\epsilon]}$ is linear. If $\wbar(\beta)<2$ or $\wbar(\beta)=2$ and $\op{ldeg}([G/2^{\wbar(\beta))}])$ is defined and negative, its slope is 
		%		$\op{ldeg}([G/2^{\wbar(\beta))}])>0$. Otherwise, it has slope $0$. 
	\end{Lem}
	\begin{Proof}
		Apply Lemma \ref{slopeLemLeftSide} to $F(x^{-1})$. 
	\end{Proof}

\subsection{Models of curves}\label{ModelsofCurvesSubSection}

Let $Y$ be a connected smooth proper algebraic curve over $K$. By a \emph{model of~$Y$} we mean a flat and finitely presented curve $\CY$ over $R$ with generic fiber~$Y$. We call such a model \emph{semistable} if the special fiber $Y_0$ is smooth except possibly for finitely many ordinary double points. Every double point $p\in Y_0$ then possesses an \'etale neighborhood in $\CC$ which is \'etale over $\Spec R[x,y]/(xy-a)$ for some nonzero $a\in \Fm$, such that $p$ corresponds to the point $x=y=0$. Here the valuation $v(a)$ depends only on the local ring of $\CY$ at~$p$, for instance by Liu \cite[\S10.3.2 Cor.\,3.22]{LiuAlgGeo2002}. Following Liu 
%\cite[\S5 Def.\,2]{Liu1993} 
\cite[\S10.3.1 Def.\,3.23]{LiuAlgGeo2002} we call $v(a)$ the \emph{thickness of~$p$}. 

%%%%%%%%%%%%%%%%%%%%%
\medskip
Any model is an integral separated scheme. Thus for any two models $\CY$ and $\CY'$ over~$R$, the identity morphism on $Y$ extends to at most one morphism $\CY\to\CY'$. If this morphism exists, we say that $\CY$ \emph{dominates}~$\CY'$. This defines a partial order on the collection of all models of $Y$ up to isomorphism. By blowing up one model one can construct many other models that dominate it. 

\subsection{Hyperelliptic curves}\label{SubSectionHyperellipticCurves}
%\label{HyperellCurves}

%ToDo: Clean this up and emphasize what we did in the other article and how we need to solve the local thing. 

Now let $C$ be a \emph{hyperelliptic curve} of genus $g$ over~$K$. Thus $C$ is a connected smooth proper algebraic curve which comes with a double covering $\pi\colon C \onto \bar C$ of a rational curve $\bar C\cong\BP^1_K$. Often the genus $g$ is required to be $\ge2$, but in this article we only assume $g\ge1$.

%%%%%%%%%%%%%%%%%%%%%
%\medskip
%
% By a \emph{model of~$C$} we mean a flat and finitely presented curve $\CC$ over $R$ with generic fiber~$C$. We call such a model \emph{semistable} if the special fiber $C_0$ is smooth except possibly for finitely many ordinary double points. Every double point $p\in C_0$ then possesses an \'etale neighborhood in $\CC$ which is \'etale over $\Spec R[x,y]/(xy-a)$ for some nonzero $a\in \Fm$, such that $p$ corresponds to the point $x=y=0$. Here the valuation $v(a)$ depends only on the local ring of $\CC$ at~$p$, for instance by Liu \cite[\S10.3.2 Cor.\,3.22]{LiuAlgGeo2002}. Following Liu 
%%\cite[\S5 Def.\,2]{Liu1993} 
%\cite[\S10.3.1 Def.\,3.23]{LiuAlgGeo2002} we call $v(a)$ the \emph{thickness of~$p$}. 

%Conversely, one can construct the \emph{contraction} of an irreducible component with the following properties:

\medskip
Consider a model $\bar\CC$ of $\bar C$ over~$R$. We say that a model $\CC$ of $C$ over $R$ \emph{dominates} $\bar\CC$ if and only if $\pi$ extends to a morphism $\CC\onto\bar\CC$.
% Since any model is an integral separated scheme, this morphism is then unique.
%By applying Proposition \ref{RelStabMod} in the case that $\CX$ is the normalization of $\bar\CC$ in the function field of~$C$, 
By general theory, there exists a minimal semistable model of $C$ that dominates~$\bar\CC$, and it is unique up to unique isomorphism.

\medskip
%Throughout the rest of this article we assume that $K$ has characteristic~$0$.
 As $K$ has characteristic $0$, the covering $\pi$ is tamely ramified, and by the Hurwitz formula it is ramified at precisely $2g+2$ closed points, namely, at the Weierstrass points of~$C$. Let $P_1,\ldots,P_{2g+2} \in C(K)$ denote these points and $\bar P_1,\ldots,\bar P_{2g+2} \in \bar C(K)$ their images under~$\pi$. Since $2g+2\ge4$, both $(C,P_1,\ldots,P_{2g+2})$ and $(\bar C,\bar P_1,\ldots,\bar P_{2g+2})$ are stable marked curves. 

\medskip
Let $(\CC, \CP_1, \dots, \CP_{2g+2})$ and $(\bar \CC, \bar \CP_1, \dots, \bar \CP_{2g+2})$
be the stable  models of  $(C, P_1, \dots, P_{2g+2})$   and $(\bar C, \bar P_1, \dots, \bar P_{2g+2})$, respectively. The next proposition is a special case of \cite[Prop. 2.2.8]{GP24}. 

\begin{Prop}\label{MinMarkModelStab}
	The model $\CC$ is the minimal semistable model of $C$ dominating $\bar \CC$. 
\end{Prop} 

%\medskip
Next let $\sigma$ denote the covering involution of $\pi\colon C\onto\bar C$. By the uniqueness of the minimal semistable model, this extends uniquely to an automorphism of $\CC$ of order~$2$. We denote this extension again by $\sigma$ and consider the quotient $\bbar\CC := \CC/\langle\sigma\rangle$. Since $\CC$ dominates $\bar\CC$, it follows that $\bbar\CC$ dominates~$\bar\CC$. Also $\sigma$ fixes each ramification point~$P_i$ and therefore each section~$\CP_i$. 

\medskip

In \cite{GP24}, we developed an algorithm to compute the model $\CC$. 
The primary goal in this article is to build on this work and simplify the explicit computation of 
%the stable model of $(C,P_1,\ldots,P_{2g+2})$. 
the special fiber of $\CC$. Proposition \ref{MinMarkModelStab} turns this into a local problem over the stable model of $(\bar C,\bar P_1,\ldots,\bar P_{2g+2})$. After this reformulation, the stability condition becomes irrelevant. The same method therefore solves the slightly more general problem of computing the minimal semistable model of $C$ that dominates~$\bar\CC$. Hence from now on, by $(\bar\CC,\bar\CP_1,\ldots,\bar\CP_{2g+2})$ we denote an arbitrary semistable model of $(\bar C,\bar P_1,\ldots,\bar P_{2g+2})$. Accordingly, by $(\CC,\CP_1,\ldots,\CP_{2g+2})$ we denote the minimal semistable marked model of $(C, P_1, \dots P_{2g+2})$ dominating $(\bar\CC,\bar\CP_1,\ldots,\bar\CP_{2g+2})$ and let $(\bbar\CC,\bbar\CP_1,\dots,\bbar\CP_{2g+2})$ the quotient of $\CC$ under $\sigma$ with the corresponding sections. 
%Throughout the following we therefore only work with the semistable marked models introduced above.

%%%%%%%%%%%%%%%%%%%%%
\medskip
For reference we collect the schemes and sections we have introduced in the following diagram. Recall that we have natural morphisms $\CC\onto\bbar\CC\onto\bar\CC$ that are compatible with the given sections. We let $(C_0, p_1, \dots, p_{2g+2})$ and $(\bbar C_0,\bbar p_1, \dots,\bbar p_{2g+2})$ and $(\bar C_0,\bar p_1, \dots,\bar p_{2g+2})$ denote the special fibers of $(\CC,\CP_1,\dots,\CP_{2g+2})$ and $(\bbar\CC,\bbar\CP_1,\dots,\bbar\CP_{2g+2})$ and $(\bar\CC,\bar\CP_1,\dots,\bar\CP_{2g+2})$, respectively.

\UseTheoremCounterForNextEquation
\begin{equation}\label{AllCPDiagram}
	\vcenter{\xymatrix{
			\ C\  \ar@{^{ (}->}[r] \ar@{->>}[d]_-\pi & 
			\ \CC\ \ar@{->>}[d]_-{} & 
			\ C_0\ \ar@{_{ (}->}[l] \ar@{->>}[d]_-{} &&
			\ P_i\ \ar@{^{ (}->}[r] \ar@{|->}[d] & 
			\ \CP_i\ \ar@{|->}[d] & 
			\ p_i\ \ar@{|->}[d] \ar@{_{ (}->}[l] \\
			\ \bbar C\ \ar@{^{ (}->}[r] \ar@{=}[d] & 
			\ \bbar\CC\ \ar@{->>}[d] & 
			\ \bbar C_0\ \ar@{->>}[d] \ar@{_{ (}->}[l] &&
			\ \bbar P_i\ \ar@{^{ (}->}[r] \ar@{=}[d] & 
			\ \bbar \CP_i\ \ar@{|->}[d] & 
			\ \bbar p_i\ \ar@{|->}[d] \ar@{_{ (}->}[l] \\
			\ \bar C\ \ar@{^{ (}->}[r] \ar@{->>}[d] &
			\ \bar\CC\ \ar@{->>}[d] & 
			\ \bar C_0\ \ar@{_{ (}->}[l] \ar@{->>}[d] &&
			\ \bar P_i\ \ar@{^{ (}->}[r] & 
			\ \bar\CP_i\ & 
			\  \bar p_i\ \ar@{_{ (}->}[l]\\
			\ \Spec K\ \ar@{^{ (}->}[r] & 
			\ \Spec R\  & 
			\ \Spec k\ \ar@{_{ (}->}[l] \\}}
\end{equation}

%%%%%%%%%%%%%%%%%%%%%
\medskip
To describe the relation between the special fibers $\bbar C_0$ and~$\bar C_0$, we use the terminology concerning the type 
%(a), (b), (c) or (d) 
of an irreducible component of $\bbar C_0$ from \cite[Def. 2.3.2]{GP24}. In particular, an irreducible component of $\hat C_0$ is called 
\begin{itemize}
	\item of type (a) if it maps isomorphically to an irreducible component of $\bar C_0$;
	\item of type (b) if it lies between irreducible components of type (a);
	\item of type (c) if it is not of type (a) or (b) and is not a leaf;
	\item of type (d) if it is not of type (a) or (b) and is a leaf.
\end{itemize} 
%ToDo: Move this proposition. 
%ToDo: The Proposition is due to R., even though he did not state it in this generality his proof still works. However, we can give a short proof with our methods.

We also divide the double points of $\bar C_0$ into two classes. For this recall that $\bar C_0$ is marked with $2g+2$ distinct points in the smooth locus. As the complement of a double point consists of two connected components, this divides the $2g+2$ marked points into two groups.
\begin{Def}\label{EvenOddDef}
	A double point $\bar p$ of $\bar C_0$ is called \emph{even} if each connected component of $\bar C_0 \setminus \{\bar p\}$ contains an even number of the points $\bar p_1,\ldots,\bar p_{2g+2}$. Otherwise, it is called \emph{odd}. 
\end{Def}

\subsection{Square defect of irreducible components}\label{SubSectionSDI}
Let $\CX$ be a semistable model of $\bar C$ and let $X$ be an irreducible component of its special fiber. Moreover, let $x$ be a global coordinate of $\CX$ along $X$ and consider a hyperelliptic equation $z^2=F(x)$ of $C$ with  $F\in R[x^{\pm1}]$ and $v(F)=0$. 
\begin{ThmDef}\label{Square-defect}
The square defect of $X$, denoted by $\wbar(X)$, is defined as  $$\wbar(X):=\wbar(F).$$ It is independent of the choice of $x$ and  $F$. 
\end{ThmDef}
\begin{Proof}
First, we show that $\wbar(F)$ is invariant under the action of $\PGL_2(R)$ which implies that the square defect does not depend on the choice of $x$.
  Let $F=H^2+G$ be an optimal decomposition of  $F$.
  
Consider the transformations $x \mapsto ax+b$ for $a \in R^\times$ and $b \in R$. Under such a transformation, we obtain
$$F(ax+b)=H(ax+b)^2+G(ax+b).$$
This is a decomposition of $F(ax+b)$ with $v(G(ax+b))\geq v(G)$, hence $\wbar(F(ax+b))\geq \wbar(F).$
Conversely, we can write $F(x)=F(a(a^{-1}(x-b))+b)$, which yields $\wbar(F)\geq \wbar(F(ax+b)).$ Therefore $\wbar(F(ax+b))=\wbar(F)$.

Next, consider the transformation $x \mapsto 1/x$. As before, the decomposition $F(1/x)=H(1/x)^2+G(1/x)$ is optimal, hence $\wbar(F)= \wbar(F(1/x)).$
Since $\PGL_2(R)$ is generated by these transformations, we conclude that $\wbar(X)$ does not depend on the choice of $x$. 

Now let $c\in R^\times$. As $K$ is algebraically closed, there exists $d\in R^\times$ with $d^2=c$.  The Laurent polynomial $cF$ admits a decomposition 
$$cF=(dH)^2+cG $$
with $v(cG)=v(G)$, hence $\wbar(cF)\geq \wbar(F)$. By symmetry, we get equality and the square defect of $X$ does not depend on the choice of  $F$. 
\end{Proof}
%\begin{Proof}
%	We have to show that $\wbar(F)$ is invarinat under $\PGL_2(R)$. For this, let  $F=G+H^2$ be an optimal decomposition. 
%	For all $a\in R^\times$ and $b\in R$ we get the decomposition
%	 $F(ax+b)=H(ax+b)^2+G(ax+b)$  of $F(ax+b)$ which satisfies $v(G(ax+b))\geq v(G)$. Hence $\wbar(F(ax+b))\geq \wbar(F)$. 
%	 As $F(x)=F(a(a^{-1}(x-b))+b)$, we get  $\wbar(F(ax+b))\leq \wbar(F)$ and hence equality holds. 
%	 Similarly, the decompsition $F(1/x)=H(1/x)^2+G(1/x)$ is optimal. As  $\PGL_2(R)$ is generated by these transformations, we are done. 
%\end{Proof}
\begin{Rem}{\rm
The notion of  the square defect of $X$ is compatible with the function $\wbar$ from \eqref{WBarDef} in the following way: Let $\CY_\lambda$ be the smooth model of $\BP^1_K$ with coordinate $u=2^\lambda x$, and let $Y_\lambda$ be its special fiber. Then $\wbar(\lambda)=\wbar(Y_\lambda)$. }
\end{Rem}
%\begin{Ex}\label{SquareDefectExample}
%%	\begin{enumerate}
%	 Let $F:=x^5+2x^4+x^2+x$ and consider any decomposition $F=H^2+G$ with $H\in R[x^{\pm 1}]$. As $H^2$ is a square, the valuation of its linear coefficient is at least $1$. Since the linear coefficient of  $F$ has valuation $0$, we get that $v(G)=0$. Hence $\wbar(X)=0$. 
%%		\item Let 
%%	\end{enumerate}
%\end{Ex}

From now on, additionally assume that $\CX$ is a marked model of $(\bar C, \bar P_1, \dots, \bar P_{2g+2})$. 
%It later turns out that it is often fruitful to track when the square defect of a component is trivial as in Example \ref{SquareDefectExample}.
\begin{Prop}\label{SquaredefectZeroProp}
	The square defect of $X$ is zero if and only if $X$ contains a marked point or an odd double point. 
\end{Prop}
\begin{Proof}
	First, assume that $X$ contains a marked point. By the proof of \cite[Prop. 4.1.1]{GP24}, there exists a coordinate $x$ such that $F=x\tilde F$ with $\tilde F \in R[x]$ with $\tilde F(0)\in R^\times$.  Now consider any decomposition  $F=G+H^2$ with $H\in R[x^{\pm 1}]$. As $H^2$ is a square, the valuation of its linear coefficient is at least $1$. Since the linear coefficient of  $F$ has valuation $0$, we get that $v(G)=0$. Hence $\wbar(X)=0$. 
	
	Next, assume that $X$ contains an odd double point. By \cite[Prop. 4.3.1]{GP24}, there exists a coordinate $x$, an element $a\in \Fm$ and $F_1, F_2\in R[x]$ with constant terms $1$ such that an equation for $C$ is given by  $z^2=F_1(x)F_2(a/x)$. In particular the constant coefficient of $F:=F_1(x)F_2(a/x)$ has valuation $0$ and by the same argument as in the first case, the square defect of $X$ is $0$. 
	
	Now, assume that $X$   contains neither marked points nor odd double points. This implies that the ramification points of $C$ reduce to $X$ in pairs. Choosing a coordinate $x$ such that no such pair reduces to $x=\infty$, there exists $\xi_1, \dots, \xi_{g+1}\in R$, possibly not all distinct, and $m_1, \dots, m_{g+1}\in \Fm$ such that $z^2=F$ with
	$$F:=\prod_{i=1}^{g+1}(x-\xi_i)(x-\xi_i-m_i) $$
	is a hyperelliptic equation of $C$. Now define $H:=\prod_{i=1}^{g+1}(x-\xi_i)$ and $G:=F-H^2$. Then  $F=G+H^2$ is a decomposition with $v(G)>0$, thus $\wbar(X)>0$.
\end{Proof}

\begin{Prop}\label{OnlyOneDoublePointProp}
	Assume that $\CX=\hat \CC$ and $\wbar(X)<2$.
	Let
	$\hat p \in X$ be a double point. Then the curve $C_0$ possess a unique double point over $\hat p$. 
\end{Prop}
\begin{Proof}
		If $\hat p$ is an odd double point, the proposition follows from \cite[Prop. 4.4.1]{GP24}. Otherwise let $\alpha>0$ be the thickness of $\bar p$ and identify a neighborhood of $\bar p$ with an open subscheme of $\op{Spec}(R[x,y]/(xy-2^\alpha))$ for $\alpha>0$. Then by \cite[Prop. 4.3.1]{GP24} there exists a Laurent polynomial $F\in R[x^{\pm 1}]$ with $v(F)=0$ such that the function field of $C$ is $K(x,z)$ with $z^2=F(x)$. Moreover, we have $w(F)=\wbar(X)<2$ by Theorem-Definition \ref{Square-defect}. 
		This implies that for the function $\wbar$ from \eqref{WBarDef} we have $\wbar(0)=\wbar(X)<2$. Hence $\wbar$ has a non-horizontal segment by \cite[Prop. 3.4.4]{GP24}. As $\CX=\hat \CC$, there are no components of type (b) above $\hat p$ and thus $\wbar$ is linear by \cite[Prop. 4.5.1]{GP24}. So $\wbar$ is linear with nonzero slope, which implies that the curve $C_0$ possess a unique double point over $\hat p$ by \cite[Prop. 4.5.12]{GP24}. 
%, which implies the proposition by \cite[Prop 3.3.2 (b)]{GP24}. 
%
\end{Proof}

The following Proposition was first proved by Raynaud in \cite[Th. 2, Prop. 2 (iii)]{Raynaud1990} and \cite[Lem. 3.1.2]{Raynaud1999} under the additional assumption that $\bar C_0$ is smooth, although his proof works in full generality. We present a short proof using our methods. 
\begin{Prop}\label{RaynaudProp}
	Let $\bar p\in \bar C_0$ be a smooth point. Then the dual graph of the inverse image of $\bar p$ in $C_0$ is a tree. Moreover,
	let $\bbar T$ be an irreducible component of $\hat C_0$ above~$\bar p$. Then its inverse image $T$ in $C_0$ is irreducible and 
	\begin{itemize}
		\item either $\bbar T$ is of type (d) and $T$ has genus $>0$ and $2$-rank $0$ and is separable over~$\bbar T$,
		\item or $\bbar T$ is of type (c) and $T$ has genus $0$ and is purely inseparable over~$\bbar T$.
	\end{itemize}
\end{Prop}
\begin{Proof}
	We identify a neighborhood of $\bar p$ with an open subscheme of $\op{Spec}(R[x])$. Then there exists a separable polynomial $F\in R[x]$ with $v(F)=0$ such that the function field of $C$ is $K(x,z)$ with $z^2=F(x)$. By \cite[Prop. 4.1.1, Prop. 3.3.2]{GP24}, the existence of $\hat T$ implies that $\bar p$ does not lie in the zero locus of  $F$ and that $\gamma:=w(F)<2$. Moreover the component $\hat T$  is given by a coordinate $y$ with $x=2^\epsilon y+ \xi_0$ for a rational $\epsilon>0$ and $\xi_0\in R$. After replacing $x$ by $x-\xi_0$, the coordinate has the form $x= 2^\epsilon y$. Moreover, the inverse image $T$ of $\hat T$ in $C_0$ does not contain a marked point.

	First, assume $w(F(2^\epsilon y))<2$ and let $\hat p\in \hat T$ be a double point. By Proposition \ref{OnlyOneDoublePointProp} the curve $C_0$ possesses a unique double point over $\hat p$. 
	Moreover, the component $T$ has genus $0$ by \cite[Prop. 3.3.2]{GP24} and is purely inseparable over $\hat T$.
	If $\hat T$ is of type (d) this implies that $T$ is smooth, thus blowing it down yields a semistable marked model of $C$, contradicting the minimality of $\hat \CC$. 
	
	Second, assume that $w(F(2^\epsilon y))=2$ and let $\hat p$ be the double point defined by $y=0$. Then $T$ is smooth outside the preimage of $\hat p$ by \cite[Prop. 3.3.2 (a)]{GP24}. Hence $\hat T$ is a component of type (d).  
	Let $\alpha$ be the thickness of $\hat p$. Consider the function $\wbar$ for the polynomial $F(2^{\epsilon-\alpha}y)$. As there are no components of type (b) over $\hat p$, the function $\wbar$ does not have any break points by \cite[Prop. 4.5.1]{GP24}. 	Since  $F$ is a polynomial, we have $w(F^\lambda y))<2$ for $\lambda\in[0,\epsilon)\cap \BQ$, which implies that $\wbar$ cannot be constant. Hence $\wbar$ is linear with nonzero slope, and  thus $T$ has a single double point by \cite[Prop. 4.5.12]{GP24}. 
	As $C_0$ is stable, the component $T$ must have genus $>0$. Finally, as no Weierstrass points of $C$ reduce to $T$, it has $2$-rank $0$. 
	
	Third, assume that $w(F(2^\epsilon y))>2$. Then $T$ is the disjoint union of two irreducible components of genus $0$ by \cite[Prop. 3.3.2]{GP24}, which are both leaf components. Blowing these down again yields a semistable marked model of $C$, contradicting the minimality of $\hat \CC$. 
	
	These three cases show that for every double point in the preimage of $\bar p$ in $\hat C_0$, the curve $C_0$ possesses a unique double point, which together with the facts about $T$ implies that  the inverse image of $\bar p$ in $C_0$ is a tree.
\end{Proof}

\subsection{Square defect function at double points} \label{SubSectionSDF}
Let $\bar p$ be an even double point of~$\bar C_0$ and identify a neighborhood $\bar\CU \subset \bar \CC$ with an open subscheme of $\Spec R[x,y]/(xy-a)$ for some nonzero $a\in\Fm$, such that $\bar p$ is given by $x=y=0$. After multiplying $a$ by a unit we may assume that $a=2^\alpha$ for some $\alpha>0$, where $\alpha$ is the thickness of~$\bar p$.
Identifying $y$ with $\frac{2^\alpha}{x}$, the ring in question is isomorphic to the subring $R[x,\tfrac{2^\alpha}{x}]$ of the ring of Laurent polynomials $K[x^{\pm1}]$. 
Moreover, there exists $F  \in R[x,\frac{2^\alpha}{x}]$ with constant coefficient $1$ such that function field of $C$ is $K(x,z)$ with $z^2=F (x)$. Let $X$ be the irreducible component of $\bar C_0$ for which $x$ is a coordinate. Similarly, let $Y$ be the component for which $y$ is a coordinate. 

\begin{ThmDef}\label{SquareDefectFunctionTheoremDefinition}
	The square defect function at $\bar p$ with respect to $X$ is defined as 
	\UseTheoremCounterForNextEquation
\begin{equation}\label{SquareDefectFunctionDef}
	\wbar_{\bar p, X}\colon\ \BQ\,\cap\, [0, \alpha] \longto\BR,\ \ \lambda\mapsto 
	\wbar_{\bar p, X}(\lambda) := \wbar(F(2^\lambda u)).
	% \min \{ 2, w(F(2^\lambda u)) \}.
\end{equation}
It is independent of the choice of $x$ and  $F$. Moreover, for all $\lambda \in \BQ \cap [0,\alpha]$ we have $\wbar_{\bar p,X}(\lambda)=\wbar_{\bar p, Y}(\alpha-\lambda)$. 
\end{ThmDef}
\begin{Proof} By preceding definitions, we have $\wbar_{\bar p, X}(0)=\wbar(X)$ and $\wbar_{\bar p, X}(\alpha)=\wbar(Y).$ Moreover, by \cite[Prop. 4.5.1]{GP24}, the substitution $x=2^\lambda u$ yields an irreducible component $X_\lambda$ of $\hat C_0$ if and only if $\lambda$ is a break point of $\wbar_{\bar p,X}$. Here $\wbar_{\bar p, X}(\lambda)=\wbar(X_\lambda)$ and the $\lambda$ is determined by the thicknesses of the double points of $\hat C_0$ over $\bar p$. 
As $\wbar_{\bar p, X}$ is a piecewise linear function it is fully determined by these values, which are all independent of the choice of  $F$ and $x$. Hence $\wbar_{\bar p, X}$ is as well. The last statement follows from repeating the computation with interchanged roles of $x$ and $y$. 
\end{Proof}

\begin{Rem}{\rm
All properties previously established for $\wbar$ in \cite{GP24} are inherited by $\wbar_{\bar p,X}$. We will apply this  throughout the remainder of this article without explicit mention.}
\end{Rem}

%\begin{Prop}
%	Let $\CX$ be the 
%\end{Prop}

\begin{Prop}\label{WbarinCohatProp} Let $\hat X$ be the proper transform of $X$ in $\hat C_0$ and let $\hat p$ be the unique double point over $\bar p$ contained in $X$. Let $\beta$ be the thickness of $\hat p$. Then $\beta\leq \alpha$ and we have 
	$$\wbar_{\hat p, \hat X}(\lambda) \ = \ \wbar_{\bar p, X}(\lambda) $$
	for all $\lambda \in \BQ \cap [0,\beta]$.   
\end{Prop}
\begin{Proof}
First, assume that there are no components of type (b) above $\bar p$. Then $\hat \CC \to \bar \CC$ restricts to an isomorphism in a neighborhood of $\bar p$, from which the claim follows. 

Second, assume that there are components of type (b) above $\bar p$ and let $\hat T$ be the unique such component intersecting $\hat X$. Then $\hat T$ has coordinate $t=x/2^{\lambda_1}$ for some $\lambda_1>0$ and $\hat p$ is a double point of thickness $\beta=\alpha-\lambda_1$.  Moreover, there exists an open neighborhood of $\hat p$ isomorphic to an open subscheme of $R[x,t]/(xt-2^{\beta}).$ Furthermore,  note that we have $F\in R[x,\frac{2^\beta}{x}]$. For $\lambda \in \BQ \cap [0, \beta]$, we get
$$\wbar_{\hat p, \hat X}(\lambda)\ =\ \wbar(F(2^\lambda u)) \ = \ \wbar_{\bar p, X}(\lambda), $$
as desired. 
\end{Proof}

\begin{Cor}\label{UniqueDoublePointOverDoublePoint}
	%	ToDo
	The preimage of $X$ in $C_0$ has a unique double point above $\bar p$ if and only if the slope of $\wbar_{\bar p, X}$ to the right of $0$ is non-zero. Otherwise, it contains two double points above $\bar p$. 
\end{Cor}
\begin{Proof} Let $\hat X$ be the proper transform of $X$ in $\hat C_0$ and let $\hat p$ the unique double point of $\hat X$ over $\bar p$. 
	 By Proposition \ref{WbarinCohatProp}, the unique slope of $\wbar_{\hat p, \bar X}$ is equal to the slope of $\wbar_{\bar p, X}$ to the right of $0$. 
	 
If the slope of $\wbar_{\bar p, X}$ to the right of $0$ is non-zero, the model $\CC$ has a unique double point above $\hat p$, which is then the only double point in the preimage of $X$ in $C_0$ above $\bar p$ 	 by \cite[Prop. 4.5.12]{GP24}. 

If the slope of $\wbar_{\bar p, X}$ to the right of $0$ is $0$, there are two double points in $\CC$ above $\hat p$  by \cite[Prop. 4.5.10]{GP24}. 
\end{Proof}

%ToDo: Draw picture of special fiber of $\CX$. 
\begin{Prop}\label{WbarExtraComponentProp} 
		Let $ \CX$ be the semistable model of $\bar C$ obtained from $\bar \CC$ by adjoining an irreducible component $\hat T$ with coordinate  $t=x/2^{\lambda_1}$ for $\lambda_1\in \BQ \cap (0,\alpha)$. Let $\hat p$ be the unique double point of $\hat T$ above $\bar p$ and  contained in the proper transform of $Y$. Denote the  thickness of $\hat p$ by $\beta$. Then $\beta=\alpha-\lambda_1$ and we have
			$$\wbar_{\hat p, \hat T}(\lambda) \ = \ \wbar_{\bar p, X}(\lambda+\lambda_1) $$
		for all $\lambda \in \BQ \cap [0,\beta]$.   
\end{Prop}
\begin{Proof} 
	Let $\hat Y$ be the proper transform of $Y$ in $\bar \CX$. 
As $t$ is a coordinate for $\hat T$ and $y$ is a coordinate for $\hat Y$ and we have $ty=2^{\alpha-\lambda_1}$ and the thickness  of $\hat p$ is $\beta=\alpha-\lambda_1$.
Moreover, there exists an open neighborhood of $\hat p$ isomorphic to an open subscheme of $R[t,y]/(ty-2^{\beta}).$ Furthermore, for $F_t:=F(2^{\lambda_1} t)$ we have $F_t\in R[t,\frac{2^\beta}{t}]$. For $\lambda \in \BQ \cap [0, \beta]$, we compute
$$\wbar_{\hat p, \hat T}(\lambda)\ =\ \wbar(F_t(2^\lambda u)) \ = \ \wbar(F_t(2^{\lambda+\lambda_1})) \ = \ \wbar_{\bar p, X}(\lambda + \lambda_1), $$
as desired. 
\end{Proof}

%\begin{Prop}
%	Let $X$ be an irreducible component of $\bar C_0$ and let $\bar p\in \bar X$ be an ordinary double point. Assume that the morphism $\hat \CC \to \bar \CC$ restricts to an isomorphism in a neighborhood of $\bar p$. 
%\end{Prop}

%%%%%%%%%%%%%%%%%%%%%%%%%%%%%%%%%%%%%%%%%%
%%%%%%%%%%%%%%%%%%%%%%%%%%%%%%%%%%%%%%%%%%
\section{Approximation}\label{ApproximationSection}

\subsection{Finding optimal decompositions}\label{SubSectionOptimalDecompositions}
%We use the setting and notation of ToDo: 
%More specific \cite{GP24}. 
Let $\alpha>0$ and  fix Laurent polynomials $F,\tilde F \in R[x,\frac{2^\alpha}{x}]$. Fix optimal decompositions  $F=G+H^2$ and $\tilde F=\tilde G+ \tilde H^2$.  
%Note that this does not require the ramification points to specialize to different points. 

\begin{Lemma}\label{FirstApproximationLemma}
	Assume that $v(F-\tilde F) >\wbar(F)$ or $v(F-\tilde F)\geq 2$. 
	Then   $\wbar(F)=\wbar(\tilde F)$.
\end{Lemma}
\begin{Proof}
As $F=\tilde H^2 +(F-\tilde H^2)$ is a decomposition of $F$, we have 
	\UseTheoremCounterForNextEquation
	\begin{equation} \wbar(F)\ \geq  \ \min\{2,v(F-\tilde H^2)\}\label{FirstApproximationLemmaFirstEq}\end{equation}
by Definition \ref{WBarDef}.
	Moreover, we compute 
	$$F-\tilde H^2=(F-\tilde F) + (\tilde F-\tilde H^2)=(F-\tilde F) + \tilde G,$$  which  yields 
	$v(F-\tilde H^2)\geq \min\{v(F-\tilde F), v(\tilde G)\}$.
	Hence 
	\UseTheoremCounterForNextEquation
	\begin{equation}\min\{2,v(F-\tilde H^2)\}\geq \min \{2, v(F-\tilde F), v(\tilde G)\} \label{FirstApproximationLemmaExtraEq}.\end{equation}
	As $\tilde F=\tilde G+ \tilde H^2$  is optimal,  we have $\wbar(\tilde F) = \min\{2, v(\tilde F-\tilde H^2)\}$, which together with Equation \eqref{FirstApproximationLemmaExtraEq} implies 
		\UseTheoremCounterForNextEquation
	\begin{equation}  \min\{2,v(F-\tilde H^2)\} \ \geq \ \min\{v(F-\tilde F), \wbar(\tilde F) \}. \label{FirstApproximationLemmaSecondEq}\end{equation}
	Using that $v(F-\tilde F)>\wbar(F)$ or $v(F - \tilde F)\geq 2$ and combining equations \eqref{FirstApproximationLemmaFirstEq} and \eqref{FirstApproximationLemmaSecondEq} yields 
	$\bar w(F) \geq \bar w(\tilde F)$.
	
	 If $v(F-\tilde F)\geq 2$, exchanging  $F$ and $\tilde F$ yields the other inequality.
	 
	 Otherwise, we have  $v(F-\tilde F)>\wbar(F)\geq \wbar(\tilde F)$ and exchanging  $F$ and $\tilde F$ again yields the other inequality. 
%	
%	 the second inequality is derived writing $F-\tilde F + \tilde F-\tilde H^2$. 
%	In both cases, this yields $\wbar(F)\geq \wbar(\tilde F)$. Exchanging the roles of  $F$ and $\tilde F$ yields the other inequality. 
\end{Proof}

\begin{Lem}\label{ApproximateDecompositionIsGoodEnoughLemma}
	Assume that $v(F-\tilde F) > \wbar(F)$. Then  $\tilde F= H^2 +(\tilde F- H^2)$ is an optimal decomposition with $[(\tilde F- H^2)/2^{\wbar(F)}]=[G/2^{\wbar(F)}]$.
\end{Lem}
\begin{Proof}
	We have 
	$$\tilde F -H^2 -G\ =\ (\tilde F -F)+ (F-H^2+G) \ = \ \tilde F -F$$ which yields
	$$v(\tilde F -H^2 -G) \ = \  v(\tilde F- F) \ > \ \wbar(F).$$
	The last equation implies  $[(\tilde F- H^2)/2^{\wbar(F)}]=[G/2^{\wbar(F)}]$. As $\wbar{F}=\wbar(\tilde F)$ by  Lemma \ref{FirstApproximationLemma}, we get that $\tilde F= H^2 +(\tilde F- H^2)$ is an optimal decomposition.
\end{Proof}

\begin{Ex}{\rm
	Let $F=2^{1/2}x^5+x^{4}+2 x^{3}+3 x^{2}+2 x$ and $\tilde F =2^{1/2}x^5+x^{4}+2 x^{3}+5 x^{2}+2 x$. Then $ F =(x^2+x+1)^2+2^{1/2}x^5$ is an odd decomposition, hence $\wbar(F)=1/2.$ As $v(F-\tilde F)=v(2 x^2)=1$, we have $\wbar(\tilde F)=1/2.$ }
\end{Ex}

\subsection{Changing  $F$ near a smooth point}\label{SubSectionSmoothPointApproximation}
%ToDo: Change from $\bar \CC$ to $\CX$ where needed. 
%In the following, let $ \CX$ be a semistable model of $\bar C$ and fix a closed smooth point $\bar p$ of the special fiber $X_0$ of $\CX$.
%Choose a coordinate $x$ of $\bar \CC$ such that $\bar p$ is defined by $x=0$ and let $ \CU$ be an open neighborhood of $\bar p$ isomorphic to an open subscheme of $\Spec R[x^{}]$. 
Identify a neighborhood $\bar\CU \subset \bar\CC$ with an open subscheme of $\Spec R[x^{}]$.
%Let $U_0$ be the special fiber of $\bar \CU$.
Furthermore, let $F\in R[x]$ with $v(F)=0$ be such that an equation for $C$ is given by $z^2=F$. Now fix $\tilde F\in R[x^{}]$ %with $v(F-\tilde F)>2$ 
%the field $K(x,\tilde z)$ with $\tilde z^2=\tilde F$. % 
and let $\tilde C$ be the curve  defined by $z^2=\tilde F(x)$. 
Assume that $v(F-\tilde F)>\wbar (F)=:\gamma$ and let $\bar U_0$ be the reduction of $\bar \CU$ in $\bar C_0$. Let  $F=G+H^2$ be an optimal decomposition.

%%%%Old description:
% OLD:
%In the following, let $\bar \CC$ be any marked semistable model of $\bar C$ that dominates its stable marked model. 
%From now on, we use the notation of Sections 3.6 and 4.2 of \cite{GP24}. Let $\bar p$ be a smooth unmarked point in the special fiber of any unmarked point. 
%%
%%we assume that $\bar p$ is a smooth point of~$\bar C_0$ that does not lie in any of the sections~$\bar\CP_i$.
%We identify a neighborhood $\bar\CU \subset \bar\CC$ with an open subscheme of $\Spec R[x^{}]$.
%%such that $\bar p$ is defined by $x=0$. 
%Suppose that the function field of $C$ is $K(x,z)$ with $z^2=F (x)$ for a separable polynomial $F  \in R[x]$ with $v(F )=0$. Then by assumption $\bar p$ does not lie in the zero locus of~ $F$.  
%%Consider any $\xi_0\in R$ with $F(\xi_0)\in R^\times$. This is equivalent to $ \tilde F(\xi_0)\in R^\times$.  Take a new variable~$y$. Then any component of type (d) over $\bar p$ arises with the substitution \dots. 
%Now fix $\tilde F\in R[x^{}]$ %with $v(F-\tilde F)>2$ 
%%the field $K(x,\tilde z)$ with $\tilde z^2=\tilde F$. % 
%and let $\tilde C$ be the curve  defined by $z^2=\tilde F(x)$. 
% Assume that $v(F-\tilde F)>\wbar (f):=\gamma$. 
%%If there are components of type (d) over $\bar p$, assume that  $v(F-\tilde F)>2$ .

%\begin{Proposition} \label{stableReductionOverSmoothUnmarkedPoint}
%	The stable reduction of $C$ over $\bar p$ is isomorphic to the stable reduction of $\tilde C$ over $\bar p$. 
%\end{Proposition}
\begin{Prop}\label{BareNormalizationSpecialFiberIsIsomorphic}
	The preimage of $\bar U_0$ in the normalization of $ \bar \CU$ in $K(C)$ is isomorphic to the preimage of $\bar U_0$ in the normalization of $ \bar \CU$ in $K(\tilde C)$. 
\end{Prop}
\begin{Proof} 
%	Let  $F=G+H^2$ be an optimal decomposition of  $F$.
	 Then $\tilde F=(\tilde F-H^2)+H^2$ is an optimal decomposition by Lemma \ref{ApproximateDecompositionIsGoodEnoughLemma}.
	By  \cite[Prop. 3.3.2]{GP24}, the normalization of $\CU$ in $K(C)$ is isomorphic to $R[x][t_1]$ with
	\UseTheoremCounterForNextEquation
	\begin{equation}2^{1-\gamma/2}Ht_1+t_1^2=G/2^\gamma.\label{FirstNormalizationEquation}\end{equation}
	Similarly, the normalization of $\CU$ in $K(\tilde C)$ is isomorphic to $R[x][\tilde t]$ with 
	\UseTheoremCounterForNextEquation
		\begin{equation}2^{1-\gamma/2}Ht_2+t_2^2=(\tilde F -H^2)/2^\gamma.\label{SecondNormalizationEquation}\end{equation}
		By Lemma  \ref{ApproximateDecompositionIsGoodEnoughLemma}, we get $[(\tilde F- H^2)/2^{\gamma}]=[G/2^{\gamma}]$. This implies that the reductions modulo $\Fm$ of the equations \eqref{FirstNormalizationEquation} and \eqref{SecondNormalizationEquation} have the same form and hence the special fibers of the normalization are isomorphic. 
\end{Proof}
From now on, additionally assume that  $v(F-\tilde F)>2$. Let $\tilde \CU$ be the minimal semistable model of the generic fiber of $\bar \CU$ dominating $\bar \CU$ such that the normalization of $\tilde \CU$ in $K(\tilde C)$ is semistable. 
% and that $\CX$ is a marked model of $(\bar C, \bar P_1, \dots, \bar P_{2g+2})$. %ToDo: Maybe strengthen to $\geq 2$. 
%Let $\hat \CX_1$ be the minimal semistable marked model of $(\bar C, \bar P_1, \dots, \bar P_{2g+2})$ dominating $\CX$ such that the normalization in $K(C)$ above the preimage of $\bar p$ is semistable. 
%The fact that $\CX$ is a marked model of $(\bar C, \bar P_1, \dots, \bar P_{2g+2})$ implies that $\CX_1$ has only components of type (b) and (c) 
%Similarly, let $\hat \CX_2$ be the minimal semistable marked model of $(\bar C, \bar P_1, \dots, \bar P_{2g+2})$ dominating $\CX$ such that the normalization in $K(\tilde C)$ above the preimage of $\bar p$ is semistable. 
%Moreover, let $\hat \CC\to \hat \CC' $ be the contraction of all irreducible components of the special fiber that are not  

\begin{Proposition} \label{stableReductionOverSmoothUnmarkedPoint}	
The preimage of $\bar \CU$ in $\hat \CC$ is equal to $\tilde \CU.$ 
	Moreover, the preimages of the components of the special fiber of $\tilde \CU$ of type (c) and (d)  in  the normalizations in $K(\tilde C)$ are isomorphic to their preimages in $\CC$.

%The normalization of $\hat \CC$ in $K(\tilde C)$ is semistable above $\bar p$. Moreover, locally $\hat \CC$ is the minimal model dominating $\bar \CC$ with that property. 
% For any model $\hat \CC'$ that is the contraction of  irreducible components above $\bar p$, the normalization of $\hat \CC'$ in $K(\tilde C)$ is not semistable.
%   Moreover, the preimages of the components of $\hat C_0$ above $\bar p$ in the normalizations in $K(C)$, respectively $K(\tilde C)$ are isomorphic. 
%
%ToDo: Maybe split this into two proposition, one for no components of type (d). 
%	 
%	
%	The stable reduction of $C$ over $\bar p$ is isomorphic to the stable reduction of $\tilde C$ over $\bar p$. 
%	%	The irreducible components of type (d) over $\bar p$ are in $1:1$ correspondence to the components of type (d) when we replace   $F$ by $\tilde F$. Under this correspondence, the substitutions yielding them are the same. Moreover, the special fibers of the models are isomorphic. 
\end{Proposition}

\begin{Proof} 
%		Let $F=H^2+G$ be an optimal decomposition. 
	By Lemma \ref{ApproximateDecompositionIsGoodEnoughLemma} the decomposition $\tilde F=(\tilde F-H^2)+H^2$ is optimal  and we have   $[(\tilde F- H^2)/2^{\gamma}]=[G/2^{\gamma}]$.  %As $\CX$ is a marked model of $(\bar C, \bar P_1, \dots, \bar P_{2g+2})$, there are only components of type (c) and (d) above $\bar p$ 
		
		First, assume that there are no components of type (d) above $\bar \CU$. Then $\hat \CC\to \bar \CC$ is an isomorphism above $\bar \CU$ and our claim follows directly from  Proposition \ref{BareNormalizationSpecialFiberIsIsomorphic}.
		
		Second, assume that there are components of type (d) above $\bar \CU$. Then $\gamma<2$ and
%		 and by Proposition 3.3.2 in \cite{GP24}, the normalization of $R[x]$ in $K(C)$ modulo $\Fm$ has the form $t^2=[g/2^\gamma]$. 
	each component of type (c) or (d) over $\bar p$ is given by a coordinate $y$ with $x=\xi_0+2^\epsilon y$, where $\xi_0$ is an element of $R$ and $\epsilon>0$ is a rational number. 
		Then $$v(F(\xi_0+2^\epsilon y)-\tilde F(\xi_0+2^\epsilon y))\geq v(F-\tilde F)>2.$$ Let $\BP^1(y)$ be the smooth model of $\BP^1_K$ with coordinate $y$. By Proposition 3.3.2 in \cite{GP24} and Lemma \ref{ApproximateDecompositionIsGoodEnoughLemma}, the special fiber of the normalization of $\BP^1(y)$ in $K(C)$ is isomorphic to the special fiber of its normalization in $K(\tilde C)$. In particular, a component of type (d) appears in $\tilde \CU$ if and only if its normalization in $K(C)$ has positive genus. As the components of type (c) appearing in $\tilde \CU$ are implied by this datum, our claim follows. 
\end{Proof}
%\newpage 
\begin{Ex}
\begin{enumerate}[(a)]\label{ExampleApproxSmoothPoint}
	\item
	Let $a\in R$ and $L\in R[x]$ whose constant coefficient is a unit. Define $\tilde F:=F(x)\cdot L(a/x)$. If $v(a)>\gamma$, the Laurent polynomial $\tilde F$ satisfies the prerequisites of Proposition \ref{BareNormalizationSpecialFiberIsIsomorphic}. If additionally $v(a)>2$, the prerequisites of Proposition \ref{stableReductionOverSmoothUnmarkedPoint} are satisfied as well. 
	\item Let 
	$$F=128x^7 + 24x^6 + 257x^5 + 176x^4 + 154x^3 + 25x^2 + x$$ and assume that $\bar \CU$ is isomorphic to $\Spec R[x]$. 
	Set $\tilde F:=x^5 + 2x^3 + x^2 + x$.  Applying the genus 2 classification from \cite{GP24}, we get that $\tilde F$ defines a genus $2$ curve of type (A3), hence by Proposition \ref{stableReductionOverSmoothUnmarkedPoint} there is one component of type (d) above a smooth point of $\bar \CU$ whose preimage in the normalization in $K(C)$ has genus $2$. 
%	
%	Then 
% $F\equiv \tilde F \mod 8$, hence by Proposition \ref{stableReductionOverSmoothUnmarkedPoint} the normalization 
%  Applying the genus 2 classification from \cite{GP24}, we get that $\tilde F$ defines a curve of type (A3), hence there is one component of type (d) above a smooth point of $X$ whose preimage in the normalization in $K(C)$ has genus $2$. 
\end{enumerate}
	\end{Ex}

\subsection{Changing  $F$ near a double point}\label{SubSectionDoublePointApprox}
%Recall that $\bar \CC$ is a semistable marked model of $(\bar C, \bar P_1, \dots, \bar P_{2g+2})$. 
Consider an ordinary double point $\bar p$ of the special fiber $\bar C_0$ of $\bar \CC$.
We identify a neighborhood $\bar{\mathcal{U}} \subset \bar{\mathcal{C}}$ with an open subscheme of Spec $R[x, y] /(x y-a)$ for some nonzero $a \in \mathfrak{m}$, such that $\bar{p}$ is given by $x=y=0$. After multiplying $a$ by a unit we may assume that $a=2^\alpha$ for some $\alpha>0$, where $\alpha$ is the thickness of $\bar{p}$. Identifying $y$ with $\frac{2^\alpha}{x}$, the ring in question is isomorphic to the subring $R\left[x, \frac{2^\alpha}{x}\right]$ of the ring of Laurent polynomials $K\left[x^{ \pm 1}\right]$.

%We use the notation from Sections 3.5 and 4.5 of \cite{GP24}. 

%Old:
%We use the notation from Sections 3.5 and 4.5 of \cite{GP24}. Moreover, we specify for which polynomial $p\in R[x^{\pm1}]$ we compute $\wbar$ by writing $\wbar_p$. 
%Let $F\in R[x, \frac{2^\alpha}{x}]$ be as in Proposition 3.5.1 of [GP24]. 
%%\begin{Cor} Assume that for $\lambda \in [0,\alpha]\cap \BQ$, we have $v(F(2^\lambda x)-\tilde F(2^\lambda x)) > \wbar_F(\lambda)$ or $v(F(2^\lambda x)-\tilde F(2^\lambda x))\geq 2$.  Then $\wbar_{f}(\lambda)=\wbar_{\tilde F}(\lambda)$. 
%%\end{Cor}
%%\begin{Proof}
%%	Apply Lemma  \ref{FirstApproximationLemma} to $F(2^\lambda x)$ and $\tilde F(2^\lambda x)$.
%%\end{Proof}
%Moreover, let $\bar p\in \bar C_0$ be an even double point. Assume that the function field of $C$ is $K(x,z)$ with $z^2 = f$.  Fix $\tilde F$ with $v(F(2^\lambda x)-\tilde F(2^\lambda x)) > \wbar_F(\lambda)$ for all $\lambda \in [0,\alpha]\cap \BQ$. Let  $F=G+H^2$ and $\tilde F=\tilde H^2 +\tilde G$ be optimal decompositions.
%Finally, let $\tilde C$ be the curve  defined by $ z^2=\tilde F(x)$. ToDo: Complete reformulation. 
\medskip 
By Proposition 4.3.1 in \cite{GP24}, there exist $F_1\in R[x]$ and $F_2\in R[y]$ with constant coefficients $1$ such that an equation for $C$ is given by $z^2=F$, where
$$F\ :=\ \biggl\{\!\!\begin{array}{cl}
	F_1(x)F_2(\tfrac{2^\alpha}{x}) &\hbox{if $\bar p$ is even,}\\[3pt]
	xF_1(x)F_2(\tfrac{2^\alpha}{x}) &\hbox{if $\bar p$ is odd.}
\end{array}$$
Consider $\tilde F\in R\left[x, \frac{2^\alpha}{x}\right]$ with $v(F-\tilde F)>0$ and let $\tilde C$ be the curve defined by $\tilde z^2=\tilde F$. 

\medskip 
\begin{Prop}\label{Odddoublepointapprox}
	Assume that $\bar p$ is an odd double point in $C_0$.  The normalization of $\bar \CU$ in the function field  $K(\tilde C)$ is semistable and possesses a unique double point $p$ over $\bar p$ with thickness $\alpha/2$. 
\end{Prop}
\begin{Proof}
	As $v(F-\tilde F)>0$, we have $\frac{\tilde F}{x}\in R^\times +R[x]x+R[2^{\alpha}/x]$. By  \cite[Prop. 3.5.1]{GP24}, there exist $\tilde F_1\in R[x]$ and $\tilde F_2\in R[y]$ with constant coefficients $1$ and $c\in R^\times$ such that $$\tilde F= cx\tilde F_1(x) \tilde F_2(\tfrac{2^\alpha}{x}).$$ 
	By  \cite[Prop. 4.4.1]{GP24}, the claim follows from this factorization of $\tilde F$. 
\end{Proof}
\medskip
From now on, we assume that $\bar p$ is an even double point. Let $X$ be the irreducible component of $\bar C_0$ with coordinate $x$. 
%
%For a new variable $u$, a Laurent polynomial $L \in R[x,\frac{2^\alpha}{x}]$ and $\lambda\in[0,\alpha]\cap \BQ$, define 
%$$\wbar_{L}(\lambda):=\wbar(L(2^\lambda u)).$$ This extends the notion \eqref{WBarDef}, and by \cite[Prop. 3.4.4]{GP24}, the function $\wbar_f$ is continuous and piecewise linear concave. Moreover by \cite[Prop. 4.5.1 ]{GP24} the substitution $x=2^\lambda u$ defines a component of type (b) over $\bar p$ if and only if $\lambda$ is a break point of $\wbar_f$.  
Assume that for $\lambda \in [0,\alpha]\cap \BQ$ and a new variable $u$, we have 
$$v(F(2^\lambda u)-\tilde F(2^\lambda u)) > \wbar_{\bar p, X }(\lambda) \; \; \text{ or } \; \; v(F(2^\lambda u)-\tilde F(2^\lambda u))\geq 2.$$

%ToDo: Mention how to compute $\wbar$. 
%Second, we assume that $\bar p$ is an odd double point. Then there exists $F_1\in R[x]$ and $F_2\in R[y]$ with constant coefficients $1$ such that an equation for $C$ is given by $z^2=F_1(x)F_2(2^\alpha(x))$. 
\begin{Prop}\label{OnlyComponentsOfTypeBApproxProp}
	 For $\lambda \in [0,\alpha]\cap \BQ$, we have
	 \UseTheoremCounterForNextEquation
	 \begin{equation}\wbar_{\bar p, X}(\lambda)\ =\ \wbar(F( 2^\lambda u))\ = \ \wbar(\tilde F( 2^\lambda u)).\label{wbarapproxEx}\end{equation}
% Moreover,  the components of type (b) above $\bar p$ and their preimages in the normalizations in $K(C)$, respectively $K(\tilde C)$ are isomorphic. 
\end{Prop}
\begin{Proof}
The first equality in Equation \ref{wbarapproxEx} is true by definition of $\wbar_{\bar p, X}$. To get the second equality, apply  Lemma  \ref{FirstApproximationLemma} to $F(2^\lambda x)$ and $\tilde F(2^\lambda x)$ for all $\lambda \in [0,\alpha]\cap \BQ$.
%
% Now let  $F=G+H^2$ be an optimal decomposition. By Lemma \ref{ApproximateDecompositionIsGoodEnoughLemma} the decomposition $\tilde F= H^2 +(\tilde F- H^2)$ is optimal  with $[(\tilde F- H^2)/2^{\wbar(F)}]=[g/2^{\wbar(F)}]$. yields equations for the preimages of the components of type (b) that are equal modulo $\Fm$ to those obtained using the decomposition $F=H^2+G$. 
\end{Proof}
Additionally, assume that  $v(F(2^\lambda u)-\tilde F(2^\lambda u)) > 2$ for all $\lambda \in [0,\alpha]\cap \BQ$. 
 Let $\tilde \CU$ be the minimal semistable model of the generic fiber of $\bar \CU$ dominating $\bar \CU$ such that the normalization of $\tilde \CU$ in $K(\tilde C)$ is semi-stable. 
%Let $\hat \CX_1$ be the minimal marked model of $(\bar C, \bar P_1, \dots, \bar P_{2g+2})$ dominating $\bar \CC$ such that the normalization in $K(C)$ above the preimage of $\bar p$ is semistable. 
%Similarly, let $\hat \CX_2$ be the minimal marked model of $(\bar C, \bar P_1, \dots, \bar P_{2g+2})$ dominating $\bar \CC$ such that the normalization in $K(\tilde C)$ above the preimage of $\bar p$ is semistable. 
%Moreover, let $\hat \CC\to \hat \CC' $ be the contraction of all irreducible components of the special fiber that are not  
\begin{Proposition} \label{ComponensofTypeCDApproxProp}	
	The preimage of $\bar \CU$ in $\hat \CC$ is isomorphic to $\tilde \CU.$ 
	Moreover, the preimages of the components of the special fiber of $\tilde \CU$ of type (b), (c) and (d) above $\bar p$ in  the normalizations in $K(\tilde C)$ are isomorphic to their preimages in $\CC$.  
%	Moreover, the preimages of the components of the special fiber of $\hat \CX_1$ above $\bar p$ in  the normalizations in $K(C)$, respectively $K(\tilde C)$ are isomorphic. 
\end{Proposition}
%
%
%\begin{Prop}
%	Further assume that $v(F(2^\lambda x)-\tilde F(2^\lambda x)) > 2$ for all $\lambda \in [0,\alpha]\cap \BQ$. The stable reduction of $C$ over $\bar p$ is isomorphic to the stable reduction of $\tilde C$ over $\bar p$. 
%\end{Prop}
\begin{Proof}
	By Proposition \ref{OnlyComponentsOfTypeBApproxProp} and \cite[Prop. 4.5.1]{GP24}, 	the components of type (b) in the preimage of $\bar \CU$ in $\hat \CC$ and $\tilde \CU$ are the same.  
	Let $\lambda \in (0,\alpha)$ be a break point of $\wbar_{\bar p, X}$. Then $y=x/2^\lambda$ is a coordinate for the component of type (b) corresponding to $\lambda$. Moreover, we have $K(C)=K(y,z')$ with $z'=z y^{ \left \lceil{\op{ldeg}(F)/2}\right \rceil }$. Then 
	$$z'^2= y^{ 2\left \lceil{\op{ldeg}(F)/2}\right \rceil }F(2^\lambda y)$$
	 and $y^{ 2\left \lceil{\op{ldeg}(F)/2}\right \rceil }F(2^\lambda y)\in R[y]$. Moreover, the Laurent polynomials  $y^{ 2\left \lceil{\op{ldeg}(F)/2}\right \rceil }F(2^\lambda y)$ and $y^{ 2\left \lceil{\op{ldeg}(F)/2}\right \rceil }\tilde F(2^\lambda y)$ satisfy the conditions of Proposition  \ref{stableReductionOverSmoothUnmarkedPoint}, hence the stable reduction of $C$ and of $\tilde C$ over all the closed points of the component of type (b) are isomorphic.
	Moreover, by Proposition \ref{BareNormalizationSpecialFiberIsIsomorphic}, the preimage of the component of type (b) in the normalization in $K(\tilde C)$ is isomorphic to its preimage in $\CC$,
	 which  concludes the proof. 
\end{Proof}

\begin{Ex}
	\begin{enumerate}[(a)]\label{ExampleApproxDoublePoint} {\rm
		\item 		Fix an element $a\in R$ with $v(a)>2$ and let $L\in R^\times +xR[x]$. Define $$\tilde F:=F(x)\cdot L(ax).$$ If $\bar p$ is an even double point, the Laurent polynomial $\tilde F$ satisfies the prerequisites of Proposition \ref{ComponensofTypeCDApproxProp}. Likewise, if $\bar p$ is an odd double point,  the Laurent polynomial $\tilde F$ satisfies the prerequisites of Proposition \ref{Odddoublepointapprox}.
		\item Assume that $\bar p$ is an even double point of thickness $\alpha=2$ and let $$F=\left(2 x^{5}+2 x^{3}+x^{2}+x+1\right) \left(\frac{4}{x}+\frac{512}{x^{3}}+1\right).$$ Define $\tilde F:=x^2+x+1+\tfrac{4}{x}$. Then 
		$$F-\tilde F=2x^5+8 x^{4}+2 x^{3}+1032 x^{2}+4 x+1028+\frac{512}{x}+\frac{512}{x^{2}}+\frac{512}{x^{3}} $$
		and hence $v(F(2^\lambda u)-\tilde F(2^\lambda u))\geq \min\{1, 9-3\lambda\}\geq 1$. Moreover, as $$\tilde F=(x+1)^2-x+\frac{4}{x} $$ is an odd decomposition, we have 
		$\wbar(\tilde F(2^\lambda u))=\min\{\lambda, 2-\lambda\}.$ By Proposition \ref{OnlyComponentsOfTypeBApproxProp} there is exactly one component of type (b) over $\bar p$. }
	\end{enumerate} 
	\end{Ex}

%Todo: Make this into a nice proposition or get rid of it. 
%\begin{Prop}
%	Let $X$ be an irreducible component for $X$. 
%Let $x$ be a coordinate for $X$ and consider a hyperelliptic equation $z^2=F(x)$ of $C$ with $F\in R[x]$ and constant coefficient $1$. Let $S$ be the set of all roots of  $F$ that correspond to Weierstrass points of $C$ specializing to components which have distance less than $2$ to $X$ in the dual graph of $\bar C_0$. Here the distance of two components is the sum of the thicknesses of the double points along a shortest path in the dual graph of $C_0$.  Then we may replace  $F$ by
%$$\tilde F:=\prod_{\xi\in S} \left(1-\frac{x}{\xi}\right) $$  for the computation of $C_0$ above $X$. 
%\end{Prop}

%\begin{Rem}\label{DistanceRemark} { \rm
%Taking Examples \ref{ExampleApproxSmoothPoint} (a) and  \ref{ExampleApproxDoublePoint} (a) together, we see that for the local  calculation of $C_0$ over a component of $\bar C_0$ with component $X$, the following statements holds true:
%
%}
%
%
%Weierstrass points of $C$ that specialize to components with  distance larger than $2$ in the dual graph of $\bar C_0$ do not have to be taken into account. More precisely, let $x$ be a coordinate of $X$ and let $z^2=F(x)$ be a Weierstrass equation with $F\in R[x]$. Assume that $F=F_1 F_2$ with $F_1$   }
%\end{Rem}

\section{The local genus over a point}\label{SectionLocalGenus}
Recall that $\bar \CC$ is a semistable marked model of $(\bar C, \bar P_1, \dots, \bar P_{2g+2} )$ and that $\hat \CC$ is the minimal model of $\bar C$ whose normalization in $K(C)$ is semistable and which dominates $\bar \CC$. We denote the induced morphism on the special fibers by   $\Psi: \hat C_0\to \bar C_0$.

%be a semistable model of $\bar C$ that dominates the stable marked model of $\BP^1_K$. Let  $\Psi: \hat \CC\to \bar \CC$ be the model of $\bar C$ minimal with respect to dominating $\bar \CC$ and having its normalization $\CC$ in $K(C)$ being semistable. 
\subsection{Local genus}\label{SubSectionLocalGenus}

%ToDo: Change to local genus over point everywhere. Perhaps split this Definition and first define it for points of $\hat C_0$?
\begin{Def}\label{relativeGenusDefSemistable}
	Let $\hat p\in \hat C_0$ be a not necessarily closed point. Then the local genus $g_{\hat p}$ of $C_0$ over $\hat p$ is defined as follows.
	\begin{enumerate}[(a)]
		\item If $\hat p$ is the generic point of an irreducible component whose preimage in $C_0$ is reducible, set $g_{\hat p}:=-1$. 
		%normalization in $K(C)$ is split, set $g_{\hat p}:=-1$. 
		\item If $\hat p$ is the generic point of an irreducible component whose preimage in $C_0$ is  irreducible of genus $g'$, set $g_{\hat p}:=g'$.
		\item If $\hat p$ is a smooth closed point, set $g_{\hat p}:=0$. 
		%whose preimage in the normalization in $K(C)$ is smooth, set $g_{\hat p}:=0$. 
		\item If $\hat p$ is an ordinary double point whose preimage in the normalization in $K(C)$ are $n$ ordinary double points, set $g_{\hat p}:=n-1$. %Note that here $n\in \{1,2\}$. 
	\end{enumerate}
\end{Def}
	The curve $\hat C_0$ has only finitely many generic points and ordinary double points; hence the relative genus is $0$ for all but finitely many points. 	
%The local genus over a point of $\hat C_0$ is $0$ for all but finitely many points. 
\begin{Def}\label{relativeGenusDef}
	Let $\bar p\in \bar C_0$ be a not necessarily closed point. Then the local genus $g_{\bar p}$ of $C_0$ over $\bar p$ is 
	$$g_{\bar p}\ :=\ \sum_{\hat p\in \Psi^{-1}(\bar p)}g_{\hat p}.$$ 
\end{Def}

\begin{Prop}\label{LocalGenusSumIsGFormula}
	We have 
	$$g \ = \ \sum_{\bar p\in \bar C_0}g_{\bar p} .$$
\end{Prop}
\begin{Proof} 
%	 If the formula is true for points of $\hat C_0$, \dots  
	First, assume that $\bar \CC=\hat \CC$. 
%
%	Moreover, there are only finitely many points whose preimage in the normalization in $K(x,z)$ is singular, hence the relative genus is $0$ for all but finitely many points. 
%	Next, for a closed point $\bar p$ which falls under case (e), all the points in $\Psi^{-1}(\bar p)$ fall under the cases (a)-(d) of \ref{relativeGenusDef} as the normalization of $\hat \CC$ is semistable. 
%	Together with the previous argument applied to $\hat C_0$, we get that the relative genus of $C_0$ over $\bar p$ is a finite sum of integers and as such an integer itself.
%	By definition, the sum of all 
%	As $\Psi$ is a sequence of blow ups centered in points whose preimage in the normalization in $K(C)$ is singular, the sums of all relative genera of $\bar \CC$ and $\hat \CC$ are the same. Thus it is sufficient to prove the formula for $\hat \CC$. 
	By assumption,  the normalization of $\bar \CC$ in $K(C)$  is a semi-stable model $\CC$ of $C$. Hence the arithmetic genus of its special fiber $C_0$ is $g$. Let $\beta(C_0)$ be the first Betti number of the dual graph of $C_0$ and let $g_1$, \dots,  $g_n$ be the geometric genera of the irreducible components of $C_0$. Then 
	$$g(C_0) \ = \ g \ = \ \beta(C_0)+\sum_{i=1}^n g_i. $$
			Furthermore, let $\eta_1, \dots, \eta_m$ be the generic points of the irreducible components of $\bar C_0$ whose normalization are split and let $\bar p_1, \dots, \bar p_L$ be the double points of $\bar C_0$.  By  \cite[Prop. 3.48]{LiuAlgGeo2002}, the preimages of each of the points $\bar p_1, \dots, \bar p_L$ are one or two double points.
			
%			. Assume that the preimage in $C_0$ of a smooth unmarked point $p\in \hat C_0$ contains a singularity. Then by \cite[Prop. 3.3.3]{GP24}, there is at least one component of type (d) above $\bar p$, contradicting $\hat \CC= \bar \CC$. 
%			Together with \cite[Prop. 4.1.1]{GP24}, we get that the preimage of the smooth locus of $\bar C_0$ does not contain double points. 
%			    Moreover, by   \cite[Prop. 3.3.3]{GP24}, every double point of $C_0$ maps to ${\bar p_1, \dots, \bar p_L}$.   
%	Note that $\CC \to \bar \CC$ is a double covering and the ordinary double points of $C_0$ are mapped to ordinary double points of $\tilde C_0$.

	 By \cite[Prop. 4.6.6]{GP24},  the first Betti number of $C_0$ is equal to the number of double points of $\bar C_0$ whose preimage in $C_0$ are two double points minus the number of irreducible components whose normalization in $K(C)$ is split. In other words
	$$\beta(C_0) \ = \  \sum_{i=1}^m g_{\eta_i}+\sum_{i=1}^L g_{\bar p_i}.$$
	Finally, let $\eta_{m+1}, \dots, \eta_{n}$ be the generic points of the remaining irreducible components of $\tilde C_0$. As the normalization of an irreducible component of $\tilde C_0$ can only have positive genus if it is not split, we get 
	$$\sum_{i=1}^n g_i\ = \ \sum_{i=m+1}^n g_{\eta_i}.$$
	 Together
	 $$\sum_{\bar p\in \bar C_0}g_{\bar p} \ = \ \sum_{i=1}^m g_{\eta_i}+\sum_{i=1}^L g_{\bar p_i} + \sum_{i=m+1}^n g_{\eta_i} \ = \ \beta(C_0)+ \sum_{i=1}^n g_i = g,$$
	 	as desired. 
	 	
	 Second, assume that $\bar \CC \neq \hat \CC$. 
%	 As the local genus over a point of $\hat C_0$ is $0$ for all but finitely many points, there can only be finitely $\bar p\in \bar C_0$ for which the sum of the local genera of their preimages is non-trivial, proving the first assertion for $\hat C_0$. Moreover, w
We have
	 $$\sum_{\bar p\in \bar C_0}g_{\bar p}\ =\ \sum_{\bar p\in \bar C_0}\sum_{\hat p\in \Psi^{-1}(\hat p)}g_{\hat p}\ =\ \sum_{\hat p\in \hat C_0} g_{\hat p}=g.$$
	\end{Proof}

%	\begin{Ex}
%	For any odd double point $\bar p\in \bar C$, we have $g_{\bar p}=0$ by Prop. 4.4.1 in \cite{GP24}.
%	\end{Ex}
\begin{Prop}\label{LocalGenusTrivialForMarkedAndOddDP}
	Let $\bar p \in \bar C_0$ be a marked point or an odd double point. Then $g_{\bar p}=0$. 
\end{Prop}
\begin{Proof}
	If $\bar p$ is an odd double point, we know that the morphism $\Psi: \hat C_0\to \bar C_0$ restricts to an isomorphism in a neighborhood of $\bar p$ and that $C_0$ has a unique double point above $\bar p$ by \cite[Prop. 4.4.1]{GP24}. 
	If $\bar p$ is a marked point, we know that $\bar p$ is smooth and that $\Psi$ restricts to an isomorphism in a neighborhood of $\bar p$ by \cite[Prop. 4.1.1]{GP24}.
	In both cases, this implies $g_{\bar{p}}=0$. 
\end{Proof}

In the following subsections, we develop methods for computing the local genus of $C_0$ that do not require the computation of a model of $C$.  
%	For a given curve $C$ and a point $\bar p\in \bar C$, it is a priori not clear whether 
%To motivate the 
%	In the following sections, we will develop algorithms for the computation of the local genus which avoid computing the models $\hat \CC$ and $\CC$. 
	
\subsection{Smooth unmarked points}\label{SubSectionsLocalGenusUnmarkedPoints}
Let $\bar p\in \bar C_0$ be a closed smooth unmarked point contained in the irreducible component $X$. We identify a neighborhood $\bar\CU \subset \bar \CC$ with an open subscheme of $\Spec R[x]$. %such that $\bar p$ is given by $x=0$. 
Moreover, suppose that the function field of $C$ is $K(x,z)$ with $z^2=F (x)$ for a separable polynomial $F  \in R[x]$ with $v(F )=0$. 
%
%Then by assumption $\bar p$ does not lie in the zero locus of~ $F$.  By Proposition 3.3.2 of GP24, the normalization of $R[x]$ is smooth if $w(F)\geq 2$.  Hence assume that $\gamma:=\wbar(F)<2$. 
Let  $F=G+H^2$ be an odd decomposition.

\begin{proposition}
	Assume that $\wbar(F)=2$. Then $g_{\bar p}=0$. 
\end{proposition}
\begin{Proof}
	 By \cite[Prop. 3.3.2]{GP24}, the normalization of $R[x]$ in $K(C)$ is smooth. Hence the morphism $\Psi$ is an isomorphism in a neighborhood of $\bar p$ and so the local genus is $0$.   
\end{Proof}

%ToDo: Make this more like this: In the following, let $ \CX$ be a semistable model of $\bar C$ and consider a closed smooth point $\bar p$ of the special fiber $X_0$ of $\CX$.
%%Choose a coordinate $x$ of $\bar \CC$ such that $\bar p$ is defined by $x=0$ and let $ \CU$ be an open neighborhood of $\bar p$ isomorphic to an open subscheme of $\Spec R[x^{}]$. 
%We identify a neighborhood $\bar\CU \subset \bar\CC$ with an open subscheme of $\Spec R[x^{}]$.
%%Let $U_0$ be the special fiber of $\bar \CU$.
%Furthermore, let $F\in R[x]$ with $v(F)=0$ be such that an equation for $C$ is given by $z^2=F$. Now fix $\tilde F\in R[x^{}]$ %with $v(F-\tilde F)>2$ 
%%the field $K(x,\tilde z)$ with $\tilde z^2=\tilde F$. % 
%and let $\tilde C$ be the curve  defined by $z^2=\tilde F(x)$. 
%Assume that $v(F-\tilde F)>\wbar (f):=\gamma$ and let $\bar U_0$ be the reduction of $\bar \CU$ in $X_0$. 

%ToDo: reminder of where this equation comes from, explain that it is independent of choice. 
\begin{Thm}\label{LocalGenusUnmarkedPointProp}
	 Assume that $\gamma:=\wbar(F)<2$. 
	The local genus $g_{\bar p}$ of $C_0$ over $\bar p$ is equal to half the multiplicity of the root of $[\frac{dG}{dx}/2^\gamma]$ at $\bar p$. 
	%	$g_{\bar{p}}=\frac{1}{2}({\op{ldeg}([G/2^\gamma])-1})$.
\end{Thm}
\begin{Proof} %ToDo: reordering
	%%%New:
	As $v(G/2^\gamma)=0$, all odd terms of $\frac{dG}{dx}/2^\gamma$ have valuation at least $1$ and so $[\frac{dG}{dx}/2^\gamma]$ is a polynomial in $x^2$. 
	 Let $n\in \BZ$ be half the multiplicity of the root of $[\frac{dG}{dx}/2^\gamma]$ at $\bar p$. %Hence $n\in \BZ_{\geq 0}$. 
 We perform induction on $n$. 
 
  For $n=0$, the normalization of a neighborhood of $\bar p$ in $K(C)$ is smooth by \cite[Prop. 3.3.2]{GP24}, hence $g_{\bar p}=0$. 
  
 For $n>0$, the same proposition implies that there exists an irreducible component  of type (d) over $\bar p$. We fix one such component $\bbar T$ of type (d). It is given by a coordinate $y$ with $x=2^\epsilon y+ \xi_0$ for a rational $\epsilon>0$ and $\xi_0\in R$. After replacing $x$ by $x-\xi_0$, the coordinate has the form $x= 2^\epsilon y$. Moreover, as now $\bar p$ is given by $x=0$ and $G$ is odd, the multiplicity of the root of   $[\frac{dG}{dx}/2^\gamma]$ is equal to $\frac{1}{2}({\op{ldeg}([G/2^\gamma])-1})$. In other words, we need to prove $g_{\bar{p}}=\frac{1}{2}({\op{ldeg}([G/2^\gamma])-1})$.
% Moreover by Lemma \ref{ExchangeOptimalForOddLem}, we can choose the decomposition to be odd. 

 	Let $ \CY$ be the semistable model of $\bar C$ obtained from $\bar \CC$ by adjoining an irreducible component of coordinate $u=x/2^{\epsilon}$. Denote the special fiber of $\CY$ by $Y_0$ and let $\bar q\in Y_0$ be the unique ordinary double point over $\bar p$. Let $\hat X$ be the proper transform of $X$ in $\CY$.  	 As $\hat T$ is a component of type (d), its preimage  in $\CC$ is an irreducible component $T$ of positive genus and thus we have $\wbar_{\bar q, \hat X}(\epsilon)=2$. Moreover, by Proposition \ref{RaynaudProp}, there is exactly one double point in $\CC$ over the double point of $\hat T$, hence the slope of $\wbar_{\bar q, \hat X}$ to the left of $\epsilon$ is non-zero by Corollary \ref{UniqueDoublePointOverDoublePoint}. Together with $\wbar_{\bar q, \hat X}$ being concave with maximal possible value $2$, we get that it is a strictly increasing function.

% 	Consider the function $\wbar: [0,\epsilon]\cap \BQ \to \BR$ from \eqref{WBarDef}.

 	First, assume that the function $\wbar_{\bar q, \hat X}$ has no break point. By \cite[Prop. 4.5.1]{GP24}, this implies that there are no components of type (b) over $\bar q$, which is equivalent to $\bar T$ being the only component over $\bar p$. It remains to prove that $\frac{1}{2}({\op{ldeg}([G/2^\gamma])-1})=g(T)$.  	
 	As not having a break point means that $\wbar_{\bar q, \hat X}$ is linear it has a single slope $s_1>0$.  Now Lemma \ref{slopeLemLeftSide} yields that $\op{ldeg}([G/2^\gamma]) =   s_1 $ and Lemma 
 \ref{slopeLemRightSide} implies that $s_1  =  \op{deg}([G(2^\epsilon y)/4]).$ Together, we get
 	\UseTheoremCounterForNextEquation
\begin{equation}\op{ldeg}([G/2^\gamma])\ = \   \op{deg}([G(2^\epsilon y)/4]).\label{LeftSideRightSideSlopeEquation}\end{equation}
	By \cite[Prop. 3.3.2]{GP24}, an equation for $T$ in $C_0$ is given by 
	$$[H(2^\epsilon y)]t+t^2=[G(2^\epsilon y)/4].$$ 
	As $\epsilon>0$, we have $[H(2^\epsilon y)]=1$ and hence the genus of $T$ is
$g(T)=\frac{1}{2}(\op{deg}([G(2^\epsilon y)/4])-1)$ by the genus formula for Artin-Schreier curves. Inserting equation \eqref{LeftSideRightSideSlopeEquation}, we get $g(T)=\frac{1}{2}({\op{ldeg}([G/2^\gamma])-1})$, as desired. 

Second, assume that $\wbar_{\bar q, \hat X}$ has a break point.  We call the smallest such point $\lambda_1$. Let $s_1$ be the slope of $\wbar_{\bar q, \hat X}$ to the right of $0$. As $\lambda_1$ is the smallest break point, the slope to the left of $\lambda_1$ is $s_1$. 
Now Lemma \ref{slopeLemLeftSide} yields that $\deg([G(2^{\lambda_1})/2^{\wbar_{\bar q, \hat X}(\lambda_1)}])=s_1 $ and Lemma 
\ref{slopeLemRightSide} yield that $s_1=\op{ldeg}([G/2^\gamma]).$  Together, we get 
	\UseTheoremCounterForNextEquation
\begin{equation}\op{ldeg}([G/2^\gamma])\ = \ s_1 \ = \   \deg([G(2^{\lambda_1})/2^{\wbar_{\bar q, \hat X}(\lambda_1)}]).\label{LeftSideRightSideFirstSlopeEquation}\end{equation}
Moreover, as $\wbar_{\bar q, \hat X}$ is strictly increasing  and $\wbar_{\bar q, \hat X}(\epsilon)=2$, we have $\wbar_{\bar q, \hat X}(\lambda_1)<2$.
Hence by Lemma 
\ref{slopeLemRightSide}  the slope $s_2$ to the right side of the break point is 
	\UseTheoremCounterForNextEquation
\begin{equation}s_2=\op{ldeg}([G(2^{\lambda_1})/2^{\wbar_{\bar q, \hat X}(\lambda_1)}]).\label{RightSideRightSideFirstBreakPointEquation}\end{equation}

% By Lemmata \ref{slopeLemLeftSide} and \ref{slopeLemRightSide} yield that the slope $s_1$ to the left side of the break point is  $$s_1=\deg([G(2^{\lambda_1})/2^{\wbar(\lambda_1)}])=\op{ldeg}([G/2^\gamma]).$$ 

Let $\bar \CC '$ be the semistable model of $\bar C$ obtained from $\bar \CC$ by adjoining an irreducible component $X$ of coordinate $u=x/2^{\lambda_1}$. We have 
$$\wbar( X)=\wbar_{\bar q, \hat X}(\lambda_1)<2.$$
%As $\wbar$ is strictly monotone  and $\wbar(\epsilon)=2$, we have $\wbar(\lambda_1)<2$.  
Moreover,  Equations \ref{LeftSideRightSideFirstSlopeEquation} and \ref{LeftSideRightSideSlopeEquation} yield
 $$[G(2^{\lambda_1}u)/2^{\wbar_{\bar q, \hat X}(\lambda_1)}]=P u^{s_2},$$ 
 where $P\in k[u]$ is a polynomial of degree $s_1-s_2$ with non-vanishing constant coefficient. Now by \cite[Prop. 3.3.2]{GP24}, the preimage of $X$ in the normalization of $\bar \CC '$ in $K(C)$ is smooth outside finitely many points where $d(P u^{s_2})/du=0$. Counted with multiplicities, the polynomial $\frac{d}{du} P u^{s_2}=u^{s_2-1}(P+uP')$ has $s_1-s_2$ non-zero roots. Each of these roots has multiplicity at most $s_1-s_2<n$ and we apply the induction hypothesis to each of them. This yields that the sum of the relative genuses of all closed points of this component is $\frac{s_1-s_2}{2}$ without the point $u=0$. Over the point $u=0$, we again apply the induction hypothesis to find that the relative genus over it is $\frac{s_2-1}{2}$. Together, we get that the relative genus over $\bar p$ is $\frac{s_1-s_2}{2}+\frac{s_2-1}{2}=\frac{s_1-1}{2}$. Inserting Equation \ref{LeftSideRightSideFirstSlopeEquation} yields the desired equality.

\end{Proof}

%\subsection{Marked smooth points} 
%\subsection {Odd Double point}

\subsection{Even double points}\label{SubSectionEvenDoublePointsLocalGenus}
Now we assume that $\bar p$ is an even double point of~$\bar C_0$. Let $X$ be an irreducible component containing $\bar p$. For each pair $m,n$ of integers write $\delta_{mn}$ for the Kronecker delta 
$$\delta_{m n}:=\begin{cases}
	1 & \text{ for } m=n \\
	0 & \text{ otherwise. }
\end{cases}$$
%
%%old:
%Now we assume that $\bar p$ is an even double point of~$\bar C_0$ and identify a neighborhood $\bar\CU \subset \bar \CC$ with an open subscheme of $\Spec R[x,y]/(xy-a)$ for some nonzero $a\in\Fm$, such that $\bar p$ is given by $x=y=0$. After multiplying $a$ by a unit we may assume that $a=2^\alpha$ for some $\alpha>0$, where $\alpha$ is the thickness of~$\bar p$.
%Identifying $y$ with $\frac{2^\alpha}{x}$, the ring in question is isomorphic to the subring $R[x,\tfrac{2^\alpha}{x}]$ of the ring of Laurent polynomials $K[x^{\pm1}]$. 
%Moreover, there exists $F  \in R[x,\frac{2^\alpha}{x}]$ with constant coefficient $1$ such that function field of $C$ is $K(x,z)$ with $z^2=F (x)$ . %Let $F=H^2+g$  be an optimal decomposition. 

\begin{Thm}\label{RelativeGenusOverDoublePoint}
Let $s_1$  and $s_2$ be the largest, respectively the smallest slope of $\wbar_{\bar p, X}$. Then the local genus of $C$ over $\bar p$ is
$$g_{\bar p}=\frac{s_1-s_2+\delta_{0s_1}+\delta_{0s_2}}{2} .$$ 	
\end{Thm}
\begin{Proof}
%	Identify a neighborhood $\bar\CU \subset \bar \CC$ with an open subscheme of $\Spec R[x,y]/(xy-a)$ for some nonzero $a\in\Fm$, such that $\bar p$ is given by $x=y=0$. After multiplying $a$ by a unit we may assume that $a=2^\alpha$ for some $\alpha>0$, where $\alpha$ is the thickness of~$\bar p$.
%	Identifying $y$ with $\frac{2^\alpha}{x}$, the ring in question is isomorphic to the subring $R[x,\tfrac{2^\alpha}{x}]$ of the ring of Laurent polynomials $K[x^{\pm1}]$. 
%	Moreover, there exists $F  \in R[x,\frac{2^\alpha}{x}]$ with constant coefficient $1$ such that function field of $C$ is $K(x,z)$ with $z^2=F (x)$ . Let $F=H^2+G$  be an odd decomposition. 
%	
	As in Section \ref{SubSectionDoublePointApprox}, identify a neighborhood $\bar\CU \subset \bar \CC$ with an open subscheme of $\Spec R[x,y]/(xy-2^{\alpha})$, where $\alpha$ is the thickness of $\bar p$ and  $\bar p$ is given by $x=y=0$. Moreover, fix $F  \in R[x,\frac{2^\alpha}{x}]$ with constant coefficient $1$ such that function field of $C$ is $K(x,z)$ with $z^2=F (x)$ . Let $F=H^2+G$  be an odd decomposition. 	
	Then $\wbar_{\bar p, X}(\lambda)=\wbar(F(2^\lambda u))$ for all $\lambda \in \BQ\cap [0,\alpha]$ and a new variable $u$.
	
		We proceed by induction on the number of slopes of $\wbar_{\bar p, X}$, denoted by $n$. 
		\begin{Lem}
			Theorem \ref{RelativeGenusOverDoublePoint} holds true when $\wbar_{\bar p,X}$ has  $n=1$ slopes.
		\end{Lem}
		\begin{Proof}
Assume that $\wbar_{\bar p, X}$ has only one slope. By \cite[Prop. 4.5.1]{GP24} the morphism $\Psi$ restricts to an isomorphism in a neighborhood of $\bar p$. If $s_1=s_2=0$, the only slope of $\wbar_{\bar p, X}$ is $0$. Then there are two double points in $\hat \CC$ over $\bar p$ by \cite[Prop. 4.5.10]{GP24}. Hence by definition $g_{\bar p}=1$, as desired. If $\wbar_{\bar p, X}$ has non-zero slope,  by \cite[Prop. 4.5.12]{GP24} there is exactly one double point in $C_0$ over $\bar p$. This yields $g_{\bar p}=0$, as desired. 
	\end{Proof}
 Now assume $n>1$  and let $\lambda_1$ be the smallest break point of $\wbar_{\bar p, X}$. Let $s_3$ be the second largest slope of $\wbar_{\bar p, X}$. 
 As $\wbar_{\bar p, X}$ is a concave function, its slope to the left of $\lambda_1$ is $s_1$ and its slope to the right of $\lambda_1$ is $s_3$. 
 
	Let $\bar \CC '$ be the semistable model of $\bar C$ obtained from $\bar \CC$ by adjoining an irreducible component $\bar T$ of coordinate $u=x/2^{\lambda_1}$. Then the model $\bar \CC '$ dominates $\bar \CC$ and is dominated by $\hat \CC$.
	Furthermore, let $T$ be the inverse image of $\bar T$ in $C_0$. 
	
	\begin{Lem}\label{FirstLemmaLongProof}
		The induction step is valid if $s_1=0$. 
	\end{Lem}
	\begin{Proof}
If $s_1=0$, we have $\wbar_{\bar p, X}(\lambda_1)=2$ and  $s_3<0$, hence $\wbar_{\bar p,X}(\lambda)<2$ for $\lambda>\lambda_1$. This implies  $w(F(2^{\lambda_1} u))=2$. 
Moreover, as $F\in R[x,\frac{2^\alpha}{x}]$ has constant coefficient $1$ and $\lambda_1 \in(0,\alpha)$, we have  $[H(2^{\lambda_1}u)]=1$.
	 By \cite[Prop. 3.3.2]{GP24},  the equation for $T$ modulo $\Fm$ has the form 
%		$$[H(2^{\lambda_1}u)]t+t^2=[G(2^{\lambda_1}u)/4]$$
$$t+t^2=[G(2^{\lambda_1}u)/4]$$
		and is smooth outside two double points at $u=0,\infty$ by \cite[Lemma 4.5.2 (b)]{GP24}.
	As ${F\in R[x,\frac{2^\alpha}{x}]}$ with constant coefficient $1$ and $0<\lambda_1<\alpha$, we have 
	
	By Lemma \ref{slopeLemLeftSide} the equality $s_1=0$  implies that $\op{deg}([G(2^{\lambda_1}u)/4])\leq 0$. 
	 Similarly, by Lemma  \ref{slopeLemRightSide} $s_3$ being strictly negative implies $s_3=\op{ldeg}([G(2^{\lambda_1}u)/4])$. Together with the genus formula for Artin-Schreier curves, we get
	 $$g(T)\ =\ \frac{-\op{ldeg}([G(2^{\lambda_1}u)/4])-1}{2} \ =\ \frac{-s_3-1}{2}.$$ 
	 At the double point $\bar q$ defined by $u=\infty$, the function $\wbar_{\bar q, \bar T}$ has $n-1$ slopes, where $s_3$ is its largest and $s_2$ is its smallest slope. By induction hypothesis, the relative genus of $C_0$ over the point defined by $u=\infty$ is $\frac{s_3-s_2+\delta_{s_2=0}}{2}$. Moreover, the relative genus at the double point at $u=0$ is $1$, as the square defect function of that point is linear of slope $0$. Together, we get 
	$$g_{\bar p}= 1+ \frac{-s_3-1}{2}+\frac{s_3-s_2+\delta_{0s_2}}{2}=\frac{s_1-s_2+\delta_{0s_1}+\delta_{0s_2}}{2},$$
	as desired.
	\end{Proof}
	
		\begin{Lem}\label{SecondLemmaLongProof}
		The induction step is valid if $s_1 \neq 0$ and $s_3\neq 0$ and $\wbar_{\bar p, X}(\lambda_1)<2$.
	\end{Lem}
	\begin{Proof}
 By \cite[Prop. 3.3.2]{GP24}, the equation for $T$ modulo $\Fm$ has the form 
	$$t^2=[G(2^{\lambda_1}u)/2^{\wbar_{\bar p, X}(\lambda_1)}].$$
	Moreover by Lemmata \ref{slopeLemLeftSide} and \ref{slopeLemRightSide}, the slopes are given by $s_1=\deg([G(2^{\lambda_1}u)/2^{\wbar_{\bar p, X}(\lambda_1)}])$ and $s_3=\op{ldeg}([G(2^{\lambda_1}u)/2^{\wbar(\lambda_1)}])$, thus we can write $[G(2^{\lambda_1}u)/2^{\wbar_{\bar p, X}(\lambda_1)}]=P u^{s_2}$, where $P\in k[u]$ is a polynomial of degree $s_1-s_3$ with non-zero constant coefficient.  Set ${d_1:= \left \lceil{ \op{ldeg}(F)/2}\right \rceil}$ and substitute $z$ by $\frac{w} { u^{-d_1 }}$ . Then 
	$w^2=F(2^{\lambda_1}u)u^{2d_1}$ is an equation defining $C$ and we have $\tilde F:= F(2^{\lambda_1}u)u^{2d_1}\in R[x]$ with $v(\tilde F)=0$. 
	Moreover $$\tilde F(2^{\lambda_1}u)=u^{2 d_1}H^2(2^{\lambda_1}u)+ u^{2d_1}G((2^{\lambda_1}u))$$ is an odd decomposition and hence also optimal. Again by \cite[Prop. 3.3.2]{GP24}, this yields an equation modulo $\Fm$ for $T$ of the form 
	$$r^2=Pu^{s_2+2d_2}. $$
	Now by \cite[Prop. 3.3.2]{GP24}, the normalization of this component in $K(C)$ is smooth outside finitely many points where $d(P u^{s_2+2d_1})/du=0$. Counted with multiplicities, the polynomial $\frac{d}{du} P u^{s_2+2d_1}=u^{s_2+2d_1-1}(P+uP')$ has $s_1-s_3$ non-zero roots. By Theorem \ref{LocalGenusUnmarkedPointProp} this yields that the sum of the relative genera at all closed points of $\bar T$ without the points defined by $u=0,\infty$ is $\frac{s_1-s_3}{2}$. At the double point $\bar q$ defined by $u=\infty$, the function $\wbar_{\bar q, \bar T}$ has $n-1$ slopes, where $s_3$ is its first and $s_2$ its last slope. By induction hypothesis, the relative genus of $C_0$ over this  point is $\frac{s_3-s_2+\delta_{s_2=0}}{2}$. Moreover $\wbar_{\bar r,T}$ is linear with  non-zero slope $-s_1$. By \cite[Prop. 4.5.12]{GP24}, this implies that there is a unique double point of $C_0$ over $\bar r$. Hence by definition $g_{\bar r}=0$.  
Together, we get 
	$$g_{\bar p}= \frac{s_1-s_3}{2}+\frac{s_3-s_2+\delta_{s_2=0}}{2}=\frac{s_1-s_2+\delta_{0s_1}+\delta_{0s_2}}{2},$$
	as desired.
	\end{Proof}
	\begin{Lem}\label{ThirdLemmaLongProof}
	The induction step is valid if $s_1\neq 0$ and $s_3 \neq 0$ and $\wbar_{\bar p, X}(\lambda_1)=2$
	\end{Lem}
	\begin{Proof}
	As in the proof of Lemma \ref{FirstLemmaLongProof}, we get $w(F(2^{\lambda_1} u))=2$ and $[H(2^{\lambda_1}u)]=1$.   By \cite[Prop. 3.3.2]{GP24}, the equation for $T$ modulo $\Fm$ has the form 
		$$t+t^2=[G(2^{\lambda_1}u)/4]$$ 
		and is smooth outside two double points at $u=0,\infty$  by by \cite[Lemma 4.5.2 (b)]{GP24}.
		
		By concavity of $\wbar_{\bar p, X}$, we have $s_1>0$ and $s_3<0$.  
		Hence by Lemmata \ref{slopeLemLeftSide} and \ref{slopeLemRightSide}, the slopes are given by $s_1=\deg([G(2^{\lambda_1}u)/2^{\wbar(\lambda)}])$ and $s_3=\op{ldeg}([G(2^{\lambda_1}u)/2^{\wbar(\lambda_1)}])$. By the genus formula for Artin-Schreier curves, the genus of $T$ is
		$$g(T)\ = \ \frac{1}{2}(\deg([G(2^{\lambda_1} u)/4])-\op{ldeg}([G(2^{\lambda_1} u)/4]))\ =\ \frac{s_1-s_3}{2}.$$
For the double point $\bar q$ at $u=\infty$, the function $\wbar_{\bar q, \bar T}$ has $n-1$ slopes, where $s_3$ is its largest and $s_2$ its smallest slope. By induction hypothesis, the relative genus of $C_0$ over this  point is 
$$g_{\bar q}\ = \ \frac{s_3-s_2+\delta_{0s_2}}{2}.$$ 
For the double point $\bar a$ at $u=0$, the function $\wbar_{\bar a,T}$ is linear with  non-zero slope $-s_1$. So by \cite[Prop. 4.5.12]{GP24},  there is a unique double point of $C_0$ over $\bar a$. Hence by definition $g_{\bar a}=0$.  
Together, we get 
$$g_{\bar p}= \frac{s_1-s_3}{2}+\frac{s_3-s_2+\delta_{0s_2}}{2}=\frac{s_1-s_2+\delta_{0s_1}+\delta_{0s_2}}{2},$$
as desired.
\end{Proof}

\begin{Lem} \label{FourthLemmaLongProof}
	The induction step is valid if $s_1\neq 0$ and $s_3=0$.
\end{Lem}
\begin{Proof}
Let $Y$ be the irreducible component of $\bar C_0$ with coordinate $y$. The slopes of $\wbar_{\bar p, Y}$ are in reversed order to those of $\wbar_{\bar p, X}$ by Theorem-Definition \ref{SquareDefectFunctionTheoremDefinition}. As $\wbar_{\bar p, X}$ is concave, this reduces the situation to Lemma \ref{FirstLemmaLongProof} if $n=2$ or to Lemma \ref{SecondLemmaLongProof} or \ref{ThirdLemmaLongProof}  otherwise. 
\end{Proof}
Together, the last four Lemmata show that the induction step is valid in all cases, which concludes the proof.
		\end{Proof}

\section{Computational criteria}\label{SectionComputationalCriterions}
\subsection{Thickness bound over smooth unmarked points}\label{SubSectionThicknessBound}

Let $X$ be an irreducible component of $\bar C_0$ and consider a closed smooth unmarked point $\bar p\in X$.  Fix an irreducible component $\bbar T$ of $\hat C_0$ above~$\bar p$ of type (d). Let $T$ be its inverse image in $C_0$ and write $g(T)$ for the genus of $T$. Furthermore, let $\CY$ be the semistable model of $\bar C$ obtained from blowing down all components of type (c) and (d) above $\bar p$ other than $\bbar T$. Then $\CY$ has a unique double point $\bar q$ above $\bar p$.
%By Theorem 4.2.2 of GP24, there exits $\xi_0\in R$ with $S_F(\xi_0)=0$ and $\epsilon>0$ such that $\bbar T$ is given by the substitution $x=\xi_0+2^\epsilon y$. 
\begin{Prop}\label{ComputeThickness}
	Let $\epsilon$ be thickness of $\bar q$. Then $$\epsilon\leq \frac{2-\wbar(X)}{2g(T)+1}.$$ Moreover, equality holds if and only if $\bbar T$ is the only component of $\hat C_0$ above $\bar p$.  
\end{Prop}
\begin{Proof}
	Let $\hat X$ be the proper transform of $X$ in $\CY$. By definition, we have $$\wbar_{\bar q, \hat X}(0)=\wbar(\hat X)=\wbar(X).$$ As $T$ has positive genus, we have $\wbar_{\bar q, \hat X}(\epsilon)=2$. Write $s_1$ for the largest slope of $\wbar_{\bar q, \hat X}$. 
	As the function $\wbar_{\bar q, \hat X}$ is concave, we get 
	\UseTheoremCounterForNextEquation
	\begin{equation}\label{FirstEpsilonInequality}
\epsilon \leq \frac{2-\wbar(X)}{s_1}
	\end{equation}
	Moreover, in the proof of Theorem \ref{LocalGenusUnmarkedPointProp} we have showed that the slope of $\wbar_{\bar q, \hat X}$ to the left of $\epsilon$ is $2g(T)+1$, hence $s_1\geq 2g(T)+1$. Together with \eqref{FirstEpsilonInequality}, this yields the desired inequality.

Furthermore, if equality holds, we must have $s_1=2g(T)+1$. Hence the largest and the smallest slopes of $\wbar_{\bar q, \hat X}$ are equal. As the function is concave and piece wise linear, this already implies that it is linear. By \cite[Prop. 4.5.1]{GP24}, we get that $\hat T$ is the only component of type (d) above $\bar p$. 

Finally, if there is only one component of type (d) above $\bar p$, the function $\wbar_{\bar q, \hat X}$ must be linear,  from which equality follows. 
\end{Proof}

\begin{Ex}\label{ComputeThicknessExercise}{\rm
	\begin{enumerate}[(a)]
		\item Consider the curve defined by $z^2=x^4+3x^3+3x^2+4x+1+8x^{-1}$. Then $\bar C_0$ has two irreducible components meeting in a double point $\bar p$ of thickness $\alpha=3$ defined by $x=0$. By \cite[Ex. 3.4.7]{GP24}, there are two components $\hat T_{1}, \hat T_{2}$ of  type (b) above the double point with $\wbar(\hat T_1)=\frac{3}{2}$ and $\wbar(\hat T_2)=2$. 
		
		By \cite[Ex. 5.2.3 (c)]{GP24} the stable reduction of the curve is of type (B5) of the genus $2$ classification. Moreover, the preimage of $\hat T_2$ in $C_0$ is a component of geometric genus $1$ and there is a component of type (d) over a smooth point of $\hat T_1$. 
		By Proposition \ref{ComputeThickness}, the double point in which the proper transform of $\hat T_1$ intersects the component of type (d) has thickness 
		$$\frac{2-\frac{3}{2}}3 \ = \ \frac{1}{6}. $$		
		\item Assume that $C$ is any genus $2$ curve whose stable reduction is of type (B5) and let $\alpha$ and $\delta$ be as in \cite[Section 5.2]{GP24}. Then $\hat C_0$ has one component of type (d), which meets a component $\hat T_1$ of type (b). The computations in  \cite{WorksheetsGP24} show that $\wbar(T_1)=\frac{3\delta}{2}$, hence by Proposition \ref{ComputeThickness}, the double point in which $\hat T_1$ intersects the component of type (d) has thickness
		$$\frac{4-3\delta}{6} $$
	\end{enumerate}}
\end{Ex} 	

\begin{Rem} {\rm Using Proposition \ref{ComputeThickness} as in Example \ref{ComputeThicknessExercise} one can compute the thicknesses of all double points of $C_0$ purely in terms of $\alpha, \beta, \gamma, \epsilon, \delta$ for all curves of genus $2$. In an upcoming article \cite{G25}, we will use this technique to enrich the genus $2$ classification by the values of thicknesses.}
\end{Rem}

\subsection{Grounded double points}\label{SubSectionGroundedDoublePoints}
%In this section, let $\bar p$ be an even double point o
%Let $\bar p\in \bar C_0$ be an even double point. Moreover, let $X,Y$ be the irreducible components of $\bar C_0$ containing $\bar p$. 
\begin{Def}\label{GroundedDef}
	An even double point $\bar p\in \bar C_0$ is called \emph{grounded} if for the irreducible components $X,Y$ of $\bar C_0$ containing it, we have  $\wbar(X)=\wbar(Y)=0$. 
\end{Def}
In the following section, we keep the notation of Definition \ref{GroundedDef} and denote the thickness of $\bar p$ by $\alpha$. Furthermore, we assume that $\bar p$ is grounded. 

\begin{Prop} \label{GroundedDPAtlestgenus1}
 We have $g_{\bar p}>0$.
% and there exists at least one component of type (b) above $\bar p$.  
\end{Prop}
\begin{Proof}
%	Let $\alpha$ be the thickness of $\bar p$. 
	By assumption, we have $\wbar_{\bar p, X}(0)=\wbar(X)=0$ and $\wbar_{\bar p, X}(\alpha)=\wbar(Y)=0$. As $\wbar_{\bar p, X}$ is a piecewise linear function and only has horizontal segments of value $2$ by \cite[Prop. 3.4.4]{GP24}, its first slope is positive and is last slope is negative. 
%	Hence $\wbar_{\bar p, X}$ has at least one break point and thus there exists a component of type (b) by \cite[Prop 4.5.1]{GP24}. 
%	Furthermore, 
	By Theorem \ref{RelativeGenusOverDoublePoint}, the local genus over $\bar p$ is strictly positive. 
\end{Proof}

\begin{Cor}
	Let $\hat p\in \hat C_0$ be a double point (odd or even). Then $\hat p$ is odd if and only if for the irreducible components $\hat X, \hat Y$ of $\hat C_0$ containing it, we have $\wbar(\hat X)=\wbar(\hat Y)=0$. 
\end{Cor}
\begin{Proof}
	First, assume that $\hat p$ is odd. Then $\wbar(\hat X)=\wbar(\hat Y)=0$ by Proposition \ref{SquaredefectZeroProp}. 
	Second, assume that $\hat p$ is even and $\wbar(\hat X)=\wbar(\hat Y)=0$. Then $\hat p$ is a grounded double point, hence  by Proposition \ref{GroundedDPAtlestgenus1} we have $g_{\hat p}>0$. Using the definition of local genus for double points and Proposition \ref{OnlyOneDoublePointProp}, we get $g_{\hat p}=0$, which is a contradiction. 
\end{Proof}

\begin{Prop}\label{Genus1BehaviouroverGroundedDP}
	Assume  $g_{\bar p}=1$. Then the stable reduction above $\bar p$ has the same structure as in the $g=1$ case. More specifically: 
%	Denote the thickness of $\bar p$ by $\alpha$ and let $X,Y$ be the irreducible comopnents o
	\begin{enumerate}[(a)]
		\item If $0<\alpha\leq 4$, there is unique component $\hat T$ of type (b) above above $\bar p$, which meets the proper transforms of $X$ and $Y$ in double points of thickness $\alpha/2$. Moreover, we have $\wbar(\hat T)=\alpha/2$. 
		\begin{enumerate}[(i)]
			\item If $\alpha<4$, the component $\hat T$ intersects a component of type (d) in a double point of thickness $\frac{4-\alpha}{6}$.  These are all components over  $\bar p$. 
		\item  If $\alpha=4$, the preimage of $\hat T$ is irreducible and has genus $1$. There are no other components over $\bar p$. 
				\end{enumerate}
				\item If $\alpha>4$, there are two components $\hat T_{1}, \hat T_{2}$ of type (b) above $\bar p$. The component $\hat T_{1}$ intersects the proper transform of $X$ in a double point of thickness $2$. Similarly,  the component $\hat T_{2}$ intersects the proper transform of $X$ in a double point of thickness $2$. Furthermore, the components  $\hat T_{1}$ and $ \hat T_{2}$  intersect in a double point $\bar q$ of thickness $4-\alpha$. The preimage of $\bar q$ in $C_0$ consists of two double points. There are no further components over $\bar p$. 
	\end{enumerate}
\end{Prop}
\begin{Proof}
By assumption, we have $\wbar_{\bar p, X}(0)=\wbar(X)=0$ and $\wbar_{\bar p, X}(\alpha)=\wbar(Y)=0$. As $\wbar_{\bar p, X}$ is a piecewise linear function and only has horizontal segments of value $2$ by \cite[Prop. 3.4.4.]{GP24}, its first slope $s_1$ is positive and is last slope is negative. Hence by Theorem \ref{RelativeGenusOverDoublePoint}, we have $s_1=1$ and $s_2=-1$. As $\wbar_{\bar p, X}$ is concave, this implies that for all  $\lambda \in \BQ \cap [0,\alpha]$, we have
\UseTheoremCounterForNextEquation
\begin{equation}\ \wbar_{\bar p, X}(\lambda)\ = \ \min\{2,\alpha-\lambda, \lambda \}. \label{WbarIdentityEquationProofGroundedDP}\end{equation}
%for $\lambda \in \BQ \cap [0,\alpha]$. 
 First, assume $\alpha\leq 4$. Then the function $\wbar_{\bar p, X}$ has one break point of value $\alpha/2$ at $\lambda=\alpha/2$. Hence by \cite[Prop. 4.5.1]{GP24}, there is a unique component $\hat T$ of type (b) above $\bar p$ with the desired square defect $\wbar(\hat T)=\wbar(\alpha/2)=\alpha/2$. Moreover, it intersects the proper transforms of $X$ and $Y$ in double points of thickness $\alpha/2$. 
 
 If $\alpha<4$, the component $\hat T$ meets at least one component of type (c) or (d) by \cite[Prop. 4.5.8]{GP24}. As $g_{\bar p}=1$ and any component of type (d) contributes at least genus $1$, there is precisely one component of type (d) intersecting $\hat T$. The thickness of the corresponding double point is $\frac{2-\alpha/2}{3}=\frac{4-\alpha}{6}$ by Proposition \ref{ComputeThickness}. 
 
 If $\alpha=4$, the component  $\hat T$ meets no components of type (c) or (d) by \cite[Prop. 4.5.8]{GP24}. Furthermore, both double points of $\hat T$ have a unique double point above them by Corollary \ref{UniqueDoublePointOverDoublePoint}. As $C_0$ is stable, the preimage of $\hat T$ has genus $1$. 
 
 Second, assume $\alpha>4$. Then the function $\wbar_{\bar p, X}$ has two break points, one at  $\lambda_1:=2-\alpha$ and one at $\lambda_2:=\alpha-2$. Both break points have value $2$. Thus by 
 \cite[Prop. 4.5.1]{GP24}, there exist  two components $\hat T_{1}, \hat T_{2}$ of type (b) above $\bar p$ such that 
 \begin{itemize}
 	\item  $\hat T_{1}$ intersects the proper transform of $X$ in a double point of thickness $2$;
 	\item  $\hat T_{2}$ intersects the proper transform of $Y$ in a double point of thickness $2$;
 	\item $\hat T_{1}$ and $ \hat T_{2}$  intersect in a double point $\bar q$ of thickness $4-\alpha$. 
 \end{itemize}
By Propositions \ref{WbarExtraComponentProp} and \ref{WbarinCohatProp} and Equation  \eqref{WbarIdentityEquationProofGroundedDP}, the function $\wbar_{\bar q, \hat T_{1}}$ is constant with value $2$, hence the preimage of $\bar q$ in $C_0$ consists of two double points by \cite[Prop. 4.5.10]{GP24}. 
\end{Proof}

\begin{Rem}{\rm
	If $C$ has exactly $g$ grounded double points, Proposition \ref{GroundedDPAtlestgenus1} together with the total genus being $g$  imply that the local genus at any such point is $1$, hence we can use Proposition \ref{Genus1BehaviouroverGroundedDP} repeatedly to fully determine the stable reduction. }
\end{Rem}

\begin{Cor}
	Assume that $C$ has $g$ grounded double points. Then the reduction of the Jacobian of $C$ has toric rank equal to the number of grounded double points with thickness strictly larger than $4$. 
\end{Cor}
\begin{Proof}
Any irreducible component $X\subset \bar C_0$ whose normalization in $K(C)$ is split has $\wbar(X)=2$, and for each double point $\bar p\in X$ we have $g_{\bar p}>1$. This double point can't be a grounded double point, since otherwise $\wbar(X)=0$.
As the total genus of $C_0$ is $g$ and we already have $g$ grounded double points with local genus at least $1$, there exists no component $X\subset \bar C_0$ whose normalization in $K(C)$ is split. Applying Proposition \ref{Genus1BehaviouroverGroundedDP} for every grounded double point and using \cite[Prop. 4.6.7]{GP24} yields the claim. 
\end{Proof}

\begin{figure}[h]
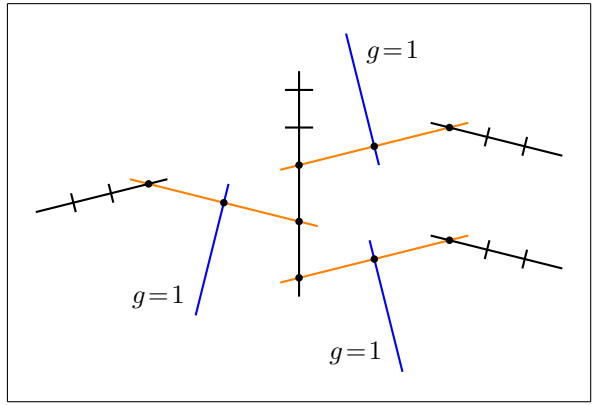
 \centering \FigGenusThreeExampleWithThreeGroundedDoublePoints 
	\caption{The stable marked reduction $ C_0$  in Example \ref{ExampleWithThreeGroundedDoublePoints} \label{FigGenusThreeWithThreeGroundedDP}}
\end{figure}

\begin{Ex}\label{ExampleWithThreeGroundedDoublePoints}{\rm
	Consider the curve $C$ defined by $z^2=F$ for
	$$F\ :=\ -7 x - 12 x^2 - 6 x^3 + 6 x^4 + 12 x^5 + 6 x^6 + x^7 \ = \ x(x^3-1)((x+2)^3-1). $$ As the degree of  $F$ is 7, we have $g(C)=3$. Moreover, since the zeros of  $F$ modulo $\Fm$ reduce to 4 distinct points,  the special fiber $C_0$ has one irreducible component $X$ with coordinate $x$. The points defined by $x=0, \infty$ are marked. Moreover, the component $X$ intersects three components $X_1, X_2, X_3$, each in a double point of thickness $1$. The components $X_1, X_2, X_3$ contain two marked points each, hence $\bar C_0$ contains 3 grounded double points.  Since $g(C)=3$, we get that the local genus over each double point is $1$. Applying Proposition \ref{Genus1BehaviouroverGroundedDP} For each double point yields that over each double point, there is a unique component of type (b) and a unique component of type (c). The stable reduction $C_0$ is depicted in Figure \ref{FigGenusThreeWithThreeGroundedDP}. Note that all double points of $C_0$ have thickness $1/4$. }
\end{Ex}

\subsection{Leaf components} \label{SubSectionLeafCompoenents}

Let $X$ be an irreducible component of $\bar C_0$ which corresponds to a leaf of its dual graph, i.e. $X$ contains exactly one double point $\bar p$ and at least two marked points. 
\begin{Prop}\label{EvenLeavGenusProp} Let $N>0$ be an integer and assume that $X$ contains exactly $2N$ marked points. Then 
	$$\sum_{\bar q \in X\setminus\{\bar p\}} g_{\bar q}\leq N-1. $$
\end{Prop}
\begin{Proof}
	Let $x$ be a coordinate of $X$ such that $\bar p$ is defined by $x=0$ and such that one marked point is defined by $x=\infty$. As in the proof of \cite[Prop. 4.3.1]{GP24}, there exist $F_1\in R[x]$ of degree $2N-1$ and with constant coefficient $1$ and $F_2\in 1+\Fm[x^{-1}] $ such that for $F:=F_1 F_2$, the equation $z^2=F$ defines $C$. Since $X$ contains a marked point, we have $\wbar(X)=0$ by Proposition \ref{SquaredefectZeroProp}, hence $F=0^2+F$ is an optimal decomposition. As $[F_2]=[1]$, we compute  
	$$[\tfrac{d}{dx} F]\ = \ [\tfrac{d}{dx} F_1]$$ 
	Here $[\tfrac{d}{dx} F_1]$ is a polynomial in $x^2$ of degree $N-1$, hence the sum of the local genera $ g_{\bar p}$ at closed smooth points is at most $N-1$ by Theorem \ref{LocalGenusUnmarkedPointProp}. Furthermore, as $\wbar(X)=0$, the local genus over the generic point of $X$ is $0$ as well by \cite[Prop. 3.3.2]{GP24}.
\end{Proof}

\begin{figure}[h]
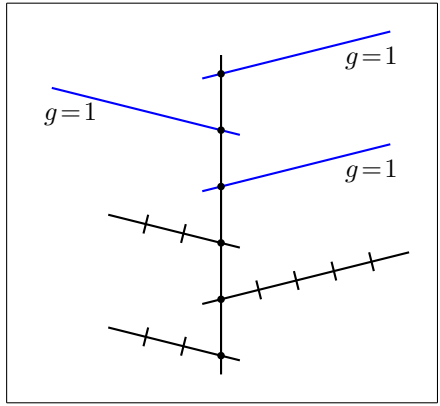
 \centering \FigGenusThreeExampleNotingOverDP 
		\caption{The stable marked reduction $ C_0$  in Example \ref{ExampleEvenLeaves} (b) \label{FigGenusThreeNothingOverDP}}
	\end{figure}

The inequality in Proposition \ref{EvenLeavGenusProp} is sometimes strict. Furthermore, there exist curves for which $g_{\bar p}+\sum_{\bar q \in X\setminus\{\bar p\}} g_{\bar q} >N-1$. We provide examples below:

\begin{Ex}\label{ExampleEvenLeaves}{\rm 
	\begin{enumerate}[(a)]
		\item Let $C$ be a curve of genus $2$ whose reduction type is of type (B4). Then $\bar C_0$ has exactly one double point $\bar p$ and two components $X_1, X_2$, one marked by $4$ and the other one marked by $2$ points. We have $g_{\bar p}=2$ and 
		$$\sum_{\bar q \in X_1\setminus\{\bar p\}} g_{\bar q} \ = \  \sum_{\bar q \in X_2\setminus\{\bar p\}} g_{\bar q} =0$$		
		\item 	Consider the curve $C$ defined by $z^2=F$ for 
		$$ F:=
%		2 x^{7}+7 x^{6}+2 x^{5}+9 x^{4}+50 x^{3}+40 x^{2}-10 x-12 \ = \ 
		(x^4 + 6x + 4)(x + 1)(x + 3)(2x - 1).$$ 
		
		As the degree of  $F$ is 7, we have $g(C)=3$. Moreover, since the zeros of  $F$ modulo $\Fm$ reduce to 2 distinct points and the degree of $F$ is odd,  the special fiber $\bar C_0$ has one irreducible component $X$ with coordinate $x$. The component $X$ intersects three components $X_1, X_2, X_3$. The components $X_1, X_2$ contain two marked points each, and the component $X_3$ contains $4$ marked points.  Let $c\in R$ be a square root of $-1$ and let $H:=3x^2+cx^3.$ Then
		$$G:=F-H^2=12+2 x^{7}+8 x^{6}+\left(2-6 \, c\right) x^{5}+50 x^{3}+40 x^{2}-10 x$$
		is odd modulo $2\Fm[x]$, hence $F=H^2+G$ is optimal and $\wbar(F)=1$ by Lemma \ref{FirstApproximationLemma}. We compute   
		$$[\frac{dG}{dx}/2]=[x^6 + x^2 + 1].$$
		Hence there are three smooth points of $X$ at which the local genus of $C_0$ is $1$ by Theorem \ref{LocalGenusUnmarkedPointProp}. This implies that over these points there are components of type (d), whose preimages in $C_0$ have genus $1$. The stable reduction $C_0$ is depicted in Figure \ref{FigGenusThreeNothingOverDP} and for each $1\leq i \leq 3$, we have
		$$\sum_{\bar q \in X_1} g_{\bar q} =0.$$
	\end{enumerate}}
	\end{Ex}

\begin{Prop}\label{OddLeavGenusProp}
	Let $N$ be an integer and assume that $X$ contains exactly  $2N+1$ of marked points. Then
	$$\sum_{\bar q \in X} g_{\bar q}=N. $$
\end{Prop}
\begin{Proof}
	Let $x$ be a coordinate of $X$ such that $\bar p$ is defined by $x=0$ and  no marked point is defined by $x=\infty$. By \cite[Prop. 4.3.1]{GP24}, there exists a separable $F_1\in R[x]$ with constant coefficient $1$ and $F_2\in 1+\Fm[x^{-1}] $ such that for $F:=xF_1 F_2$, the equation $z^2=F$ defines $C$. Moreover, as no marked point is defined by $x=\infty$, all $2N+1$ marked point are defined by roots of $xF_1$, hence $\deg(xF_1)=2N+1$ and the leading coefficient of $F_1$ is a unit. Since $X$ contains a marked point, we have $\wbar(X)=0$ by Proposition \ref{SquaredefectZeroProp}, hence $F=0^2+F$ is an optimal decomposition. As $[F_2]=[1]$, we compute  
	$$[\tfrac{d}{dx} F]\ = \ [\tfrac{d}{dx} (x F_1)].  $$ 
	Since $xF_1$ is a polynomial of odd degree whose leading coefficient is a unit in $R$, the residue class of its derivative $P:=[\frac{d}{dx} xF_1]$ is a polynomial in degree ${2N}$. Counted with multiplicity, the polynomial $P$ has $2N$ roots. Furthermore, the constant coefficient of $P$ is $1$, hence $P(0)\neq 0$. 
	By Theorem \ref{LocalGenusUnmarkedPointProp}, the sum of the local genera at points of $X$ is equal to half the sum of the multiplicities of the roots of $P$, so 
	$$\sum_{\bar q \in X} g_{\bar q}=N, $$ as desired. 
\end{Proof}

\begin{figure}[h]
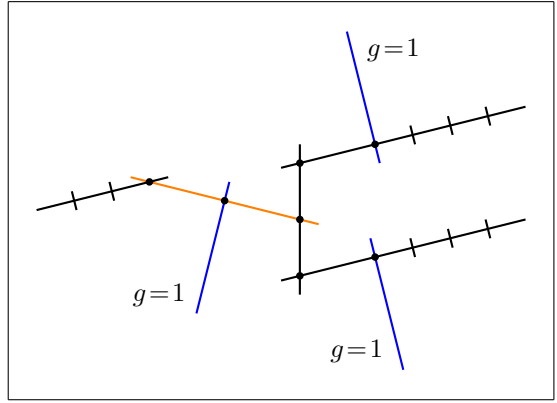
 \centering \FigGenusThreeExampleOneGroundedDPdownstairs
	\caption{The stable marked reduction $C_0$  in Example \ref{GenusThreeExampleOneGroundedDP} \label{GenusThreeExampleOneGroundedDPFigure} }
\end{figure}

\begin{Ex}\label{GenusThreeExampleOneGroundedDP}{\rm
	%\begin{enumerate}[(a)]
	%	\item Define $F:=x^7+2x^6-x^2-2x$ and let $C$ be the curve defined by the equation $z^2=F$.  The polynomial  $F$ factors as $$F=x(x^5+1)(x+2), $$ hence $\bar C_0$ has 
	%	\item 
	Define $F:=12x^7 - 18x^6 - 11x^5 + 3x^4 - 11x^3 + 3x^2 - 2x$ and let $C$ be the curve defined by the equation $z^2=F$.  The polynomial  $F$ factors as 
	$$F\ =\ (x + 1)(3x^2 + 1)(4x^2 - 2x + 1)(-2 + x)x. $$
	From the factorization of  $F$, we read off that  $\bar C_0$ has an unmarked component $X$ with coordinate $x$. Direct computations show that $\bar C_0$ has a component $X_1$, which is marked by two points and intersects $X$ in an even double point $\bar q$ of thickness 1 defined by $x=0$.  Moreover, there are further components $X_2, X_3$ of $\hat C_0$, which contain three marked points each and intersect $X$ in odd double points of thickness 1 at $x=1, x=\infty$. 
	
	By Proposition \ref{OddLeavGenusProp}, the sum of the local genera above points of $X_2$ is $\frac{3-1}{2}=1$. Hence there exists a component of type (d) above $X_2$ whose preimage in $C_0$ is an elliptic curve. By the same argument, there exists a component of type (d) above $X_2$ whose preimage in $C_0$ is an elliptic curve.
	
	Moreover, as the component $X_1$ contains a marked point and $X$ contains an odd double point, Proposition \ref{SquaredefectZeroProp} yields $\wbar(X_1)=\wbar(X)=0$. Hence $\bar q$ is a grounded double point and we have $g_{\bar q}\geq 1$ by Proposition \ref{GroundedDPAtlestgenus1}. Since the components we already found altogether have genus contribution $2$ and the total genus is $3$, we get $g_{\bar q}=1$. By Proposition \ref{Genus1BehaviouroverGroundedDP}, we get that there is a unique component of type (b) above $\bar q$, which intersects a component of type (d)  whose preimage in $C_0$ is an elliptic curve. This fully determines the stable reduction $C_0$, see  Figure \ref{GenusThreeExampleOneGroundedDPFigure}.
%	See figure \ref{GenusThreeExampleOneGroundedDPFigure} For a drawing of the $\bar C_0$. 
	%	
	%	
	%	$\bar C_0$ has components $X_1, X_\infty$, which meet $X$ in odd double points of thickness  $1$ and are each marked by three points. Moreover, there is another component $X_0$ which intersects $X$ in an even double point $\bar p$ of thickness 1 defined by $x=0$.  
%	As $X_1$ contains a marked point and $X$ contains an odd double point, we have $\wbar(X_1)=\wbar(X)=0$ by Proposition \ref{SquaredefectZeroProp}. Thus $\bar p$ is a grounded double point.
%	 Let $F/x^2=H^2+G$ be an odd decomposition. As $\wbar(X)=0$, we have  $v(G))=0$ as well. Moreover as the linear term of $H^2$ has positive valuation and
%	$$[f/x^2]=[x^5+x^3+x^2+x+1],$$
%	we have $\op{ldeg}(G)=1$. Hence the first slope of $\wbar_{\bar p, X}$ is $1$ by Lemma \ref{slopeLemRightSide}. Similarly, we get that the last slope of $\wbar_{\bar p, X}$ is $-1$. By Theorem \ref{RelativeGenusOverDoublePoint}, this implies $g_{\bar p}=1$. Hence $\bar p$ satisfies the conditions of Theorem \ref{RelativeGenusOverDoublePoint} and we get that there are two components of $\hat C_0$ above $\bar p$; one of type (b) and one of type (d).  
	%\end{enumerate}	
}
\end{Ex}
%%%ToDo: Maybe this is unneccesarry and I only need the stable reduction?

\subsection{Local genus and odd double points}\label{SubSectionConjectureOddDoublePoints}

In this subsection, assume that $\bar p$ is any odd double point of $\bar C_0$.  We conjecture the following generalization of Proposition \ref{OddLeavGenusProp} holds true:
\begin{Conj}\label{GenusOddDoublePointConjecture}
	Let $\bar p \in \bar C_0$ be an odd double point and let $T$ be one connected component of $\bar C_0\setminus\{\bar p\}$ containing $2N+1$ marked points. Then 
	$$\sum_{\bar q \in T} g_{\bar q}=N. $$
\end{Conj}
For genus $2$, Conjecture \ref{GenusOddDoublePointConjecture} follows directly from the classification of the stable marked reduction given in \cite[Section 5.2]{GP24}. Straightforward arguments using Propositions \ref{GroundedDPAtlestgenus1}, \ref{OddLeavGenusProp} and \ref{EvenLeavGenusProp} show that it holds true for $g=3$ too. In the example below, we demonstrate such an argument under the assumption that the dual graph of $\bar C_0$ takes a specific form. 

\begin{figure}[h]
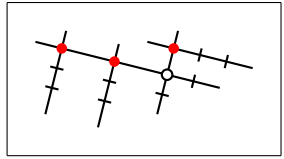
 \centering \FigMotivatingExampleFigure
	\caption{The stable marked reduction $\bar C_0$  in Example \ref{MotivatingExamplesForOddDPConjecture} \label{MotivatingExampleFigure} }
\end{figure}

\begin{Ex}\label{MotivatingExamplesForOddDPConjecture}
{\rm	Let $C$ be a curve of genus $3$ such that 
	\begin{itemize}
		\item $\bar C_0$ has $5$ components $X_1, \dots, X_5$;
		\item $X_1, X_2, X_5$ are marked by two points each;
		\item $X_3, X_4$ are marked by one point each;
		\item $X_3$ intersects $X_1, X_2, X_4$;
		\item $X_4$ intersects $X_5$. 
	\end{itemize}
	Then the double point $\bar p$ in which $X_3$ and $X_4$ intersect is odd. Moreover, the other $3$ double points are all grounded double points, hence by Proposition \ref{GroundedDPAtlestgenus1} the local genus above them must be $1$. 	See Figure \ref{MotivatingExampleFigure}, where red circles signify grounded double points and the empty circle signifies the odd double point.  }
\end{Ex}

\section{Examples}\label{SectionExamples}
%\subsection{Curves with many grounded double points}
\subsection{Curves of genus 3}\label{SubSectionGenusThreeCurves}
In Section 5.2 of \cite{GP24}, Richard Pink and the author of this article provide a classification of the stable marked reduction for genus $2$ curves into 54 cases. This classification depends only on the dual graph of $\hat C_0$ together with thicknesses of its double points and an additional parameter $\delta\in \BQ\cap [0,2]$. Moreover, for some dual graphs, the parameter $\delta$ was not necessary. In the light of the new methods introduced in this article,  these $15$ cases have in common that the dual graph of $C_0$ can be deduced using Propositions \ref{GroundedDPAtlestgenus1} and \ref{Genus1BehaviouroverGroundedDP} regarding grounded double points and Proposition \ref{OddLeavGenusProp} computing the local genus above points in odd leafs. 
%In the language of cluster pictures \cite{ToDo}, this data is equivalent to the cluster picture of $C$ together with $\delta$. 

In this section, we employ this method to determine some cases of the stable marked reduction of hyperelliptic curves of genus $3$. For the remainder of this section, assume that $g=3$. Here $C \onto \bar C$ has $8$ branch points. In total, there are 32 possibilities of the dual graph of $(\bar C_0, \bar p_1, \dots, \bar p_8)$. In \cite{WorksheetsG24}, we use the admcycles package \cite{Admcycles} to produce a list of these.

\begin{figure}[h]
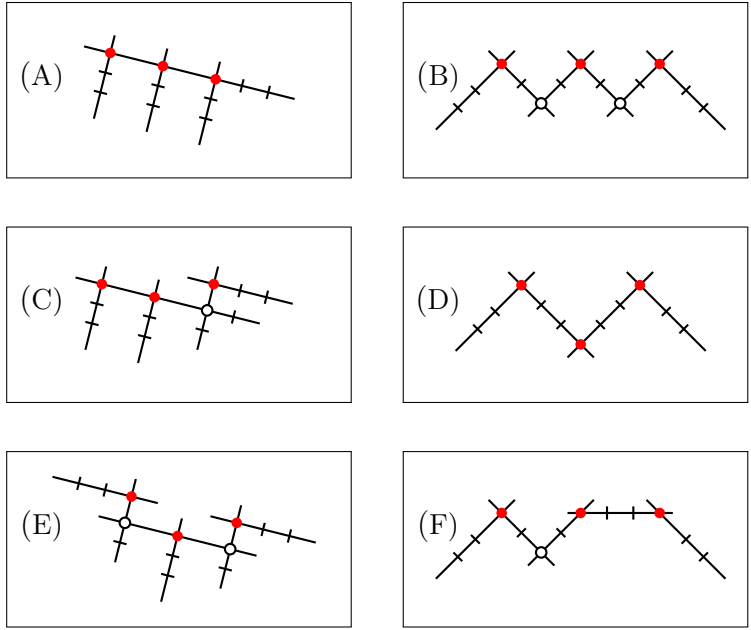
 \centering 
	\FigSomeThree
	\caption{Some of the possibilities for $(\bar C_0,\bar p_1,\ldots,\bar p_8)$ in genus~$3$.}\label{FigThreeSome}
\end{figure}

Going through the list of possibilities for the combinatorial structure of $(\bar C_0, \bar p_1, \dots, \bar p_8)$, we find exactly six cases in which the curve has three grounded double points are shown. In Figure \ref{FigThreeSome} these six cases are drawn. Filled  red circles signify grounded double points, filled black circles non-grounded even double points and empty circles odd double points. Over each grounded double point, the local genus of $C$ is $1$ and hence there are three possibilities for the preimage of each of them in $C_0$ by Proposition \ref{Genus1BehaviouroverGroundedDP} as we can choose the thicknesses of the double points freely and independently from each other. 

{\bf Case (A):} Here the combinatorial structure of $\bar C_0$ is invariant under exchanging all three grounded double points, hence this yields $10$ cases for the type of the stable marked reduction of $C$, which only depend on the thicknesses of the grounded double points. 

\medskip 

{\bf Cases (B)-(E):} Here the combinatorial structure of $\bar C_0$ is only invariant under exchanging two of the grounded double points. In each case, this yields $18$ cases for the type of the stable marked reduction of $C$, which only depend on the thicknesses of the grounded double points. 

\medskip 

{\bf Case (F):} Here there there is no symmetry in the combinatorial structure of $\bar C_0$ that exchanges grounded double points. Accordingly, there are $27$ cases for the  type of the stable marked reduction of $C$, which only depend on the thicknesses of the grounded double points. 

\medskip 

\begin{figure}[h]
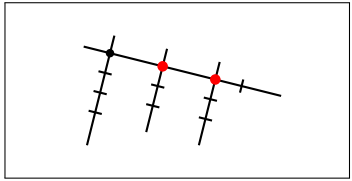
 \centering 
	\FigGenusThreeOnlyTwoGroundedDP
	\caption{One case of $(\bar C_0,\bar p_1,\ldots,\bar p_8)$ with two grounded double points.}\label{FigGThreeOnlyTwoGroundedDP}
\end{figure}

Using Proposition \ref{OddLeavGenusProp}, one may similarly treat cases where $\bar C_0$ has less than $3$ grounded double points. As an example, consider the case displayed in Figure \ref{FigGThreeOnlyTwoGroundedDP}, where we have two grounded double points. Moreover,  by Proposition \ref{OddLeavGenusProp} there is a component $\hat T$ of type (d) over the leaf component with three marked double points. As the total genus is $3$ and each grounded double point must contribute genus at least one, the preimage of $\hat T$ in $C_0$ is elliptic. Here the combinatorial structure of $\bar C_0$ is invariant under exchanging the grounded double points, hence this yields $6$ cases for the  type of the stable marked reduction of $C$, which again only depend on the thicknesses of the grounded double points.

In total the mentioned cases amount to $10+4\cdot 18+27+6=115$ cases for the combinatorial structure of $C_0$. This is more than double the total number of cases (54) for genus $2$ curves. While our current methods would in theory allow to give a full classification for higher genus, we refrain from this task at this time due to its high combinatorial complexity.

\subsection{Curves with many automorphisms}\label{CurvesWithManyAutomorphims}
In \cite{MuellerPink}, Nicolas Müller  and Richard Pink classify all hyperelliptic curves over $\mathbb{C}$ with many automorphisms. In particular, they determine which of these curves do or do not have complex multiplication. For five of the curves $X_{10},X_{11}, X_{16}, X_{17}, X_{18}$, a special argument was needed to show that the Jacobian of these curves does not have complex multiplication. If the toric rank of the reduction of one of the Jacobians was $\neq 0$, then it would immediately follow that it does not have complex multiplications. In \cite{MuellerPink}, it was proved that the toric rank of the reduction of the Jacobians is $0$ for all places of odd residue characteristic. 

One motivation for this article was to compute the stable marked reduction of the curves $X_{10},X_{11}, X_{16}, X_{17}, X_{18}$ over $2$ and to read off the toric rank of the reduction of the Jacobian. In this section, we describe the reductions of these curves. It turns out that the toric rank of the reduction of the Jacobian is $0$ in all cases, which means that we cannot deduce that the Jacobians do not have complex multiplication. To our knowledge $X_{18}$ is the curve of highest genus for which the stable reduction over $2$ has ever been explicitly computed. 

\begin{figure}[h]
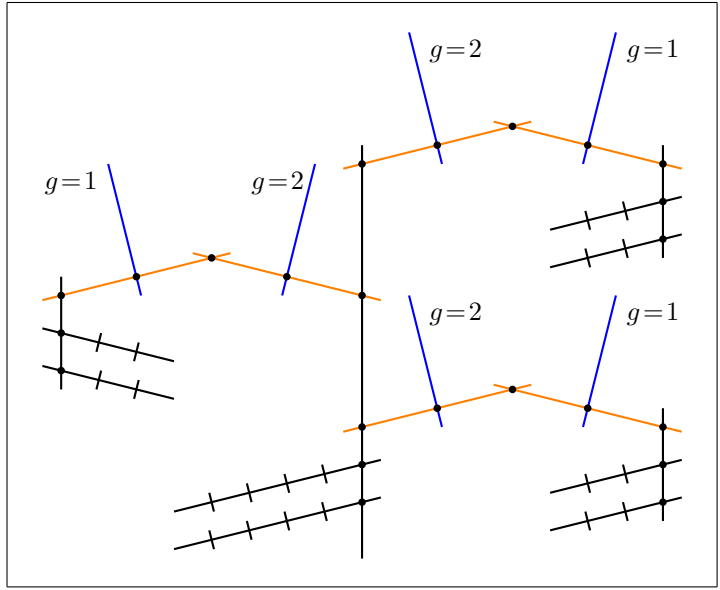
 \centering \FigXten
	\caption{Stable reduction of $X_{10}$, which has genus $9$.}\label{FigX10}
\end{figure}

We only list the final results, leaving the detailed
computations to look up in the associated computer algebra worksheets \cite{WorksheetsG24}. The equation of the curves  are given in terms of certain separable polynomials from Table \ref{TablePolynomials}.
\renewcommand{\arraystretch}{1.2}
\begin{table}[h]
	\centering
	\begin{tabular}{|c||c|}
		\hline
		$t_4$ & $x(x^4 - 1)$ \\
		\hline
		$r_4$ & $x^{12} - 33x^8 - 33x^4 + 1$ \\
		\hline
		$s_4 $ & $x^8 + 14x^4 + 1$ \\
		\hline
		$r_5$ & $x^{20} - 228x^{15} + 494x^{10} + 228x^5 + 1$ \\
		\hline
		$s_5$ & $x(x^{10} + 11x^5 - 1)$ \\
		\hline
		$t_5$ & $x^{30} + 522x^{25} - 10005x^{20} - 10005x^{10} - 522x^5 + 1$ \\
		\hline
	\end{tabular}
	\caption{Certain separable polynomials over $\mathbb{C}$}\label{TablePolynomials}
\end{table}
\renewcommand{\arraystretch}{1}
\newpage 
\begin{Ex} Consider the curve $X_{10}$ defined by
%	$$z^2= \left(x^{12}-33 x^{8}-33 x^{4}+1\right) \left(x^{8}+14 x^{4}+1\right).$$
		$$z^2= r_4 s_4.$$
	 This is a curve of genus $9$ and its stable marked reduction is sketched in Figure \ref{FigX10}.
\end{Ex}

%%%%%Comment these in when needed, slows down the document a lot. 

\begin{figure}[h]
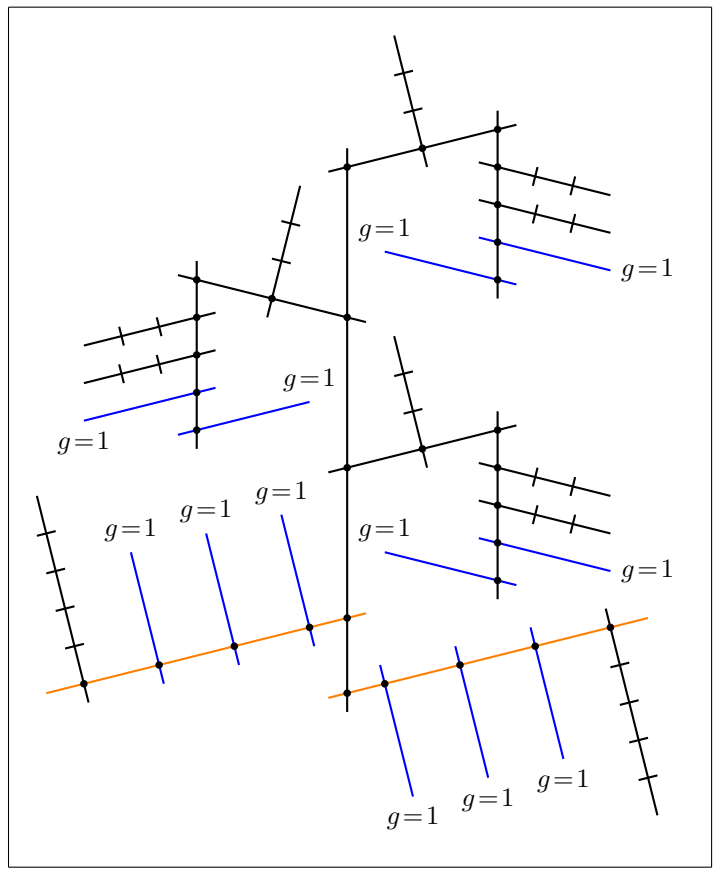
 \centering \FigXeleven
	\caption{Stable reduction of $X_{11}$, which has genus $12$.}\label{FigX11}
\end{figure}

\newpage \texttt{\texttt{\texttt{}}}
\begin{Ex} Consider the curve $X_{11}$ defined by
%	$$z^2=\left(x^{12}-33 x^{8}-33 x^{4}+1\right) \left(x^{8}+14 x^{4}+1\right) x \left(x^{4}-1\right).$$ 
	$$z^2=r_4s_4t_4.$$ 
	This is a curve of genus $12$   and its stable marked reduction is sketched in Figure \ref{FigX11}.
\end{Ex}

\begin{figure}[h]
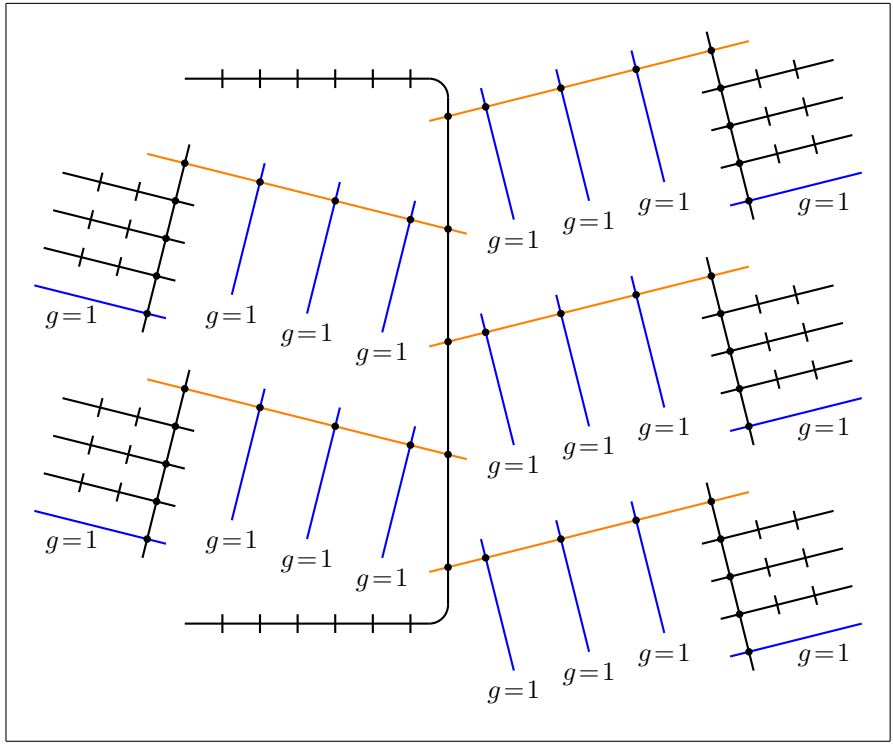
 \centering \FigXsixteen
	\caption{Stable reduction of $X_{16}$, which has genus $20$.}\label{FigX16}
\end{figure}

\newpage 
\begin{Ex} Consider the curve $X_{16}$ defined by
%	$$z^2=x \left(x^{10}+11 x^{5}-1\right) \left(x^{30}+522 x^{25}-10005 x^{20}-10005 x^{10}-522 x^{5}+1\right).$$ 
$$z^2=s_5t_5.$$ 
	This is a curve of genus $20$ and its stable marked reduction is sketched in Figure \ref{FigX16}.
\end{Ex}

\begin{figure}[h]
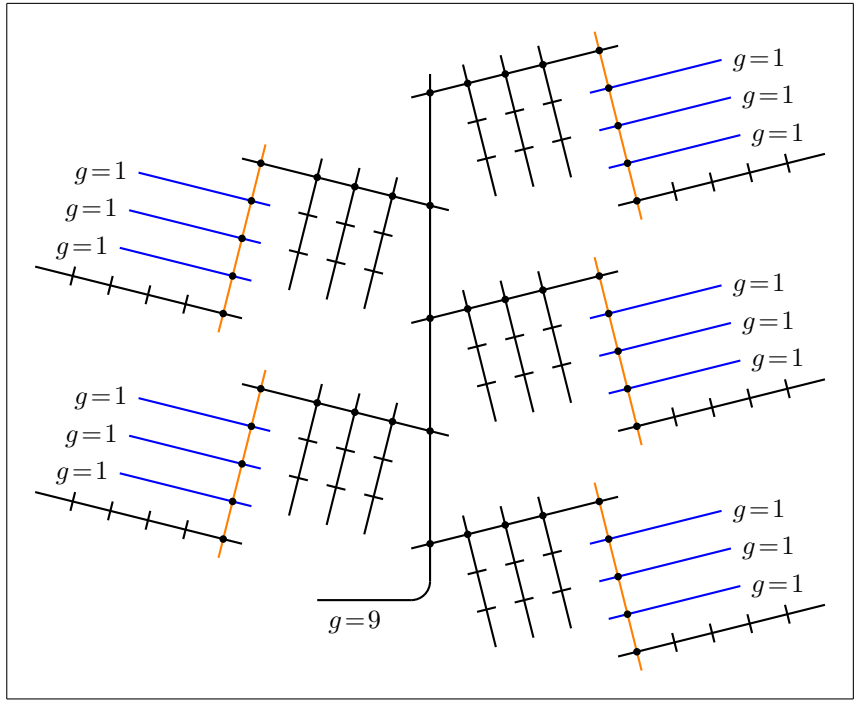
 \centering \FigXseventeen
	\caption{Stable reduction of $X_{17}$, which has genus $24$.}\label{FigX17}
\end{figure}

\newpage 
\begin{Ex} Consider the curve $X_{17}$ defined by
%	$$z^2=\left(x^{20}-228 x^{15}+494 x^{10}+228 x^{5}+1\right) \left(x^{30}+522 x^{25}-10005 x^{20}-10005 x^{10}-522 x^{5}+1\right).$$ 
		$$z^2=r_5t_5.$$ 
	This is a curve of genus $24$  and its stable marked reduction is sketched in Figure \ref{FigX17}.
\end{Ex}

\begin{figure}[h]
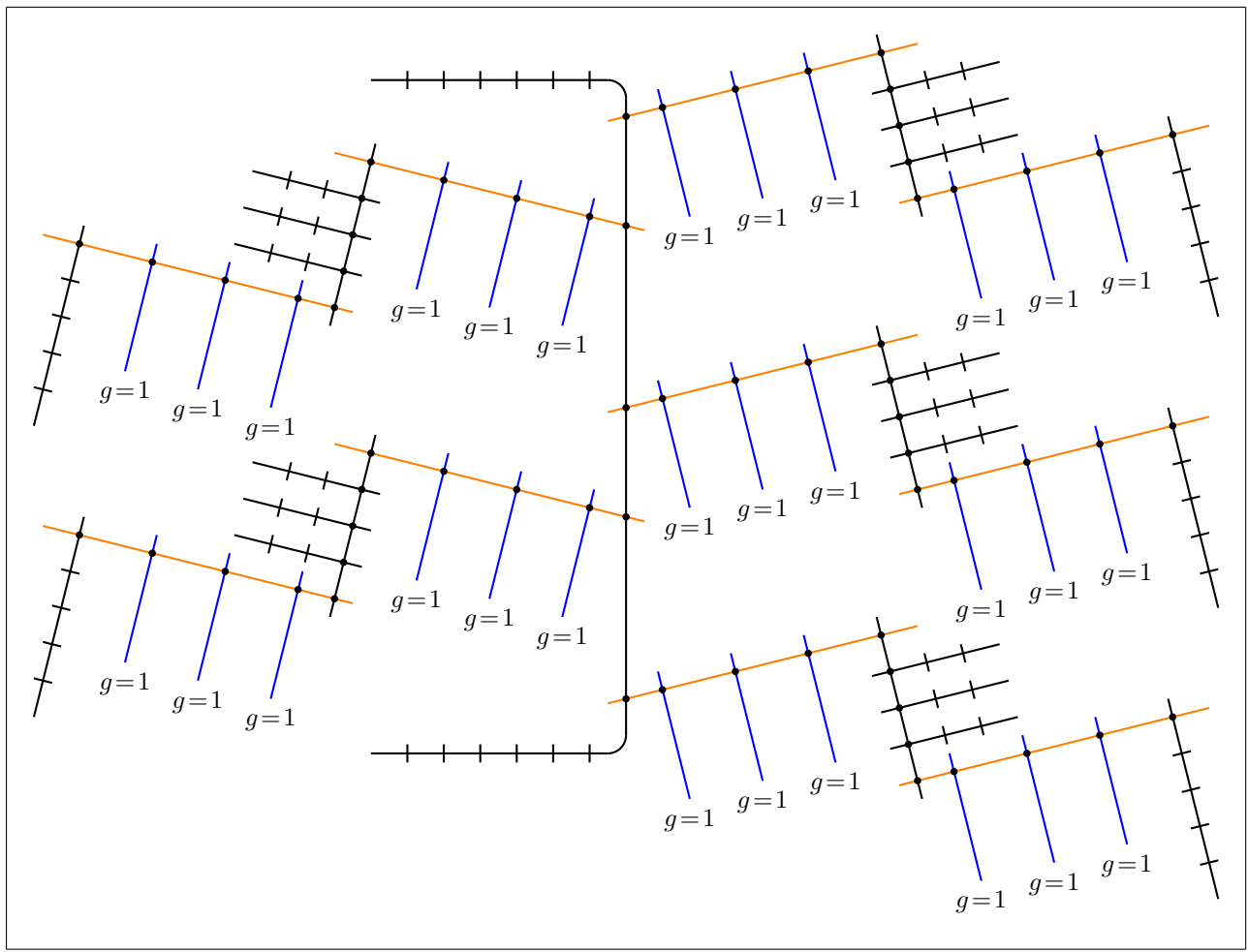
 \centering \FigXeightteen
	\caption{Stable reduction of $X_{18}$,  which has genus $30$.}\label{FigX18}
\end{figure}

\newpage 
\begin{Ex} Consider the curve $X_{18}$ defined by
%	$$z^2=\left(x^{20}-228 x^{15}+494 x^{10}+228 x^{5}+1\right) \left(x^{30}+522 x^{25}-10005 x^{20}-10005 x^{10}-522 x^{5}+1\right) x \left(x^{10}+11 x^{5}-1\right)).$$ 
	$$z^2=r_5s_5t_5.$$ 
	This is a curve of genus $30$  and its stable marked reduction is sketched in Figure \ref{FigX18}.
\end{Ex}

 \printbibliography
%%%%%%%%%%%%%%%%%%%%%%%%%%%%%%%%%%%%%%%%%%
\end{document}